\documentclass[10pt]{asl}

\usepackage{bbm,color}
\usepackage{xcolor}
\usepackage{times,enumitem,relsize}

\definecolor{blue}{rgb}{0,0,1}

\usepackage[textsize=footnotesize,color=blue!40, bordercolor=white]{todonotes}

\usepackage{verbatim}
\usepackage{lscape}
\usepackage{enumitem}
\usepackage{todonotes}
\usepackage{tikz}
\usetikzlibrary{matrix}

\newtheorem{thm}{Theorem}[section]
\newtheorem{lem}[thm]{Lemma}
\newtheorem{col}[thm]{Corollary}
\newtheorem{prop}[thm]{Proposition}

\newtheorem*{claim}{Claim}
\newtheorem*{thm*}{Theorem}
\newtheorem*{col*}{Corollary}

\theoremstyle{definition}
\newtheorem{defn}[thm]{Definition}

\newtheorem*{rem}{Remark}
\theoremstyle{remark}

\newcommand{\alp}{\alpha}
\newcommand{\bt}{\beta}
\newcommand{\gm}{\gamma}
\newcommand{\dlt}{\delta}

\newcommand{\kp}{\kappa}
\newcommand{\lmb}{\lambda}

\newcommand{\om}{\omega}
\renewcommand{\phi}{\varphi}

\newcommand{\Gm}{\Gamma}
\newcommand{\A}{\mathcal A}

\newcommand{\Cf}{\mathcal C}
\newcommand{\F}{\mathcal F}
\newcommand{\La}{\mathcal L}
\newcommand{\N}{\mathcal N}
\newcommand{\M}{\mathcal M}

\renewcommand{\P}{\mathcal P}

\newcommand{\trace}{\mathcal s}

\newcommand{\restr}{\mathord{\restriction}}
\newcommand{\lrar}{\leftrightarrow}
\newcommand{\la}{\langle}
\newcommand{\ra}{\rangle}

\newcommand{\emp}{\emptyset}
\newcommand{\sbs}{\subseteq}
\newcommand{\sps}{\supseteq}
\newcommand{\without}{\mathord{\setminus}}
\newcommand{\diag}{\bigtriangleup}
\newcommand{\diam}{\diamondsuit}

\newcommand{\qedtwo}[1]{\hfill QED (#1)}
\def\qed{\mbox{}\hfill{QED}\\}
\def\cof{\mathop{cf}}
\newcommand{\cred}[1]{{\color{red}{#1}}}
\newcommand{\ie}{{\itshape{i.e.}}{\hspace{0.25em}}}
\newcommand{\eg}{{\itshape{e.g.}}{\hspace{0.25em}}}
\newcommand{\etc}{{\itshape{etc.}}{\hspace{0.25em}}}
\newcommand{\via}{{\itshape{via }}{\hspace{0.25em}}}

\def\Bon{\mbox{{\raisebox{-0.08em}{$\mathlarger{\mathlarger{\Box}}$}}}}
\def\nod{\noindent}

\begin{document}
\title{ Generalisations of Stationarity, Closed and Unboundedness, and of Jensen's $\Bon$ }

\author{{H. Brickhill 
and P. D. Welch}}
\email{hazel.brickhill@bristol.ac.uk, \, p.welch@bristol.ac.uk}

\address{School of Mathematics,\\ University of Bristol,\\ Bristol,\\ BS8 1TW,\\ England}

{\thanks{
Most of the results in this paper formed part of the first author's PhD thesis under the second author's supervision.
This research was made possible by funding from the EPSRC, for which the first author is very grateful. The authors would also like to warmly thank Joan Bagaria for many discussions on this topic.
}}

\maketitle

\begin{abstract}

The concepts of closed unbounded (club) and stationary sets are 
generalised  to $\gamma$-club and $\gamma$-stationary sets, which are closely related to stationary reflection.  
We use these notions to define generalisations of Jensen's combinatorial principles {$\Bon$}  and $\diamondsuit$. 


We define  $\Pi^1_\gamma$-indescribability and  use the new $\Bon^{\gamma}$-sequences to extend the result of Jensen that in the constructible universe a regular cardinal  is stationary reflecting if and only if it is $\Pi^1_1$-indescribable: we  show that in $L$ 
 a cardinal is $\Pi^{1}_{\gamma}
$-indescribable iff it reflects $\gamma$-stationary sets. But more particularly
(just stating the special case of $n$ finite):



{\bf Theorem} $(V=L)$
Let $n<\omega$  and $\kp$ be $\Pi^{1}_{n}$-indescribable but not $\Pi^1_{n+1}$-indescribable, and let $A\subseteq\kp$ be $n+1$-stationary. Then there are $E_{A}\sbs A$ and a $\Bon^{n}$-sequence $S$ on $\kp$ such that $E_{A}$ is $n+1$-stationary in $\kp$ and $S$ avoids $E_{A}$.
Thus $\kp$ is not $n+1$-reflecting.

Certain assumptions on the $\gm$-club filter  allow us to prove that  $\gm$-stationarity is downwards absolute to $L$, and allows for splitting of $\gamma$-stationary sets. 
We  define $\gm$-ineffability, and look into the relation between $\gm$-ineffability and various $\diamondsuit$ principles.
\end{abstract}


\section{Introduction}

This paper is inspired by \cite{BaMS} where the authors introduced generalisations of stationarity: essentially the idea is to iterate {\em stationary reflection}: if some $S\subseteq \kappa$ is some stationary set then $\{\alpha <\kappa \mid S\cap \alpha \mbox{ is stationary}\}$ may be also stationary.  If this happens for all stationary subsets of $\kappa$ then the latter is called {\em stationary reflecting}. Of course some such ideas have been studied for a long period (see \eg \cite{magidorrss} \cite{mshelah}).

 One {conceptual} difference here with the presentation in \cite{BaMS} is the introduction of $\gm$-club sets, which are interesting combinatorial objects in their own right and also enable us to lift some arguments straight from the club and stationary set context. As combinatorial objects, $\gm$-clubs allow us to generalise the combinatorial notion of a $\Bon$ sequence, which is  one of the main results of this paper. 

 Parts of the extant literature have close connections with $\gm$-stationary and $\gm$-club sets - for example in \cite{ben-neria}, Ben-Neria uses a ``strong stationary reflection" property, which in our terminology is just the normality of the $1$-club filter.

The notions of $\gm$-club and $\gm$-stationary set are  natural, even without the study of stationary reflection. One way to see this is to think about how we might want to measure the \emph{size} or \emph{thickness} of subsets of a particular ordinal - not, of course, in terms of cardinality but rather in a way that takes into account more of the combinatorial aspects of ordinals. 
Starting with a stationary set $S\sbs\alp$, we can ``thicken" this set by requiring that $S$ is \emph{stationary closed}, i.e. whenever $S$ is stationary below any $\bt<\alp$ then we have $\bt\in S$. We shall call such sets \emph{$1$-club}.  (We recently found that this notion was defined independently by Sun in \cite{sun}. This illustrates how natural the notion is, although it was not really picked up on at the time. Sun does not define $n$-clubs for $n>1$, but a related notion called $n$-club is given in \cite{hellsten}, and very recently explored further in \cite{cody}.)  Such sets need not be club (a set can be unbounded without being stationary), so in terms of size they come between stationary and club sets. $1$-club sets share something of the structure of the club sets as they are defined by simply replacing ``closure'' with ``stationary-closure'', and the ``unbounded" with ``stationary" in the definition, but do the $1$-club sets generate a filter? This is where we see the connection with stationary reflection - if an ordinal $\alp$ reflects any two stationary sets \emph{simultaneously} i.e. if for any stationary $A,B\sbs\alp$ there is some $\bt<\alp$ such that $A\cap\bt$ and $B\cap\bt$ are both stationary in $\bt$, then the $1$-clubs do indeed generate a filter. At such ordinals, we can then define the \emph{$2$-stationary sets} as those which intersect every $1$-club (\ie are of positive measure with respect to the $1$-club filter). 
We can continue by defining a \emph{$2$-club} set to be a $2$-stationary set which is also $2$-stationary closed, and then, where an ordinal simultaneously reflects $2$-stationary sets we can show that these $2$-clubs generate a filter. Hence, we can there define $3$-stationary sets, and so on. 

These definitions were motivated by \cite{BaMS}; our notion of $\gm$-stationarity is not quite the same as is defined there, but the two definitions coincide if we are in $L$. We shall also see that our $\gm$-stationarity is exactly the $\gm$-s-stationarity Bagaria later defined in in a more recent paper \cite{B2TAMS}. (These definitions were made independently of each other.)
He gives a rather different motivation for this notion, using \emph{derived topologies} and the definition grew out of questions in proof theory and modal logic raised in \cite{provability}. We shall not discuss these connections further here (for details see \cite{B2TAMS} and \cite{provability}), but this shows that the notions we shall define have very broad range of application. There are many other possible applications to explore with these notions. In this paper we generalise the combinatorial principles of $\Bon$  (in Section 3) and $\diam$  (in Section 4) using $\gm$-club and $\gm$-stationary sets. We also look in some detail at how these properties manifest in G\"odel's constructible universe $L$. The $\Bon$ sequences constructed in Section 3 and 4 give important results about the $\gm$-club filter in $L$, which are used to show that, under certain assumptions, $\gm$-stationarity is downward absolute to $L$. 

In more detail: 
in Section 2 we define the notions of $\gm$-closed unbounded ($\gm$-club) and $\gm$-stationary set. After giving the basic properties of these sets in section 2.1 
we consider the main result of \cite{BaMS} which was that for $n<\om$ we have that in $L$ a regular cardinal $\kp$ is $n+1$-stationary in their sense, iff $\kp$ is $\Pi^1_n$-indescribable.
We give our statement of what is essentially this theorem of \cite{BaMS}:\\

{\sc Theorem 2.14} \cite{BaMS}\label{BaMSthm}
Assume $V=L$. Then a regular cardinal $\kp$ is $n$-reflecting iff $\kp$ is $\Pi^1_n$-indescribable.\\

We proceed to define the $\gm$-club filter $\mathcal{C}^{\gm}(\kappa)$ and then for finite $\gm$ we investigate the relation between the $\gm$-club filter and $\Pi^1_\gm$-indescribability. \\

{\sc Lemma 2.20}\label{filtersV1}
For any $n<\om$ and $\Pi^1_n$-indescribable cardinal $\kp$ we have
\vspace{-0.6em}
\begin{enumerate}
\item $\Cf^n(\kp)$ is contained in the $\Pi^{1}_{n}$-indescribability filter.
\item $\Cf^n(\kp)$ is normal.
\end{enumerate}

\nod A Solovay-style splitting theorem is provable:\\

{\sc Corollary 2.25}
Let $\kp$ be $\Pi^1_{n}$-indescribable, $n\geq 1$. Then any $n+1$-stationary subset of $\kp$ can be split into $\kp$ many disjoint $n+1$-stationary sets.\\


\nod This will be extended for transfinite $\gamma$ once we have definitions of $\Pi^{1}_{\gamma}$-indescribability. 




Section 3 contains one of the central results, where
we prove the existence of certain $\Bon$ sequences which give a characterisation on the $\gm$-stationary reflecting cardinals in $L$ in terms of indescribability:
\begin{thm*}[\ref{squaregm}] $(V=L)$
Let $\gm<\kp$ be an ordinal and $\kp$ be $\Sigma^1_\gm$-indescribable but not $\Pi^1_{\gm}$-indescribable, and let $A\subseteq\kp$ be $\gm$-stationary. Then there is $E_{A}\sbs A$ and a $\Bon^{<\gm}$ sequence $S$ on $\kp$ such that $E_{A}$ is $\gm$-stationary in $\kp$ and $S$ avoids $E_{A}$.
Thus $\kp$ is not $\gm$-reflecting.
\end{thm*}
Section 3.1  gives the key tool for this proof which here serves as the version of the primary club sequences derived from $\Sigma_{n}$-hulls which occur in the fine structural proof of $\Bon$: the notion of a \emph{trace}. 
In Section 3.2 we prove the existence of certain $\Bon^\gm$ sequences for finite $\gm$ (Theorem \ref{squaren}). Although this result falls out as a corollary from Theorem \ref{squaregm} in section 3.4, the proof of Theorem \ref{squaren}  gives the key moves of the later proof while in the more familiar context of $\Pi^1_n$-indescribability and without the added complication that limit cases necessitate.  
Transfinite $\Pi^1_\gm$-indescribability for infinite $\gm$, is described in Section 3.3. This is essentially that of \cite{Sharpe-Welch2012}.
We also need to change our type of $\Bon$ sequence to deal with limit cases. In section 3.4 we give the main result of Theorem \ref{squaregm}  and conclude:

\begin{col*}[\ref{Cor3.25}]
$(V=L)$ A regular cardinal is $\Pi^1_\gm$-indescribable iff it reflects $\gm$-stationary sets.\footnote{ Using an alternate definition of $\Pi^1_\gm$-indescribability,  Bagaria in \cite{B2TAMS} proves a version of this result - however these two notions of indescribability should be equivalent. Bagaria does not define or use $\Bon$-sequences, so we have a different characterisation here which then yields this corollary.}

\end{col*}

\nod Jensen's theorem is the above with $\gamma =1$. It should be noted here that the proof of Jensen's theorem required a very close analysis of condensation and satisfaction in $L$ - Jensen used fine-structure, Beller and Litman \cite{BL} used Silver machines.
However, we shall see that we do not need such a fine analysis for $n\geq 1$ constructions - essentially this will be replaced with the coarser analysis using traces and filtrations.


In the first section of Section 4
we extend the results of Section 2.2 on the relationship between $\Pi^1_\gm$-indescribability and the $\gm$-club filter to transfinite $\gm$, including the splitting of stationary sets. This is fairly straight forward, proceeding very much as in the finite case once we have the relevant lemmas. In Section 3.5.2 we generalise the notion of non-threaded square and show (Theorem \ref{nonthreaded}) that it is closely related to $\gm$-stationarity: if the $\gm$-club filter is normal and $\Bon^\gm(\kp)$ holds then $\kp$ is not $\gm+1$-reflecting. Theorem \ref{da} is the main result of Section 4 and shows that under certain assumptions, $\gm$-stationarity is downwards absolute to $L$. A simpler, although weaker, statement of downward absoluteness is the following corollary.
\begin{col*}[\ref{dasimple}]
Assume that for any ordinal $\gm$ and any $\gm$-reflecting regular cardinal $\kp$, the $\gm$-club filter on $\kp$ is normal. Then if $\kp$ is regular and $S\sbs\kp$ is $\eta$-stationary with $S\in L$ we have $(S\text{ is $\eta$-stationary in }\kp)^L$.
\end{col*}
This extends the result of Magidor in \cite{magidorrss} which essentially gives the case $\gm=2$. As well as being an interesting result in itself, Theorem \ref{da} shows (together with the main result of Section 5) that if the of existence a $\gm$-reflecting cardinal plus some mild assumptions is consistent, then so is the existence of a $\Pi^1_\gm$-indescribable cardinal.

Finally, in Section 5 we use $\gm$-stationarity to generalise the notion of ineffability and the combinatorial principle $\diam$. $\gm$-ineffability is introduced in section 5.1 and we show that many results about ineffable cardinals generalise well to this context. $\gm$-ineffable cardinals are shown to satisfy the assumptions needed for Theorem \ref{da}, and hence (Theorem \ref{ineffda}) any $\gm$-ineffable cardinal is $\gm$-ineffable in $L$. In section 5.2 we define generalisations of $\diam$ and $\diam^*$ and look into their connection with $\gm$-ineffability. In particular we show that in $L$, where $\diam^{*\gm}$ holds can be entirely characterised in terms of ineffability (Corollary \ref{ineffdiamcol}). Although many of the proofs in this chapter lift straight from the standard cases of club and stationary set for successor cases, the limit cases are generally a different matter.




\subsection{Preliminaries and Notation}

Our set theoretical notation and definitions are standard and we refer the reader to the standard texts  \cite{J03}, or \cite{kanamori} for them.  We list here a few of these to fix notation. We assume the axioms of $ZFC$ throughout. 

We use the following abbreviations for classes: $On$ denotes the class of all ordinals, $Card$, the class of cardinals, $Sing$, the singular ordinals and $Reg$, the regular cardinals. For an ordinal $\alp$, $lim(\alp)$ abbreviates the statement that $\alp$ is a limit ordinal, and $LimOrd$ is the class of such ordinals. For a set $X$ and a cardinal $\lambda$, $|X|$ is the cardinality of $X$, $[X]^\lmb=\{Y\sbs X:|Y|=\lmb\}$ and $^{\lmb}X$ is the set of functions from $\lmb$ into $X$. Similarly $[X]^{<\lmb}=\{Y\sbs X:|Y|<\lmb\}$ and $^{<\lmb}X$ is the set of functions with from some $\bt<\lmb$ into $X$. 


Angular brackets are used to denote ordered tuples or sequences, as in $\la X,\in\ra$ or $\la C_\alp:\alp<\kp\ra$. For two sequences $p$ and $q$, $p^\frown q$ denotes their concatenation, and $lh(p)$ is the length of $p$.\\ 

Indescribability is a central concern of this paper.

\begin{defn}
An uncountable cardinal $\kp$ is \emph{$\Pi^1_n$-indescribable} if for any $R\sbs V_\kp$ and $\Pi^1_n$ formula $\phi$ such that $\la V_\kp,\in, R\ra\vDash\phi$, there is some $\alp<\kp$ such that $$\la V_\alp,\in, R\cap V_\alp\ra\vDash\phi.$$
\end{defn}
\nod This turns out to be closely connected to $n$-stationarity.

Weakly compact cardinals can be defined in many equivalent ways. We list here the two relevant to this paper. For a proof of the equivalence and some further characterisations see for example \cite{De} V.1.3.

\begin{defn}
An uncountable cardinal $\kp$ is \emph{weakly compact} if any of the following hold
\begin{enumerate}
\item Whenever $f:[\kp]^2\rightarrow 2$ there is an unbounded set $X\subseteq \kp$ such that $|f``[X]^2|=1$.
\item $\kp$ is $\Pi^1_1$-indescribable.
\end{enumerate}
\end{defn}


A variation on characterisation (1) of weakly compacts gives the notion of ineffable cardinal, a definition which is easily adapted to give the new notion of $\gm$-ineffability described in Section 6.
\begin{defn}
A regular, uncountable cardinal $\kp$ is \emph{ineffable} iff whenever $f:[\kp]^2\rightarrow 2$ there is a stationary set $X\subseteq \kp$ such that $|f``[X]^2|=1$.
\end{defn}





\begin{lem}[Hanf-Scott]\label{measurables}
A measurable cardinal is $\Pi^2_1$-indescribable and hence $\Pi^1_n$-indescribable for any $n<\om$.

Furthermore, if $U$ is a (non-principal) normal ultrafilter on $\kp$, $R\sbs V_\kp$ and $\la V_\kp,\in R\ra \vDash \phi$ where $\phi$ is $\Pi^2_1$ then $$\{\alp<\kp:\la V_\alp,\in, R\cap V_\alp\ra\vDash \phi\}\in U.$$\end{lem}


\section{Generalising closed and unbounded sets}


\subsection{Definitions and Basic Properties}

We define here the central notions of $\gm$-club set and $\gm$-stationary set.  The notion of $\gm$-club is new, and helps us to generalise many of the basic properties of stationary sets. In order for $\gm$-clubs to do what we want, we also define a notion of $\gm$-reflecting ordinal which we restrict ourselves to when defining $\gm+1$-stationary sets - just as we only define stationary sets at ordinals of uncountable cofinality.

\begin{defn}
Let $\kp$ be an ordinal and $S$, $C$ sets of ordinals. We define by simultaneous induction:
\begin{enumerate}[label= (\arabic*)]
\item $S$ is \emph{0-stationary} in $\kp$ if $S$ is unbounded in $\kp$.
\item $C$ is \emph{$\gm$-stationary-closed} if for any $\alp$ such that $C$ is $\gm$-stationary in $\alp$, we have $\alp\in C$.
\item $C$ is \emph{$\gm$-club} in $\kp$ if $C$ is $\gm$-stationary-closed below $\kp$ and $\gm$-stationary in $\kp$.
\item $\kp$ is \emph{$\gm$-reflecting} if for any $A, B\subseteq \kp$ with $A$ and $B$ $\gm$-stationary in $\kp$ there is some $\alp<\kp$ such that $A$ and $B$ are $\gm$-stationary in $\alp$.
\item $S$ is \emph{$\gm$-stationary} in $\kp$ if for all $\eta<\gm$ we have that $\kp$ is $\eta$-reflecting and for every $C$ which is $\eta$-club in $\kp$, $S\cap C\neq \emp$.
\end{enumerate}
\end{defn}

It is easy to see that our ordinary notions of club and stationary sets are the $0$-clubs and $1$-stationary sets. The $1$-clubs in $\kp$ are then stationary-closed sets which are also stationary in $\kp$. We shall see shortly that for $\gm=\eta+1$, (5) reduces to the usual definition of stationary sets in terms of clubs: $S$ is $\eta+1$ stationary if it intersects every $\eta$-club. For limit $\gm$ it is easy to see that $S$ is $\gm$-stationary if for every $\eta<\gm$, $S$ is $\eta$-stationary. The requirement of (4), a simple reflection property, is needed to ensure the $\eta$-clubs form a filter - for $0$-reflecting this is just having uncountable cofinality.
\\

The following notation from \cite{BaMS} will be very useful in exploring these concepts:
\begin{defn}
For a set $S\subseteq\kp$ we set $$d_\gm(S)=\{\alp<\kp:S\cap\alp\text{ is $\gm$-stationary in }\alp\}.$$
\end{defn}
This is a version of \emph{Cantor's derivative operator}, giving the limit points of a set in a certain topology (see \cite{B2TAMS}).
Thus $d_0(S)$ is the set of limit points of $S$, $d_1(S)$ is the set of points below which $S$ is stationary, \etc

\begin{rem}
Using this notation we see that (2), $\gm$-stationary-closure, is simply the condition that $d_\gm(C)\subseteq C$. Further (4) can be restated as  ``$\kp$ is \emph{$\gm$-reflecting} if for any $A, B\subseteq \kp$ with $A$ and $B$ $\gm$-stationary in $\kp$, $d_\gm(A)\cap d_\gm(B)\neq\emp$".
\end{rem}

For the rest of this section we look at some basic properties of such sets. We shall often omit the ordinal in which a set is stationary where that is obvious from the context, for example we shall say $S$ is $\gm$-stationary if $S$ is $\gm$-stationary in $sup(S)$. In a slight abuse of notation we shall generally take $d_\gm(S)$ to exclude $sup(S)$. 


\begin{prop}\label{implications}
Let $A\subseteq\kp$ and $\gm'<\gm$. Then:\\
(i) $A$ is $\gm$-stationary $\Rightarrow A$ is $\gm'$-stationary. Thus $d_\gm(A)\sbs d_{\gm'}(A)$.\\
(ii) If $\kp$ is $\eta$-reflecting for all $\eta<\gm$, then $A$ is $\gm'$-club $\Rightarrow A$ is $\gm$-club. \\
\end{prop}
{\sc Proof:}
(i) is clear by definition.
We use (i) to show (ii): fix $\gm$ and $\gm'<\gm$. By (i), $\gm'$-stationary closure implies $\gm$-stationary closure as if $C$ is $\gm'$-stationary closed then $d_\gm(C)\sbs d_{\gm'}(C)\sbs C$ . Suppose $C$ and $C'$ are $\gm'$-club. Then as $\kp$ is  $\gm'$-reflecting $d_{\gm'}(C)\cap d_{\gm'}(C')\neq\emp $ and $d_{\gm'}(C)\cap d_{\gm'}(C')\sbs C\cap C'$. As $C'$ was an arbitrary $\gm'$-club, we have $C$ is $\gm'+1$ stationary. Inductively, therefore, $C$ must be $\gm$-stationary and hence $\gm$-club.
\qed




\begin{rem}
By (ii) above to see that a set is $\gm+1$-stationary we need only check that it intersects every $\gm$-club (and that $\kp$ is $\gm$-reflecting).
\end{rem}

The following three propositions we prove together by  a simultaneous induction.

\begin{prop}\label{P3}
If $S$ is $\gm$-stationary and $C$ is $\gm'$-club in $\kp$ for some $\gm'<\gm$ then $S\cap C$ is $\gm$ stationary.
\end{prop}

\begin{prop}\label{P1}
If $\kp$ is $\gm$-reflecting and $A$, $B$ are $\gm$-stationary subsets of $\kp$ then $d_\gm(A)\cap d_\gm(B)$ is $\gm$-stationary.
\end{prop}

\begin{prop}\label{P2}
If $\kp$ is $\gm$-reflecting and $C^1$ and $C^2$ are $\gm$-club in $\kp$ then $C^1\cap C^2$ is $\gm$-club in $\kp$. Thus the $\gm$-clubs generate a filter on $\kp$.
\end{prop}

{\sc Proof:}
Fix $\gm$ and suppose Proposition \ref{P2} holds for every $\gm'<\gm$. 
First we show Proposition \ref{P3} holds for $\gm$. If $\gm=0$ this is vacuous. So suppose $\gm>0$ and let $S$ be an $\gm$-stationary subset of $\kp$ and $\gm,\eta'<\gm$. Let $C$ be $\gm'$-club in $\kp$ and let  $C'$ be $\eta'$-club in $\kp$. Set $\eta=max\{\gm',\eta'\}$. By Proposition \ref{implications} $C$ and $C'$ are both $\eta$-club. Now by Proposition \ref{P2} we have $C\cap C'$ is $\eta$-club. Thus $(S\cap C)\cap C'=S\cap (C\cap C')\neq\emp$, so as $\eta'$ and $C'$ were arbitrary $S\cap C$ is $\gm$-stationary.

Now supposing we have Proposition $\ref{P3}$ for $\gm$ we show Proposition $\ref{P1}$ holds for $\gm$.  Let $\kp$ be $\gm$-reflecting, $A,B\subseteq\kp$ be $\gm$-stationary and let $\gm'<\gm$ with $C$ $\gm'$-club. By Proposition $\ref{P3}$ we have $C\cap A$ is $\gm$-stationary, so as $\kp$ is $\gm$-reflecting $d_{\gm}(C\cap A)\cap d_{\gm}(B)\neq\emp$. But $d_{\gm}(C\cap A)\subseteq C$ as $C$ is $\gm$-stationary in $\alp$ $\Rightarrow$ $C$ $\gm'$-stationary in $\alp$ $\Rightarrow \alp\in C$,and clearly $d_{\gm}(C\cap A)\subseteq d_{\gm}(A)$. Thus $C\cap d_{\gm}(A)\cap d_{\gm}(B)\sps d_{\gm}(C\cap A)\cap d_{\gm}(B)\neq\emp$. Hence $d_\gm(A)\cap d_\gm(B)$ is $\gm$-stationary.

Finally we show Proposition \ref{P1} implies Proposition \ref{P2}.  Take $C^1$ and $C^2$ to be $\gm$-club in $\kp$. By Proposition \ref{P1} we have $d_\gm(C^1)\cap d_\gm(C^2)$ is $\gm$-stationary so as $d_\gm(C^1)\cap d_\gm(C^2)\sbs C^1\cap C^2$ we just need to show $\gm$-stationary closure. But this is simple as if $C^1\cap C^2$ is $\gm$-stationary below $\alp<\kp$ then $C^1$ and $C^2$ are both $\gm$-stationary below $\alp$ so by the $\gm$-stationary closure of $C^1$ and $C^2$ we must have $\alp\in C^1\cap C^2$.
\qed

\begin{prop}\label{limclub}
If $\kp$ is $\gm$-reflecting and $C$ is $\gm$-club in $\kp$ then $d_\gm(C)$ is also $\gm$-club in $\kp$.
\end{prop}

{\sc Proof:}
By Proposition \ref{P1} we have $d_\gm(C)$ is $\gm$-stationary. To see we have closure, if $d_\gm(C)$ is $\gm$-stationary below $\alp$ then, as $d_\gm(C)\sbs C$ we have $C$ is $\gm$-stationary below $\alp$, so by stationary closure of $C$ we have $\alp\in d_n(C)$.
\qed

The following is a useful way of building $\gm$-stationary sets.

\begin{lem}\label{nstatunion}
Let $A\sbs \kp$ be $\gm$-stationary in $\kp$ and $\la A_\alp:\alp\in A\ra$ be a sequence of sets such that each $A_\alp$ is $\gm$-stationary in $\alp$. Then $\bigcup\{ A_\alp:\alp\in A\}$ is $\gm$-stationary in $\kp$.
\end{lem}

{\sc Proof:}
This is clearly true for $\gm=0$ so suppose $\gm>0$ and let $\gm'<\gm$ and $C\sbs\kp$ be $\gm'$-club in $\kp$. Then as $A$ is $\gm$-stationary in $\kp$ we must have $\kp$ is $\gm'$-reflecting so by Proposition \ref{limclub} $d_{\gm'}(C)$ is $\gm'$-club. Then we can find some $\alp$ in $d_{\gm'}(C)\cap A$. Now $C$ is $\gm'$-club in $\alp$ so as $A_\alp$ is $\gm$-stationary in $\alp$ we must have $C\cap A_\alp\neq\emp$. Thus $C\cap\bigcup\{ A_\alp:\alp\in A\}\neq\emp$ and we're done.
\qed

We can now prove a stronger result than \ref{limclub}:
\begin{prop}\label{statclub}
If $\kp$ is $\gm$-reflecting and $S$ is $\gm$-stationary in $\kp$ then $d_\gm(S)$ is $\gm$-club in $\kp$.
\end{prop}

{\sc Proof:}
As $\kp$ is $\gm$-reflecting $d_\gm(S)$ is $\gm$-stationary by Proposition \ref{P1}. To show closure suppose $d_\gm(S)$ is $\gm$-stationary in $\alp$. Then $S\cap\alp=\bigcup\{S\cap\bt:\bt\in d_\gm(S)\cap\alp\}$ with each $S\cap\bt$ being $\gm$-stationary in $\bt$, so by Lemma \ref{nstatunion} $S$ is $\gm$-stationary in $\alp$. Hence $\alp\in d_\gm(S)$ and we have $\gm$-stationary closure.
\qed

The following shows that there are also many ordinals below which a set is not $\gm$-stationary: 
\begin{prop}\label{prop} 
If $A$ is $\gm$-stationary in $\kp$ then $$A\setminus d_\gm(A)=\{\alp\in A:A\cap\alp \text{ is not $\gm$-stationary}\}$$ is $\gm$-stationary.
\end{prop}
Note that $A\setminus d_\gm(A)$ can never be $\gm+1$-stationary as if $\kp$ is $\gm$-reflecting then Proposition \ref{statclub} gives us that $d_\gm(A)$ is $\gm$-club.

{\sc Proof:}
Let $\gm'<\gm$ and $C$ be $\gm'$-club in $\kp$. Then $d_{\gm'}(C)$ is $\gm'$-club so we can find $\alp$ minimal in $A\cap d_{\gm'}(C)$. Then $C$ is $\gm'$-club in $\alp$. If $\alp$ is not $\gm'$-reflecting then $A$ cannot be $\gm$-stationary in $\alp$ so we're done. If $\alp$ is $\gm'$-reflecting then by Proposition \ref{limclub} $d_{\gm'}(C)$ is also $\gm'$-club in $\alp$. But $A\cap \alp\cap d_{\gm'}(C)=\emp$ so $A\cap\alp$ is not $\gm$-stationary.
\qed

We can now show how closely this definition of $\gm$-stationarity is related to that in \cite{BaMS}, and to notions of reflection:

\begin{prop}\label{s-stat}
Let $\kp$ be $\eta$-reflecting for all $\eta<\gm$ and $S\sbs \kp$. Then $S$ is $\gm$-stationary iff for any $\eta<\gm$ and any $\eta$-stationary $A\sbs \kp$ we have $d_\eta(A)\cap S\neq\emp$.

\end{prop}

{\sc Proof:}
$(\Rightarrow)$ is clear as by Proposition \ref{statclub} $d_\eta(A)$ is $\eta$-club. 
$(\Leftarrow):$ if $C$ is $\eta$-club then $C$ is $\eta$-stationary so we have $d_\eta(C)\cap S\neq \emp$, but $d_\eta(C)\sbs C$.
\qed

The following is an easy consequence:
\begin{prop}
If $\kp$ is $\gm$-reflecting then $\kp$ is $\eta$-reflecting for all $\eta\leq \gm$.
\end{prop}

Comparing this characterisation to the definition in \cite{BaMS} we see the only difference between our $\gm$-stationary sets and those defined in \cite{BaMS} is that we always require this \emph{simultaneous} reflection. We do see, however, that our definition of $\gm$-stationary is equivalent to the definition $\gm$-s-stationary given in \cite{B2TAMS}:

\begin{prop}
Define $\gm$-s-stationary inductively as follows, starting with $0$-s-stationary being unbounded. Let $S$ be $\gm$-s-stationary in $\kp$ if for every $\eta<\gm$ and $A$, $B$ which are $\eta$-stationary in $\kp$ there is some $\alp\in S$ such that $A\cap\alp$ and $B\cap\alp$ are both $\eta$-s-stationary in $\alp$. 
Then a set $S$ is $\gm$-s-stationary in $\alp$ iff $S$ is $\gm$-stationary in $\alp$.
\end{prop}

{\sc Proof:}
By induction on $\gm$. Suppose we have for $\eta<\gm$ and any ordinal $\alp$ that a set is $\eta$-stationary iff $\eta$-s-stationary. 

Let $S$ be $\gm$-s-stationary in $\kp$. First we see that $\kp$ must be $\eta$-reflecting for any $\eta<\gm$: if $A$ and $B$ are $\eta$-stationary then they are $\eta$-s-stationary by the inductive hypothesis, and hence there is some $\alp\in S$ such that $A\cap\alp$ and $B\cap\alp$ are both $\eta$-s-stationary in $\alp$. But then $A\cap\alp$ and $B\cap\alp$ are both $\eta$-stationary in $\alp$, so $\kp$ is $\eta$-reflecting. Now suppose $C$ is $\eta$ club for some $\eta<\gm$. Then $C$ is $\eta$-s-stationary so there is some $\alp\in S$ such that $C\cap\alp$ is $\eta$-s-stationary in $\alp$, \ie $C\cap\alp$ is $\eta$-stationary in $\alp$. So by $\eta$-stationary closure, $\alp\in C\cap S$. As $C$ was arbitrary, we must have $S$ is $\gm$-stationary.

Now suppose $S$ is $\gm$-stationary and let $\eta<\gm$ and $A$ and $B$ be $\eta$-s-stationary. By the inductive hypothesis, $A$ and $B$ are both $\eta$-stationary, and so by Proposition \ref{statclub} $d_\eta(A)$ and $d_\eta(B)$ are both $\eta$-club. Then by Proposition \ref{P2} $d_\eta(A)\cap d_\eta(B)$ is $\eta$-club, and hence we can find $\alp\in S\cap d_\eta(A)\cap d_\eta(B)$. But $A$ and $B$ are both $\eta$-stationary below such an $\alp$, and hence $A$ and $B$ are both $\eta$-s-stationary below such an $\alp$, so we have that $S$ is $\gm$-s-stationary.
\qed

 It is easy to show (see Proposition \ref{pin}) that any $\Pi^1_n$-indescribable is $n$-reflecting. So a simple induction shows that we get this result for our definitions too:

\begin{thm}{\em\cite{BaMS}}\label{BaMSthm}
Assume $V=L$. Then a regular cardinal $\kp$ is $n$-reflecting iff $\kp$ is $\Pi^1_n$-indescribable.
\end{thm}

In Section 3
we shall give an alternative proof of this, and after defining $\Pi^1_\gm$-indescribability in section 3.3 we shall extend it to replace $n$ with any ordinal $\gm<\kp$. Nevertheless, our result is still of further interest as we are proving the existence of certain $\Bon$-sequences, which are not constructed in either \cite{BaMS} or \cite{B2TAMS}.  

One question now is where the $\gm$-stationary sets can occur if $V\neq L$. We have that if a cardinal is $\Pi^1_1$-indescribable then it is $2$-stationary (\ie $1$-reflecting), and the above shows that it is consistent for these to be the only $2$-stationary regular cardinals. However, known results show that it is also consistent, relative to certain large cardinal assumptions, that non-weakly compact cardinals  be stationary reflecting, and indeed $1$-reflecting.

It is well known that a regular cardinal cannot be stationary reflecting unless it is either the successor of a singular cardinal or weakly inaccessible. It is an observation that $E^{\lambda^{+}}_{\lambda} =_{{df}} \{\alpha < \lambda^{+}\mid cf(\alpha)=\lambda\}$ does not reflect for $\lambda \in Reg$.
Note that singular ordinals are not so interesting in this context as $\gm$-stationarity for subsets of a singular reduces to $\gm$-stationarity for subsets of its cofinality:
\begin{prop}
Let $\alp$ be a singular ordinal and $C\sbs\alp$ be club, with $\pi:\la C,\in\nolinebreak \ra\cong\la ot(C),\in\ra$. Then for $\gm\geq 1$ any $S\sbs\alp$ is $\gm$-stationary in $\alp$ iff $\pi``S\cap C$ is $\gm$-stationary in $ot(C)$. 
Hence $\alpha$ is $\gm$-reflecting if and only if  $\cof(\alpha)$ is.

\end{prop}
{\sc Proof:}
This is proven by induction. Fix $\alp$ and $\gm$ and suppose the claim is true for all $\bt<\alp$ and for all $\eta<\gm$. Fix $C$, $S$ and $\pi:\la C,\in\ra\cong\la ot(C),\in\ra$. If $D$ is $\eta$-club in $\alp$ then by the inductive hypothesis and the closure of $C$, $\pi``D\cap C$ is $\eta$-club in $ot(C)$, and hence if $\pi``S\cap C$ is $\gm$-stationary in $ot(C)$ then $S$ is $\gm$-stationary in $\alp$. Similarly if $D$ is $\eta$-club in $ot(C)$ then $\pi^{-1}``D$ is $\eta$-club in $\alp$, so if $S$ is $\gm$-stationary in $\alp$ then $\pi``S\cap C$ is $\gm$-stationary in $ot(C)$.
\qed

A cardinal can be stationary reflecting without \emph{simultaneously} reflecting stationary sets (\ie without being $1$-reflecting) as shown in \cite{mshelah}, and in fact the consistency strength of the existence of a stationary reflecting cardinal is much less than the existence of a $1$-reflecting cardinal. This follows from results of Mekler-Shelah in \cite{mshelah} and Magidor in \cite{magidorrss}. In the latter it is shown that if a regular cardinal $1$-reflecting then it is weakly compact in $L$, and thus if the existence of a $1$-reflecting cardinal is consistent then so is the existence of a weakly compact cardinal. we shall generalise this result on downward absoluteness in Section 3.6 (although we have to add some assumptions there).

Kunen has shown \cite{kunen} that it is consistent relative to the existence of a weakly compact cardinal that there is a stationary reflecting cardinal that is not weakly compact - this is proven by giving a forcing which adds a Suslin tree to a cardinal $\kp$ that is weakly compact in the ground model. It is clear from the proof that in this model $\kp$ reflects any two stationary sets simultaneously and is thus $1$-reflecting. It is also easy to see that the $1$-club filter is normal there (see Definition \ref{completedf}). 

Magidor in \cite{magidorrss} starts from the much stronger assumption that there are infinitely many supercompact cardinals, and produces a forcing model in which $\aleph_{\om+1}$ reflects stationary sets. Again, it is clear from the proof that $\aleph_{\om+1}$ is in fact $1$-reflecting. Here, however, the $1$-club filter is not even countably complete, as shown by the following easy argument. 
For each $n<\om$ let $C_n=\{\alp<\aleph_{\om+1}:cf(\alp)\geq\aleph_{n}\}$. Then each $C_n$ is $1$-club - it is clearly stationary, and cannot be stationary below any ordinal of cofinality $\leq\aleph_{n}$ so is also stationary closed. But $\bigcap_{n<\om}C_n=\emp$. 
It is also easy to see that $\aleph_{\om+1}$ cannot be $2$-reflecting, as $\aleph_{\om+1}$ is the \emph{minimum} ordinal that can be $2$-stationary so there is nowhere for $2$-stationary sets to reflect to.


\subsection{Results at Indescribables}

\subsubsection{The Club and Indescribability Filters}

In this section we shall be focusing on the relationship between $\Pi^1_n$-indescribability and our generalised notions of club and stationarity, in particular the $n$-club filter. For this reason, we shall be mostly restricting to finite levels of the club and stationary set hierarchy here. We shall extend these results to the transfinite in section 3.5, 
after we have introduced a notion of $\Pi^1_\gm$-indescribability.  
Analysing this relationship will allow us to generalise some deeper properties of stationary sets and the club filter to $n+1$-stationary sets and the $n$-club filter at a $\Pi^1_n$-indescribable cardinal. Thus by Theorem \ref{BaMSthm} the generalisation is full in $L$, but limited if we have $n$-reflecting cardinals which are not $\Pi^1_n$-indescribable. 

We first look at the relationship between the $n$-club filter and the $\Pi^1_n$-indescribability filter.
\begin{defn}\label{completedf}
For a $\gm$-reflecting cardinal $\kp$ we denote by $\Cf^\gm(\kp)$ the $\gm$-club filter on $\kp$:
$$\Cf^\gm(\kp):=\{X\sbs\kp:X \text{ contains a $\gm$-club}\}.$$
\end{defn}
\begin{defn}[Levy]
If $\kp$ is $\Pi^1_n$-indescribable let $\F^n(\kp)$ denote the $\Pi^1_n$-indescribable filter on $\kp$:
\begin{displaymath}
\begin{split}
\F^n(\kp):=\{X\sbs\kp: & \text{ for some $\Pi^1_n$ sentence $\varphi$ with parameters such that }V_\kp\vDash\varphi,\\
& V_\alp\vDash\varphi\rightarrow\alp\in X\}.
\end{split}
\end{displaymath}
\end{defn}

\begin{rem} (1) The statement "$X$ is $n$-stationary in $\kp$" is $\Pi^1_n$ expressible over $V_\kp$. This can be seen by induction - clearly ``$X$ is unbounded" is $\Pi^1_0$.  ``$X$ is $n+1$-stationary" is equivalent (by Proposition \ref{s-stat}) to "$X$ is $n$-stationary $\wedge\forall S,T(S,T$ are $n$-stationary $\rightarrow\exists\alp\in X~S\cap\alp,T\cap\alp$ are $n$-stationary in $\alp)$" - assuming $n$-stationarity is $\Pi^1_n$ expressible this is clearly $\Pi^1_{n+1}$. \\(2) By a result of Levy If $\kp$ is $\Pi^1_n$ indescribable then $\F^n(\kp)$ is normal and $\kp$-complete.
 
\end{rem}

\begin{prop}\label{pin}
If $\kp$ is a $\Pi^1_n$-indescribable cardinal then $\kp$ is $n$-reflecting and furthermore any set $X\in \F^n(\kp)$ is $n+1$-stationary.
\end{prop}

{\sc Proof:}
Suppose $\kp$ is $\Pi^1_n$-indescribable. We have seen that being $n$-stationary in $\kp$ is $\Pi^1_n$ expressible over $V_\kp$, so for any $S,T\sbs\kp$ if we have $V_\kp\vDash ``S,T$ are $n$-stationary$"$ then for some $\alp<\kp$ we have $V_\alp\vDash ``S\cap\alp,T\cap\alp$ are $n$-stationary$"$ and hence $\kp$ is $n$-reflecting. For $X\in \F^n(\kp)$ let $\varphi$ be a $\Pi^1_n$ formula such that $V_\kp\vDash\varphi$ and $V_\alp\vDash\varphi\rightarrow\alp\in X$. For an $n$-club $C$, taking $\psi=C$ is $n$-stationary $\wedge\, \phi$, we have that $\psi$ is $\Pi^1_n$ and hence for some $\alp<\kp$, $V_\alp\vDash \psi$. But then $V_\alp\vDash\varphi$ so $\alp\in X$ and $C\cap\alp$ is $n$-stationary so by $n$-stationary closure, $\alp\in C$. Thus $X$ meets any $n$-club and so $X$ is $n+1$-stationary.
\qed


First we show that at a $\Pi^1_n$-indescribable cardinal the $n$-club filter is normal.

\begin{lem}\label{filtersV1}
For any $n<\om$ and $\Pi^1_n$-indescribable cardinal $\kp$ we have
\begin{enumerate}
\item $\Cf^n(\kp)\subseteq \F^n(\kp)$,
\item $\Cf^n(\kp)$ is normal (and hence $\kp$-complete) .
\end{enumerate}
\end{lem}
{\sc Proof:}
(1) Let $C$ be $n$-club in $\kp$. Then $d_n(C)\subseteq C$ and $d_n(C)=\{\alp<\kp:V_\alp\vDash``C\text{ is $n$-stationary}"\}$. As being $n$-stationary is $\Pi^1_n$ we see $d_n(C)\in\F^n$ so $C\in\F^n$.

(2) Let $\la C_\alp:\alp<\kp\ra$ be a sequence of $n$-clubs. Defining $C:=\bigtriangleup_{\alp<\kp}C_\alp$ we have $C\in\F^n$ (using the normality of $\F^n$ and (1)). As each element of $\F^n$ is $n+1$-stationary we thus have $C$ is $n$-stationary. To show closure suppose $C$ is $n$-stationary in $\alp<\kp$. As $C\cap\alp=\{\bt<\alp:\bt\in\bigcap_{\gm<\bt}C_\gm\}$ we see that for any $\bt<\alp$ we have $(\bt,\alp)\cap C\subseteq C_\bt$ so $C_\bt\cap\alp$ includes an end-segment of an $n$-stationary set and thus $C_\bt$ is n-stationary in $\alp$. Then $\alp\in C_\bt$ as $C_\bt$ is $n$-club. Thus $\alp\in\bigcap_{\bt<\alp}C_\bt$, \ie $\alp\in C$.
\qed

\begin{col}(Fodor's Lemma for $n$-stationary sets) 
If $\kp$ is $\Pi^1_n$-indescribable and $A\subset\kp$ is $n$-stationary, then for any regressive function $f:A\rightarrow \kp$ there is an $n$-stationary $B\subseteq A$ such that $f$ is constant on $B$.
\end{col}

\begin{thm}\label{filtersL1}
At any weakly compact cardinal $\kp$, $\Cf^1(\kp)=\F^1(\kp)$ and assuming $V=L$ we have $\Cf^n(\kp)=\F^n(\kp)$ for any $n<\om$ and $\Pi^1_n$-indescribable $\kp$.
\end{thm}
{\sc Proof:}
Firstly we show $\Cf^1(\kp)=\F^1(\kp)$; the later part will be an induction using that all $n$-reflecting cardinals are $\Pi^1_n$-indescribable, which holds in $L$.\\
Let $\forall X\varphi(X)$ be a $\Pi^1_1$ sentence (possibly with parameters) such that $V_\kp\vDash \forall X\varphi(X)$. Note that $\phi(X)$ is $\Sigma^1_0$ so $\neg\phi(X)$ is $\Pi^1_0$. We define $A\subseteq\{\alp<\kp:V_\alp\vDash\forall X\varphi(X)\}$ and show $A$ is $1$-club. Set $A:=\{\alp<\kp:\alp\text{ is inaccessible }\wedge V_\alp\vDash\forall X\varphi(X)\}$. Then as inaccessibility is $\Pi^1_1$ expressible, $A\in\F^1(\kp)$ and so $A$ is stationary. To show $A$ is stationary-closed take $\alp<\kp$ a limit of inaccessibles.\\
(i) If $\alp$ is regular then $\alp$ is inaccessible. Suppose $\alp\notin A$. Then $V_\alp\vDash\neg \varphi(Y)$ for some $Y\subseteq\alp$ so by the inaccessibility of $\kp$, $\{\bt<\alp:V_\bt\vDash\neg \varphi(Y\cap\bt)\}$ is club in $\alp$. Thus if $A$ is stationary in $\alp$ then $\alp\in A$.\\
(ii) If $\alp$ is singular then $A$ is not stationary in $\alp$: set $\lmb=cf(\alp)$ and take $C$ club in $\alp$ with $ot(C)=\lmb$. Then no inaccessible above $\lmb$ can be a limit point of $C$, so $d_0(C)\cap A$ is bounded in $\alp$.
Thus $A$ is stationary-closed and hence $1$-club, so we have $\Cf^1(\kp)=\F^1(\kp)$.

Now we assume $V=L$. Suppose $\kp$ is $\Pi^1_{n+1}$ indescribable and for all $\alp<\kp$ if $\alp$ is $\Pi^1_{n}$ indescribable then $\Cf^n(\alp)=\F^n(\alp)$. Let $\forall X\varphi(X)$ be a $\Pi^1_{n+1}$ sentence (possibly with parameters) such that $V_\kp\vDash \forall X\varphi(X)$. Take $A:=\{\alp<\kp:\alp\text{ is inaccessible }\wedge V_\alp\vDash\forall X\varphi(X)\}$. As before $A$ is a subset of the points where $\forall X\varphi(X)$ is reflected, and $A$ is $n+1$-stationary as $A\in\F^{n+1}(\kp)$. Also as above, if $\alp<\kp$ is singular then $A$ is not stationary (and hence not $n+1$-stationary) in $\alp$. So suppose $\alp$ is regular and $A$ is $n+1$-stationary in $\alp$. Then $A$ must be $\Pi^1_n$-indescribable as only $n$-reflecting ordinals admit $n+1$-stationary sets, and in $L$ the $n$-reflecting regular cardinals are exactly the $\Pi^1_n$-indescribables (Theorem \ref{BaMSthm}). Suppose for a contradiction $V_\alp\vDash\neg \varphi(Y)$ for some $Y\subseteq\alp$. Now $\neg \varphi(Y)$ is $\Pi^1_n$ so setting $B=\{\bt<\alp:V_\bt\vDash\neg\varphi(Y\cap\bt)\}$ we have $B\in\F^n(\alp)=\Cf^n(\alp)$ by the inductive hypothesis, so $A$ is not $n+1$ stationary in $\alp$, contradiction. Thus $A$ is $n+1$-club, so $\F^{n+1}(\kp)=\Cf^{n+1}(\kp)$.
\qed

\subsubsection{Splitting Stationary Sets}

We show that at a $\Pi^1_n$-indescribable $\kp$ each $n+1$-stationary set can be split into $\kp$ many disjoint $n+1$-stationary sets. This is a generalisation of Solovay's  result (\cite{solovay}) that any stationary subset of a regular $\kp$ can be split into $\kp$ many disjoint stationary sets, but the proof is very different. The difficult part of this is actually to show that each $n+1$-stationary set can be split into two disjoint $n+1$-stationary sets, \ie to show that the ideal of non-$n+1$-stationary sets is atomless. $\kp$ splitting then follows by an exercise in Jech \cite{J03} (ex.13, p.124). After we have introduced $\Pi^1_\gm$-indescribability in section 3.3 we shall extend the results below to $\gm+1$-stationarity (see section \ref{sectionsplitting2}).

\begin{lem}\label{split}$(n>0)$ 
If $\Cf^{n-1}(\kp)$ is $\kp$-complete then any $n$-stationary subset of $\kp$ is the union of two disjoint $n$-stationary sets.
\end{lem}

{\sc Proof:}
Let $S$ be $n$-stationary in $\kp$ and suppose $S$ is not the union of two disjoint $n$-stationary sets. Define 
$$F=\{X\subset\kp:X\cap S \text{ is $n$-stationary in }\kp\}$$

\emph{Claim:} $F$ is a $\kp$ complete ultrafilter\\
Upwards closure is clear. If $A$, $B\in F$ then we have $S\without A$ and $S\without B$ are both non-$n$-stationary sets, as $S$ cannot be split. Thus $A\cup (\kp\without S)$ and $B\cup (\kp\without S)$ are both in the $n-1$-club filter, and hence their intersection contains an $n-1$-club $C$. But then $C\cap S$ is $n$-stationary and $C\cap S\sbs A\cap B\cap S$, so $A\cap B\in F$. 
 For maximality, if $X\in F$ then as $S$ cannot be split we must have $\kp\without X\notin F$. That $X\notin F\Rightarrow \kp\without X\in F$ follows from the fact that the $n-1$-clubs form a filter. The $\kp$-completeness of $F$ follows from the $\kp$-completeness of the $n-1$-club filter, in the same way as the intersection property. Also, $F$ is clearly non-principal as it contains all end-segments.

\emph{Claim:} $F$ is normal\\
As we have shown that $F$ is a $\kp$ complete non-princpal ultrafilter on $\kp$, we have that $\kp$ is measurable and hence by Lemma \ref{measurables}, $\kp$ is $\Pi^1_{n-1}$-indescribable. Thus $\Cf^{n-1}(\kp)$ is normal by Proposition \ref{filtersV1}. Let $\la X_\alp:\alp<\kp\ra$ be a sequence of sets in $F$. Then each $S\without X_\alp$ is in the non-$n$-stationary ideal on $\kp$, so $X_\alp\cup (\kp\without S)\in\Cf^{n-1}(\kp)$. Now by the normality of $\Cf^{n-1}(\kp)$, we have, setting $$X:=\bigtriangleup_{\alp<\kp}X_\alp\cup(\kp\without S)$$ that $X\in\Cf^{n-1}(\kp)$, and so $X\cap S$ is $n$-stationary. But 
\begin{displaymath}
\begin{split}
X &=\{\alp<\kp:\forall\bt<\kp\, \,\alp\in X_\bt\cup(\kp\without S)\}\\
&=\{\alp<\kp:\alp\in\kp\without S\vee \forall\bt<\kp\,\,\alp\in X_\bt\}\\
&=(\kp\without S)\cup\bigtriangleup_{\alp<\kp}X_\alp
\end{split}
\end{displaymath}
so $X\cap S= \bigtriangleup_{\alp<\kp}(X_\alp\cap S)$, and thus $\bigtriangleup_{\alp<\kp}(X_\alp\cap S)$ is $n$-stationary and hence in $F$.\\

As we have shown that $F$ is a normal measure, we can use the second part of Lemma \ref{measurables}: for any $R\sbs V_\kp$ and $\phi$ that is $\Pi^2_1$ we have that if $\la V_\kp,\in R\ra \vDash \phi$ then 
$$\{\alp<\kp:\la V_\alp,\in, R\cap V_\alp\ra\vDash \phi\}\in F.$$
Thus setting $R=S$ and $\phi=``S$ is $n$-stationary$"$ we have $$\{\alp<\kp:S\cap\alp\text{ is $n$-stationary in }\alp\}\in F.$$
By definition of $F$ therefore, $$A=\{\alp\in S:S\cap\alp\text{ is $n$-stationary}\}$$ is $n$-stationary. But by Proposition $\ref{prop}$ we have  $$A':=\{\alp\in S:S\cap\alp \text{ is not $n$-stationary}\}$$ is $n$-stationary. This contradicts our assumption on $S$ as $A'$ and $A$ are  two disjoint, $n$-stationary subsets of $S$.
\qed

Adding the assumption of weak compactness we can now split $S$ into $\kp$ many pieces.

\begin{thm}\label{splitting}
If $\kp$ is weakly compact and $\Cf^{n-1}(\kp)$ is $\kp$-complete then any $n$-stationary subset of $\kp$ is can be split into $\kp$ many disjoint $n$-stationary sets.
\end{thm}

{\sc Proof:} 
This is essentially exercise 13 p124 in Jech \cite{J03}, using Lemma \ref{split}. We defer the proof to that for arbitrary $\gamma$ rather than $n$ at Theorem \ref{splittinggm}.
\qed

\begin{col}
Let $\kp$ be $\Pi^1_{n}$-indescribable, $n\geq 1$. Then any $n+1$-stationary subset of $\kp$ can be split into $\kp$ many disjoint $n+1$-stationary sets.
\end{col}

{\sc Proof:}
For $n\geq1$, $\Pi^1_n$-indescribability implies weak compactness and Lemma \ref{filtersV1} gives us the $\kp$-completeness of $\Cf^n(\kp)$, so we can apply Theorem \ref{splitting}.
\qed



\section{Generalising $\Bon$-sequences}



\label{Chapter4} 

Jensen proved \cite{Je72} that in $L$, a regular cardinal being weakly compact is equivalent to being stationary-reflecting, by constructing a square sequence below a non-weakly compact $\kp$ which avoids a certain stationary set. 
We aim to generalise this, replacing the notion of stationary set with $\gm$-stationary set, and defining a new notion of $\Bon$-sequence using $\gm$-clubs.

Firstly, for finite $n$, we shall use the natural generalisation of $\Bon$ notions to the context of $\gm$-clubs that follows, a $\Bon$-sequence that will witness the failure of a $\gm+1$-stationary set to reflect.

\begin{defn}\label{squaregmdef} 
Let $\gm<\kp$ be ordinals. A $\Bon^\gm$ \emph{sequence on} $\kp$ is a sequence $\la C_\alp:\alp\in d_\gm(\kp)\ra$ such that for each $\alp$:
\begin{enumerate}
\item $C_\alp$ is an $\gm$-club subset of $\alp$ and
\item for every $\bt\in d_\gm(C_\alp)$ we have $C_\bt=C_\alp\cap\bt$.
\end{enumerate}
\end{defn}

So, for instance, Jensen's characterisation in \cite{Je72} of non-weakly compact cardinals $\kp$ is a $\Bon^0$-sequence below $\kp$ which avoids a certain stationary set. We can now state our first generalisation of Jensen's Theorem:

\begin{thm}\label{squaren} $(V=L)$
Let $n\geq 0$ and $\kp$ be a  $\Pi^1_n$- but not $\Pi^1_{n+1}$-indescribable cardinal, and let $A\subseteq\kp$ be $n+1$-stationary. Then there is $E_{A}\sbs A$ and a $\Bon^n$-sequence $S$ on $\kp$ such that $E_{A}$ is $n+1$-stationary in $\kp$ and $S$ avoids $E_{A}$.\\
Thus $\kp$ is not $n+1$-reflecting.
\end{thm}

In order to prove this theorem for $n>0$ we shall need some technical machinery, which will be introduced in the next subsection. In section 3.2 we give the proof of the above theorem. 

However, we want to generalise this result further, replacing $n$ in the above with an arbitrary ordinal $\gm<\kp$. To do this, we need a concept of $\Pi^1_\gm$-indescribability, which is introduced in section 3.3.1. We shall also need a new type of $\Bon$-sequence, because although Definition \ref{squaregmdef} makes sense for infinite ordinals $\gm$, we shall need to deal with limit cases where that definition is not available. This $\Bon^{<\gm}$-sequence is introduced in 3.3.2, and then in section 3.4 we state and prove the most general version of our theorem.


\subsection{Traces and Filtrations}


We define some notation for familiar concepts.
\begin{defn}
For a transitive set $M$ together with a well-ordering of $M$ which we fix for this purpose, if $X\sbs M$, let $M\{X\}$ be the Skolem hull of $X$ in $M$ using the Skolem functions defined from the (suppressed) well-ordering.
\end{defn}

Although stated generally, we shall just use structures $M$ which are levels of the $L$-hierarchy, possibly with additional relations. The well-ordering is then the standard one of the levels of $L$.
 
\begin{defn}[Trace and Filtration]
Let $M$ be a transitive set equipped with a Skolem hull operator $M\{.\}$ and let $p\in M^{<\om}$ (we say $p$ is a \emph{parameter} from $M$) and $\alp\in M$. The \emph{trace} of $M,p$ on $\alp$ is the set
$$\mathcal{s}(M,p,\alp):=\{\bt<\alp:\alp\cap M\{p\cup\bt\cup\{\alp\}\}=\bt\}$$

The \emph{filtration of $M,p$ in} $\alp$ is the sequence $$\la M_\bt= M\{p\cup\bt\cup\{\alp\}\}:\bt\in\mathcal{s}(M,p,\alp)\ra$$
\end{defn}

\begin{rem} 
Note that the filtration is \emph{continuous} and
\emph{monotone} increasing.
If $\alp$ is a regular cardinal we shall have as usual the filtration is unbounded, and the trace will be club in $\alp$. 
\end{rem}



Let $\Gamma$ be a class of $\La_\in$ formulae such that each $\Gamma$ formula $\varphi$ has a distinguished variable $v_0$ .

\begin{defn}
A model $\la M,\in\ra$ is \emph{$\Gamma$ correct over} $\alp$ if $\alp\in M$ and for any $\Gamma$ formula $\phi(v_0,\dots v_n)$ with all free variables displayed and $A_1,\dots A_n$ from $\P(\alp)\cap M$, $$M\vDash\phi(\alp,A_1,\dots A_n) \text{ iff  }~  \phi(\alp,A_1,\dots A_n).$$
\end{defn}

\begin{rem} Note that if we set $\neg\Gamma=\{\neg\varphi:\varphi\in\Gamma\}$ then $\Gamma$-correctness is the same as $\neg\Gamma$-correctness. We shall initially be using this for $\Gamma = \Pi^{1}_{n}$ and later for $\Pi^{1}_{\gamma}$.
\end{rem}

We can thin out the trace and filtration by requiring that the hulls collapse to transitive models which are $\Gamma$ correct over $\bar{\alp}$, where $\bar{\alp}$ is the collapse of $\alp$. More formally: 
\begin{defn}
The \emph{$\Gm$-trace of $M,p$ on} $\alp$ is denoted $\mathcal{s}^{\Gamma}(M,p,\alp)$ and consists of all $\bt<\alp$ such that $\bt\in\mathcal{s}(M,p,\alp)$ and if $\pi:M\{p\cup\bt\cup\{\alp\}\}\cong N$ is the transitive collapse then $N$ is $\Gamma$-correct over $\bt=\pi(\alp)$.
\end{defn}

We now work under the assumption $V=L$ and prove some elementary properties of traces of levels of $L$.

\begin{lem}\label{tracemono}
\begin{enumerate}
\item If $\alp<\mu<\nu$ are limit ordinals with $p$ a parameter from $L_\nu$, $q$ a parameter from $L_\mu$ and $L_\nu=L_\nu\{p\cup\alp+1\}$ then there is some $\bt_0<\alp$ such that
$$\trace(L_\mu,q,\alp)\supseteq\trace(L_\nu,p,\alp)\without\bt_0$$

\item If in addition we assume $L_\nu$ and $L_\mu$ are $\Gm$-correct over $\alp$, the same holds for  the $\Gm$-trace $\trace^{\Gm}(L_\nu,p,\alp)$.
\end{enumerate}
\end{lem}

{\sc Proof:}
(1) As $L_\nu=L_\nu\{p\cup\alp+1\}$ there is $\bt_0<\alp$ such that $L_\mu,q\in L_\nu\{p\cup\bt_0\cup\{\alp\}\}$. It is easy to see that this $\bt_0$ works as $L_\mu\vDash\phi(\bt)\Rightarrow L_\nu\vDash \phi(\bt)^{L_\mu}$.

(2) Suppose $\bt\in\trace^{\Gm}(L_\nu,p,\alp)\without\bt_0$. Let $L_{\bar{\nu}}\cong L_\nu\{p\cup\bt\cup\{\alp\}\}$ and $L_{\bar{\mu}}\cong L_\mu\{q\cup\bt\cup\{\alp\}\}$, with $\pi$ and $\pi'$ the respective collapsing maps. As $\bt>\bt_0$, we have $L_\mu\{q\cup\bt\cup\{\alp\}\}\sbs L_\nu\{p\cup\bt\cup\{\alp\}\}$ and $\alp\cap L_\mu\{q\cup\bt\cup\{\alp\}\}=\bt=\alp\cap L_\nu\{p\cup\bt\cup\{\alp\}\}$. Thus
$$\pi'\restr\P(\alp)=\pi\restr\P(\alp)\cap L_\mu\{p\cup\bt\cup\{\alp\}\}.$$
 Suppose $\phi(\bt,A)$ for some $\Gamma$ formula $\phi$ and $A\sbs\bt$ with $A\in L_{\bar{\mu}}$. Then by $\Gm$-correctness $L_{\bar{\nu}}\vDash \phi(\bt,A)$ and so $L_\nu\vDash\phi(\alp,\pi^{-1}(A))$. 
 But $\pi^{-1}(A)={\pi'}^{-1}(A)$ by the remark above, so $\pi^{-1}(A)\in L_\mu$. Then as $L_\nu$ and $L_\mu$ are both $\Gm$-correct over $\alp$ we must have $L_\mu\vDash\phi(\alp,\pi^{-1}(A))$ and hence $L_{\bar{\mu}}\vDash \phi(\bt,A)$. By essentially the same argument we have the converse: if $L_{\bar{\mu}}\vDash \phi(\bt,A)$ then $\phi(\bt,A)$. Thus $L_{\bar{\mu}}$ is $\Gm$-correct over $\bt$ so $\bt \in\trace^{\Gm}(L_\mu,p,\alp)$.
 \qed

\begin{lem}\label{coltrace}
If $p$ is a parameter from $L_\mu$ with $\mu>\alp$ and $\pi:L_\mu\{p\cup\bt\cup\{\alp\}\}\cong L_{\bar{\mu}}$ with $\pi(\alp)=\bt$ and $\pi ``p=q$ then $$\mathcal{s}(L_{\bar{\mu}},q,\bt)=\bt\cap\mathcal{s}(L_\mu,p,\alp).$$
The same holds for the $\Gm$ trace etc.
\end{lem}
{\sc Proof:}
Straightforward from the definitions.
\qed

\begin{lem}\label{tracecont}
If $p\in \P(\alp)^{<\om}$ and $\nu>\alp$ is the least limit ordinal such that

(i) $p\in L_\nu$ and 
(ii) $L_\nu$ is $\Gm$-correct over $\alp$, 
then $L_\nu=L_\nu\{p\cup\alp+1\}$.

If in addition $\mathcal{s}^{\Gm}(L_\nu,p,\alp)$ is unbounded in $\alp$
then $L_\nu$ is the union of the filtration, \ie $$L_\nu=\bigcup_{\bt\in\mathcal{s}^{\Gm}(L_\nu,p,\alp)}L_\nu\{p\cup\bt\cup\{\alp\}\}$$
\end{lem}

{\sc Proof:}
Suppose $\pi: L_\nu\{p\cup\alp+1\}\cong L_{\bar{\nu}}$. Then $\pi``p=p$ as $\alp$ is not collapsed. Also we must have $L_{\bar{\nu}}$ is $\Gm$-correct: Let $\phi$ be a $\Gm$ formula with $A\in\P(\alp)\cap L_{\bar{\nu}}$. As $\pi(\alp)=\alp$ and $\pi(A)=A$ we have: 
$$L_{\bar{\nu}}\vDash \phi(\alp,A)\,\,\,\text{ iff  }~L_{{\nu}}\vDash\phi(\alp,A)\,\,\,\text{ iff }~ \phi(\alp,A)$$
So by the minimality of $\nu$, we must have $\nu=\bar{\nu}$.
The second part is an easy consequence.
\qed

The final, and most important lemma in this section gives a close relationship between the traces for specific classes $\Gamma$ and generalised clubs. Just as, at a regular cardinal the trace will form a club, if we are at $\Pi^1_n$-indescribable cardinal then the $\Pi^1_n$-trace forms an $n$-club.  Note that the claim here is not just that the $\Pi^1_n$-trace is in the $n$-club filter, but it is actually $n$-stationary-closed. This is important, as we shall use these  $\Pi^1_n$-traces as our $n$-clubs when we define a $\Bon^n$-sequence in the next section. When in section 3.3 we define the class $\Pi^1_\gm$, we shall extend this lemma to show that the $\Pi^1_\gm$-trace is $\gm$-club.

The class of $\Pi^1_n$ formulae are defined here as the set of formulae of the form $$\forall x_1\sbs\ v_0\exists x_2\sbs v_0\dots {Q} x_{n}\psi(v_0,v_1,x_1,x_2,\dots, x_n)\wedge v_1\sbs v_0$$ where $\psi$ is $\Delta_0$ and $v_0$ is the distinguished variable. Again, we simplify our notation by only allowing 1 parameter (the assignment of $v_1$). We also only quantify over subsets of $\alp$ (rather than $V_\alp$) and as we are working in $L$ 
this does not restrict us. Note that the $\Pi^1_0$-trace is just the trace. This is because the $\Pi^1_0$ formulae have no quantification over subsets of $v_0$, so as the models in the filtration are all transitive below the ordinal for which we require $\Pi^1_0$-correctness and $\Delta_0$ satisfaction is absolute between transitive models, all the models in the filtration are $\Pi^1_0$- correct.

\begin{lem}\label{traceclub}
If $\kp$ is $\Pi^1_{n}$-indescribable then for any limit $\nu>\kp$ with $L_\nu$ $\Pi^1_n$-correct over $\kp$ we have $\mathcal{s}^{\Pi^1_n}(L_\nu,p,\kp)$ is $n$-club in $\kp$.
\end{lem}

{\sc Proof:}
For each $\bt\in\trace(L_\nu,p,\kp)$ set $N_\bt=L_\nu\{p\cup\bt\cup\{\kp\}\}$, and $\pi_\bt :N_\bt\cong L_{\nu_\bt}$. Note $\pi_\bt(\kp)=\bt$. We first show that for $n$-stationary many $\bt$ we have $L_{\nu_\bt}$ is $\Pi^1_n$-correct.

Let $A=\la A_\alp:\alp<\kp, lim(\alp)\ra$ enumerate the subsets of $\kp$ which occur in the filtration such that for any $\bt$ in the trace, $\P(\kp)\cap N_\bt$ is just some initial segment of $A$. As $\kp$ is regular we have on a club $D$ that  $\P(\kp)\cap N_\bt=\{A_\alp:\alp<\bt\}$. Fix some ordering $\la\phi_n(v_0,v_1):n\in\om\ra$ of $\Pi^1_n$ formulae. Then for each limit $\alp<\kp$ we set 
$$C_{\alp+n}=\{\bt<\kp:\phi_n(\bt,A_\alp\cap\bt) \text{ holds}\} \text{ if  }\phi_n(\kp,A_\alp)\text{ holds}$$
and $C_{\alp+n}=\kp$ otherwise.
Now setting $$C=D\cap\triangle_{\alp<\kp}C_\alp$$
we have that $C$ is $n$-stationary as it is the diagonal intersection of elements of the $\Pi^1_n$-indescribability filter, which is normal.
We claim for each $\bt\in C$ that $L_{\nu_\bt}$ is $\Pi^1_n$-correct.  So suppose $\bt\in C$. Then $\bt\in D$ so $\P(\kp)\cup N_\bt=\{A_\alp:\alp<\bt\}$ and $\bt\in \triangle_{\alp<\kp}C_\alp$ so for each $\alp<\bt$ we have if   $\phi_n(\kp,A_\alp)$ then $\phi_n(\bt,A_\alp\cap\bt) $. Now we have $\pi^{-1}:L_{\nu_\bt}\rightarrow L_\nu$ is elementary, $L_\nu$ is $\Pi^1_n$-correct over $\kp$ and $\pi(\kp)=\bt$, and for $\alp<\bt$, $\pi(A_\alp)=A_\alp\cap\bt$. Thus if $\phi$ is $\Pi^1_n$ and $Y\sbs\bt$ with $Y\in L_{\nu_\bt}$ and
$L_{\nu_\bt}\vDash \phi(\bt,Y)$ then this is just $L_{\nu_\bt}\vDash \phi_n(\bt,A_\gm\cap\bt)$ for some $n<\om$ and $\gm<\bt$. Hence we have $L_\nu\vDash  \phi_n(\kp,A_\gm)$ and by $\Pi^1_n$-correctness of $L_\nu$ we have $\phi_n(\kp,A_\gm)$ holds. Then as $\bt\in\triangle_{\alp<\kp}C_\alp$ and $\gm<\bt$ we have $\bt\in C_{\gm+n}$, i.e $\phi_n(\bt,A_\gm\cap\bt)$ holds. 

For the other direction, suppose $L_{\nu_\bt}\vDash\neg\phi_n(\bt,A_\dlt\cap\bt)$. This is $\exists x\sbs\bt\psi(\bt,x,A_\dlt\cap\bt)$ for some $\psi$ which is $\Pi^1_{n-1}$. As a $\Pi^1_{n-1}$ formula is also $\Pi^1_n$ we have by the above that if $L_{\nu_\bt}\vDash\psi(\bt,x,A_\dlt\cap\bt)$ then $\psi(\bt,x,A_\dlt\cap\bt)$ holds, and thus so does $\exists x\sbs\bt\psi(\bt,x,A_\dlt\cap\bt)$, \ie  $\neg\phi_n(\bt,A_\dlt\cap\bt)$. Thus $L_{\nu_\bt}$ is $\Pi^1_n$-correct.
So we have shown that $\mathcal{s}^{\Pi^1_n}(L_\nu,p,\kp)$ is $n$-stationary. (This part of the argument goes through in $V$. We do need $L$ to get $n$-stationary closure though in the next part.)

We now want to show that $\mathcal{s}^{\Pi^1_n}(L_\nu,p,\kp)$ is $n$-stationary closed. So suppose $\bt<\kp$ with $\mathcal{s}^{\Pi^1_n}(L_\nu,p,\kp)$ $n$-stationary below $\bt$. We must have that $\bt$ is regular: regularity of $\alp$ is a $\Pi^1_1$ over $V_\alp$ so as $\kp$ is regular so must be each $\alp\in \mathcal{s}^{\Pi^1_n}(L_\nu,p,\kp)$. Thus, as the regular cardinals do not form a stationary set below any singular, $\bt$ must also be regular. Now as  $\bt$ has an $n$-stationary subset it must be $n-1$-stationary reflecting, and as we are in $L$ and $\bt$ is regular this means $\bt$ is $\Pi^1_{n-1}$-indescribable.  As $\mathcal{s}(L_\nu,p,\kp)$ is unbounded below $\bt$ we have $\bt\in\trace(L_\nu,p,\kp)$, so if $\bt\notin \mathcal{s}^{\Pi^1_n}(L_\nu,p,\kp)$ we must have $L_{\nu_\bt}$ is not $\Pi^1_n$-correct. Thus for some $\Pi^1_n$ sentence $\varphi$ and $X\sbs \kp$ with $\varphi(\kp,X)$ holding, we have $\neg\phi(\bt,X\cap\bt)$. Fixing such $\varphi$ and $X$ we have by the $\Pi^1_{n-1}$ indescribability of $\bt$ that $\{\alp<\bt:\neg\phi(\alp,X\cap\alp)\}$ contains an $n-1$-club. But this cannot be, as $\mathcal{s}^{\Pi^1_n}(L_\nu,p,\kp)$ is $n$-stationary below $\bt$. Thus we must have $\phi(\bt,X)$ and so $L_{\nu_\bt}$ is $\Pi^1_n$-correct and $\bt \in  \mathcal{s}^{\Pi^1_n}(L_\alp,p,\kp)$. So $\mathcal{s}^{\Pi^1_n}(L_\nu,p,\kp)$ is $n$-stationary closed, and hence $n$-club.
\qed


\subsection{The Finite Case}

We can now generalise the construction of $\Bon$-sequences, constructing a $\Bon^n$-sequence below a cardinal which is $\Pi^1_n$- but not $\Pi^1_{n+1}$-indescribable. This theorem will in fact be a corollary of the more general Theorem \ref{squaregm}, but we give the proof of this first, as here we can give the essence of the construction without getting caught up in too many new concepts and tricky details. In this subsection we assume $V=L$ throughout.

\begin{defn}\label{squarendef}
A $\Bon^n$-\emph{sequence on} $\Gm\sbs\kp$ is a sequence $\la C_\alp:\alp\in \Gm\cap d_n(\kp)\ra$ such that for each $\alp$:
\begin{enumerate}
\item $C_\alp$ is an $n$-club subset of $\alp$ and
\item for every $\bt\in d_n(C_\alp)$ we have $\bt\in \Gm$ and $C_\bt=C_\alp\cap\bt$.
\end{enumerate}

We say a $\Bon^n$-sequence $\la C_\alp:\alp\in d_n(\kp)\ra$ \emph{avoids} $E\subset\kp$ if for all $\alp$ we have $E\cap d_n(C_\alp)=\emp$.

We say $S'=\la C'_\alp:\alp\in d_n(\kp)\ra$ is a \emph{refinement} of $S=\la C_\alp:\alp\in d_n(\kp)\ra$ iff for each $\alp$ we have $C'_\alp\sbs C_\alp$.
\end{defn}

\begin{thm*}[\ref{squaren}]$(V=L)$
Let $n\geq 0$ and $\kp$ be a  $\Pi^1_n$- but not $\Pi^1_{n+1}$-indescribable cardinal, and let $A\subseteq\kp$ be $n+1$ stationary.  Then there is $E_{A}\sbs A$ and a $\Bon^n$-sequence $S$ on $\kp$ such that $E_{A}$ is $n+1$-stationary in $\kp$ and $S$ avoids $E_{A}$.\\
Thus $\kp$ is not $n+1$-reflecting.
\end{thm*}
For $n>0$ we produce the $\Bon^n$-sequence in two steps. For the first step, we define $S'=\la C'_\alp:\alp\in Reg\cap d_n(\kp)\ra$ and show that it is a $\Bon^n$-sequence below $\kp$. 
This does not yet depend on the particular $A$.
Then we set
$$E_{A}=\{\alp\in A\cap Reg:C'_\alp\cap A\text{ is not $n$-stationary in $ \alp$ or }d_n(C'_\alp)=\emp\}$$
and in the second step we construct a refinement of $S'$ which avoids $E_{A}$.  If we are just looking for \emph{any} $n+1$-stationary set with a square sequence avoiding it (so we can take $A=\kp$) then this second step is superfluous.
This is because if $A=\kp$ then $E_{A}=\{\alp\in\kp\cap Reg: \alp\notin d_n(\kp)\}$ as $C'_\alp$ is always $n$-stationary by definition, and hence the coherence of $S'$ already guarantees that $S'$ avoids $E_{A}$.
We now fix \cred{$\kp$ and } $A\sbs\kp$ with $A$ $n+1$-stationary. \\

\def\td{\todo{RemoveA?} }
\emph{Constructing $S'$}:

As $\kp$ is $\Pi^1_n$- but not $\Pi^1_{n+1}$-indescribable we can fix $\phi(v_0,v_1,v_2)$ a $\Pi^1_{n}$ formula (for ease of the construction, we use a formula with three free variables)
and $Z\sbs\kp$ such that:
$$\forall X{\sbs}\kp~ \neg\phi(\kp,X,Z)$$
but for all $\alp<\kp$
$$\exists X{\sbs}\alp~ \phi(\alp,X,Z\cap\alp)$$
Let $\alp\in d_n(\kp)$. We set $X_\alp$ to be the $<_{L}$-least subset of $\alp$ such that $ \phi(\alp,X_\alp,Z\cap\alp)$ holds.
We take $\nu_\alp>\alp$ to be the least limit ordinal such that $L_{\nu_\alp}$ is $\Pi^1_n$-correct over $\alp$ and $X_\alp,
Z\cap\alp\in L_{\nu_\alp}$. We set 
$p_\alp=\{X_\alp,
Z\cap\alp\}$ and then we set:
\begin{displaymath}
C'_\alp = \left \lbrace 
\begin{array}{l}
\mathcal{s}^{\Pi^1_n}(L_{\nu_\alp},p_\alp,\alp) \text{\,\,\, if this is $n$-club in }\alp\\
\text{ an arbitrary non-reflecting $n$-stationary set otherwise}\\
\end{array} \right. 
\end{displaymath}
Note that if $\alp$ is $n$-reflecting we shall have that it is $\Pi^1_n$-indescribable and so by Lemma \ref{traceclub} we shall be in the first case. Thus $C'_\alp$ is always well defined.
We set $$S'=\la C'_\alp:\alp\in d_n(\kp)\ra$$ $$E=E_{A}=\{\alp\in A\cap Reg:C'_\alp\cap A\text{ is not $n$-stationary in $ \alp$ or }d_n(C'_\alp)=\emp\}$$ as above.

\begin{claim}
$S'$ is a $\Bon^{n}$-sequence. 
\end{claim}

{\sc Proof:}
It is immediate from the definition that each $C'_\alp$ is $n$-club (in the second case trivially so), so we just need to show that we have the coherence property. So let $\alp\in d_n(\kp)$ and suppose $C'_\alp$ is defined as in the first case (otherwise $d_n(C'_\alp)=\emp$ so coherence is trivial). Let $\bt\in d_n(C'_\alp)$. Let $N_\bt=L_{\nu_\alp}\{p_\alp\cup\bt\cup\{\alp\}\}$ and $\pi:N_\bt\cong L_{\bar{\nu}_\bt}$ be the collapsing map. We need to show $\pi``p_\alp=p_\bt$ and $\nu_\bt=\bar{\nu}_\bt$. 

Clearly  
$\pi(Z\cap\alp)=Z\cap\bt$. Let $X=\pi(X_\alp)=X_\alp\cap\bt$. Then we have 
$$L_{\nu_\alp}\vDash \forall U\sbs\alp \,[U<_L X_\alp\rightarrow\,\neg\varphi(\alp,U,Z\cap\alp)]$$
Then by elementarity we have 
$$L_{\bar{\nu}_\bt}\vDash \forall U\sbs\bt \,[U<_L X\rightarrow\neg\varphi(\bt,U,Z\cap\bt)]$$
and so by absoluteness
$$\forall U\sbs\bt\,[ U<_L X \rightarrow L_{\bar{\nu}_\bt}\vDash \,\,\neg\varphi(\bt,U,Z\cap\bt)].$$
Thus by $\Pi^1_n$-correctness of $L_{\bar{\nu}_\bt}$ we do not have have $X_\bt<_L X$. Also $$L_{\bar{\nu}_\bt}\vDash ~\varphi(\bt,X,Z\cap\alp)$$ and so by $\Pi^1_n$-correctness $X_\bt=X$.

It remains to show that $\nu_\bt=\bar{\nu}_\bt$. We already have that $p_\bt\sbs  L_{\bar{\nu}_\bt}$ and that  $L_{\bar{\nu}_\bt}$ is $\Pi^1_n$-correct over $\bt$, so we only need to show the minimality requirement.
Now for each limit ordinal $\gm>\alp$ with $\gm<\nu_\alp$ we have: 
\begin{displaymath}
\Theta(\gm): \quad(X_\alp\notin L_\gm)
\vee (Z\cap\alp\notin L_\gm) \vee
 \exists U\notin L_\gm\exists n\in\om(U\text{ is minimal with }\psi(n, U))
\end{displaymath}
where $\psi(\ldots)$ is the universal $\Pi^1_{n-1}$ sentence. A little thought will show that this is one way to formalise the requirement that $\nu_\alp$ is minimal.

By $\Pi^1_n$-correctness we have for each $\gm<\nu_\alp$ that $L_{\nu_\alp}\vDash\Theta(\gm)$ and thus $$L_{\nu_\alp}\vDash\forall\gm\Theta(\gm)\Rightarrow  L_{\bar{\nu}_\bt}\vDash\forall\gm\Theta(\gm)\Rightarrow \forall\gm<\bar{\nu}_\bt\,\,\,L_{\bar{\nu}_\bt}\vDash\Theta(\gm)$$
So by $\Pi^1_n$-correctness we have $\Theta(\gm)$ for each limit $\gm<\bar{\nu}_\bt$ , so we do indeed have $\nu_\bt=\bar{\nu}_\bt$.
\qed

We can now define a part of our system $S$. We set $$\Gm^1=\{\alp\in d_n(\kp):C'_\alp\cap A\text{ is $n$-stationary in }\alp\}.$$ If $\alp\in\Gm^1$ we set 
\begin{displaymath}
C_\alp = \left \lbrace 
\begin{array}{l}
d_n(C'_\alp\cap A)\,\,\,\text{ if $d_n(C'_\alp\cap A)$ is $n$-stationary}\\
\text{an arbitrary non-reflecting $n$-stationary set otherwise}\\
\end{array} \right. 
\end{displaymath}

$C_\alp$ is well-defined because as $\alp\in\Gamma^1$ we have $C'_\alp\cap A$ is $n$-stationary so if $\alp$ is $n$-stationary reflecting then $\alp$ is $\Pi^1_n$-indescribable and hence $d_n(C'_\alp\cap A)$ must be $n$-club. Now it is clear that this defines a $\Bon^n$-sequence on $\Gm^1$: each $C_\alp$ is $n$-club by definition and coherence follows from the coherence of $S'$, as if $\bt\in d_n(C_\alp'\cap A)$ then $C_\bt'=C_\alp\cap\bt$ and so $\bt\in \Gm^1$. Also, each such $C_\alp$ avoids $E_{A}$  as either $d_n(C_\alp)=\emp$ or $d_n(C_\alp)\sbs d_n(d_n(C_\alp'\cap A))$.

Now we need to define $C_\alp$ for $\alp\in \Gm^2$, where  $$\Gm^2=\{\alp\in d_n(\kp):C'_\alp\cap A\text{ is not $n$-stationary in }\alp\}=d_n(\kp)\without\Gm^1.$$ For such $\alp$ we shall find $C_\alp\sbs C'_\alp$ such that  (i) $C_\alp$ avoids $A$ and hence $E$ and (ii) for $\bt\in d_n(C_\alp)$ we have $C'_\bt\cap A$ is not $n$-stationary (\ie $\bt\in\Gm^2$) and  $C_\bt=C_\alp\cap\bt$.
Once we have (i) and (ii) it is easy to see that $S=\la C_\alp:\alp\in d_n(\kp)\ra$ will satisfy Theorem \ref{squaren}, and it will only remain to show that $E_{A}$ is $n+1$-stationary.

For $\alp\in\Gm^2$ set $p'_\alp=\{C'_\alp, A\cap\alp\}\cup p_\alp$,  and  take $\eta_\alp\geq\nu_\alp$ to be the minimal limit ordinal such that $C'_\alp\in L_{\eta_\alp}$ and $L_{\eta_\alp}$ is $\Pi^1_n$-correct over $\alp$. Note that if $D_\alp$ is the $<_L$ least $n-1$-club avoiding $A\cap C'_\alp$ then $D_\alp\in L_{\eta_\alp}$ by $\Pi^1_{n}$-correctness. We now set 
\begin{displaymath}
C_\alp = \left \lbrace 
\begin{array}{l}
D_\alp\cap\mathcal{s}^{\Pi^1_n}(L_{\eta_\alp}, 
p'_\alp,\alp) \text{\,\,\, if this is $n$-club in }\alp\\
\text{ an arbitrary non-reflecting $n$-stationary subset of $C'_\alp$ otherwise.}\\
\end{array} \right. 
\end{displaymath}
Then $C_\alp$ is well defined because if $\alp$ reflects $n$-stationary sets then $\alp$ is $\Pi^1_n$-indescribable. Then the $\Pi^1_n$-trace of $L_{\nu_\alp}$ is $n$-club, and hence so also is its intersection with $D_\alp$.

If $C_\alp$ is defined as in the first case, then $C_\alp\sbs d_n(C'_\alp)$ because for each $\bt\in \trace^{\Pi^1_n}(L_{\eta_\alp},
p'_\alp,\alp)$ we have $C'\cap\bt$ is $n$-stationary. Thus $C_\alp\sbs C'_\alp$. Also because $C_\alp\sbs D_\alp$ or $d_n(C_\alp)=\emp$ we have (i): $C_\alp$ avoids $A$. If $\bt\in d_n(C_\alp)$ then $C_\alp$ is defined as in the first case so $\bt\in d_n(C'_\alp)$ and hence $C'_\bt=C'_\alp\cap\bt$. Then we have $A\cap C'_\bt$ is not $n$-stationary in $\bt$ by elementarity and $\Pi^1_n$-correctness, so $\bt\in\Gm^2$. Also by elementarity and $\Pi^1_n$-correctness, $D_\bt=D_\alp\cap\bt$ and so by Lemma \ref{coltrace} 
$$\mathcal{s}^{\Pi^1_n}(L_{\eta_\bt},
p'_\bt,\bt)=\mathcal{s}^{\Pi^1_n}(L_{\eta_\alp},
p'_\alp,\alp)\cap\alp.$$ Thus $$C_\bt=D_\bt\cap\mathcal{s}(L_{\eta_\bt},
p'_\bt,\bt)=\bt\cap D_\alp\cap\mathcal{s}(L_{\eta_\alp},
p'_\alp,\alp)=C_\alp\cap\bt.$$ So we have (ii). This completes our construction of  the $\Bon^n$-sequence on the regulars avoiding $E_A$.

This defines the square sequence for regular $\alp\in d_n(\kp)$; for singular $\alp \in d_n(\kp)$ we argue as follows: as $E_A$ has been defined, it is a set of regular cardinals, and so cannot even be stationary below a singular ordinal. Thus we can simply use Jensen's global square sequence below singulars -  let $\la D_\alp: \alp\in Sing\ra$ be a global square sequence. Then for any singular $\alp$, we have $d_0(D_\alp)\sbs Sing$ and thus $D_\alp$ avoids $E_A$, and of course we have coherence.

It remains to show that $E_{A}$ is $n+1$-stationary. 
\begin{defn}
We define $H\sbs A$ by letting $\alp\in H$ iff $\alp\in A$ and there is $\mu_\alp >\alp$ and $q$  a parameter from $L_{\mu_\alp}$ such that:
\begin{enumerate}
\item  $L_{\mu_\alp}$ is $\Pi^1_n$-correct over $\alp$;
\item $\mu_\alp<\nu_\alp$;
\item $A\cap\mathcal{s}^{\Pi^1_n}(L_{\mu_\alp},q,\alp)=\emp$.
\end{enumerate} 
\end{defn}

\begin{lem}
$H\sbs E_{A}$.
\end{lem}
{\sc Proof:}
Suppose $\alp\notin E_{A}$, so we have $C'_\alp\cap A$ is $n$-stationary and $d_n(C_\alp')\neq\emp$. Thus $C'_\alp$ was defined as in the first case: $C'_\alp=\mathcal{s}(L_{\nu_\alp},p_\alp,\alp)$. Let $\mu<\nu_\alp$ and $q\in L_\mu$. Then using lemma \ref{tracemono} for some $\bt<\alp$ we have $q\in L_{\nu_\alp}\{\bt\cup p_\alp\cup\{\alp\}\}$ and so $\mathcal{s}(L_\mu,q,\alp)\without\bt\supseteq\mathcal{s}(L_{\nu_\alp},p_\alp,\alp)=C'_\alp$. Thus $\alp\notin H$.
\qed
\begin{lem}
$H$ is $n+1$-stationary.
\end{lem}
{\sc Proof:}
Let $C\sbs\kp$ be $n$-club. Take $\mu>\kp$ minimal such that $C,Z\in L_\mu$ and $L_\mu$ is $\Pi^1_n$-correct. Set $D=\mathcal{s}(L_\mu,\{C,Z\},\kp)$ and note that $D\sbs d_n(C)\sbs C$, and $D$ is $n$-club by Lemma \ref{traceclub}. Take $\dlt=min (D\cap A)$. Set $\mu_\dlt$ to be the ordinal such that  $L_\mu\{\{C,Z,\kp\}\cup\dlt\}\cong L_{\mu_\dlt}$.
 We show that $\dlt\in H$, with $\mu_\dlt$ and $\{C\cap\dlt,Z\cap\dlt\}$ witnessing this. 

{\em Claim 1:} $\mu_\dlt<\nu_\dlt$.\\
We have that for every $X\in \P(\dlt)\cap L_{\mu_\dlt}$ 
$$L_{\mu_\dlt}\vDash \,\, \neg\phi(\dlt,X,Z\cap\dlt).$$
As $L_{\mu_\dlt}$ is $\Pi^1_n$-correct this means: $$ \forall X\in \P(\dlt)\cap L_{\mu_\dlt}\,\,  \neg\phi(\dlt,X,Z\cap\dlt).$$ But we know by the choice of $\nu_\dlt$ that 
 $$ \exists X\in \P(\dlt)\cap L_{\nu_\dlt}\,\,  \phi(\dlt,X,Z\cap\dlt).$$
 Hence we must have $\nu_\dlt>\mu_\dlt$.
 
{\em Claim 2:} $A\cap\mathcal{s}(L_{\mu_\dlt},\{C\cap\dlt, Z\cap\dlt\},\dlt)=\emp$.\\
We have (Lemma \ref{coltrace}) $A\cap\mathcal{s}(L_{\mu_\dlt},\{C\cap\dlt, Z\cap\dlt\},\dlt)=A\cap\dlt\cap\mathcal{s}(L_{\mu},\{C, Z\},\kp)=A\cap\dlt\cap D=\emp$.
\qedtwo{ Theorem \ref{squaren}}




\subsection{$Q_\gm$ games, $\Pi^1_\gm$-indescribability and $\Bon^{<\gm}$}
In this section we introduce the definitions we shall need to generalise Theorem \ref{squaren} of the previous section: $\Pi^1_\gm$ and $\Bon ^{<\gm}$.

\subsubsection{$\Pi^1_\gm$-indescribability}

We now turn to a generalisation of the notions $\Pi^1_n$ and $\Sigma^1_n$, introduced in \cite{Sharpe-Welch2012}.\footnote{ 
In \cite{B2TAMS} Bagaria introduces a definition different to the one we present and show that in $L$ a cardinal is $\gm$-stationary reflecting iff $\Pi^1_\gm$-indescribable in his sense. We believe the proof of the following section constructing a $\Bon^{<\gm}$-sequence would also work with his definition.} We shall use these notions to strengthen Theorem \ref{squaren}, obtaining $\Bon^{<\gm}$-sequences below cardinals which are $\Sigma^1_\gm$- but not $\Pi^1_\gm$-indescribable. 

The following definitions are (Definitions 3.15, 3.16 and 3.21), with a slightly different numbering and presentation of the odd cases. 
\begin{defn}
The $Q^1_\alp$ {\em  game on} $\kappa$ is a two player game of finite length. With parameter $A\sbs\kp$ and a $\Delta_0$ sentence $\varphi$ with 3 free variables, the game $G_\alp(\kp,\varphi,A)$ is defined as follows. \\

\textbf{Case 1}: $\alp$ is an even ordinal (including the case $\alp$ is a limit).

In round $n$ ($n\geq1$) Player $\Sigma$ plays a pair, $(\alp_n,X_n)$ and player $\Pi$ follows,  playing $Y_n$, with the following constraints:
\begin{enumerate}
\item Each $\alp_n$ is an odd ordinal and setting $\alp_0=\alp$ we have $\alp_n<\alp_{n-1}$,
\item $X_n, Y_n\sbs\kp$
\item Player $\Pi$ must play $Y_n$ such that $\varphi(\vec{X}_n,\vec{Y}_n,A)$ holds, where $\vec{X}_n= \la X_1,\dots, X_n\ra$ and $\vec{Y}_n= \la Y_1,\dots, Y_n\ra$
\end{enumerate}
The first player to be unable to move loses.\\

\textbf{Case 2}: $\alp$ is an odd ordinal.

The game here is similar: the players switch roles but $\Sigma$ still starts.
In round $n$  for $n\geq 1$, $\Sigma$ plays $Y_{n}$ such that $\varphi(\vec{X}_{n-1},\vec{Y}_{n},A)$, (setting $\vec{X}_0=\emp$)
and then $\Pi$ plays a pair, $(\alp_n,X_n)$ with the same constraints as for $\Sigma$ above, \ie a decreasing sequence of odd ordinals.

\end{defn}

The decreasing sequence of ordinals ensures the games are always finite in length, so by the familiar Gale-Stewart argument they are determined - one player always has a winning strategy. So without ambiguity, we shall say $\Sigma$ wins $G_\alp$ to mean $\Sigma$ has a winning strategy for $G_\alp$. In the sequel we shall write as if for a pair $A,Z\sbs\kp$, $\la A,Z\ra$ were also a subset of $\kp$ \via G\"odel pairing.

\begin{rem}
It is easy to see that $\Sigma$ wins $G_{\alp+1}(\kp,\phi,A)$ iff there exists $Z\sbs\kp$ such that $\Pi$ wins $G_{\alp}(\kp,\phi',\la A,Z\ra)$ where
\begin{displaymath}
\begin{split}
\phi'(\vec{X},\vec{Y},\la A,Z\ra)\leftrightarrow \phi(Z^\frown \vec{X},\vec{Y},A) & \text{\,\, if $\alp+1$ is even and }\\
\phi'(\vec{X},\vec{Y},\la A,Z\ra)\leftrightarrow \phi(\vec{X},Z^\frown\vec{Y},A) & \text{ \,\, if $\alp+1$ is odd.}
\end{split}
\end{displaymath}
Similarly $\Pi$ wins $G_{\alp+1}(\kp,\phi,A)$ iff for all $X\sbs\kp$ we have that $\Sigma$ wins $G_{\alp}(\kp,\phi',\la A,X\ra)$. Thus we have that $``\Sigma$ wins $G_n(\kp,\phi,A)"$ is a $\Sigma^1_n$ statement about $V_\kp$ and a statement about who wins $G_{\alp+n}$ is equivalent to a statement with $n$ alternations of second order quantifiers followed by a statement about who wins $G_{\alp}$. 

For the limit levels, a statement that $\Sigma$ wins $G_{\alp}(\kp,\phi,A)$ is equivalent to the statement that for some odd $\alp'<\alp$ (\ie the $\alp_1$ chosen by $\Sigma$) and $Y\sbs\kp$ we have that $\Pi$ wins $G_{\alp'}(\kp,\phi',\la A,Y\ra)$, and the converse: $\Pi$ wins $G_{\alp}(\kp,\phi,A)$ iff for every odd $\alp'<\alp$ and $Y\sbs\kp$ we have that $\Sigma$ wins $G_{\alp'}(\kp,\phi',\la A,Y\ra)$. Note that $\phi'$ is the same formula in all these instances - it depends only on $\phi$ and whether we are in an even or odd game. 

If $\alp$ is even and $\Sigma$ wins $G_{\alp}(\kp,\phi,A)$ with first move $\la\alp_1,X\ra$ then $\Sigma$ also wins $G_{\alp_1+1}(\kp,\phi,A)$ as she can begin this game with the same move $\la\alp_1,X\ra$.
\end{rem}

This leads to the following definitions:

\begin{defn}
A formula is $\Pi^1_\gm$ (respectively $\Sigma^1_\gm$) {\em over} $V_\kp$ if it is of the form ``$\Pi$ (resp. $\Sigma$) wins the game $G_\gm(\kp,\phi,A)$" for some $\Delta_0$ formula $\phi$ with 3 free variables and $A\sbs\kp$.
\end{defn}

We could define a notion of $\Pi^1_\gm$ set here as the class of subsets of $V_\kp$ of the form $\{x:\Pi$ wins the game $G_\gm(\kp,\phi,\la A,x\ra)\}$. However we do not use it here.

The following show that this is a good candidate for a generalisation of $\Pi^1_n$ for our purposes. 

\begin{prop}\label{good}
Being $\gm$-stationary is expressible in a $\Pi^1_\gm$ way.
\end{prop}

{\sc Proof:}
We give the formula $\phi$ for the case $\gm$ is even, the other case is only superficially different. We in fact show that generalised stationarity is \emph{uniformly} expressible in the sense that there is a $\Delta_0$ formula $\Theta$ with three free variables such that for any $\kp$ and any (even) $\gm<\kp$ and $A\sbs\kp$, $\Pi$ wins $G_\gm(\kp,\Theta,A)$ iff $A$ is $\gm$-stationary in $\kp$. 
we shall have that in each round $i$, $\Sigma$ plays $\la\alp_i,X_i\ra$ with $X_i=\la\bt_i, S_i^1,S_i^2\ra$ and $\Pi$ plays $Y_i=\la\dlt_i, T_i^1,T_i^2\ra$. In round 1, $\Sigma$ tries to show that $A$ is not $\gm$-stationary by choosing $\bt_1<\gm$ and $\bt$-stationary $S_1^1,S_1^2\sbs\kp$ such that $d_{\bt_1}(S_1^1)\cap d_{\bt_1}(S_1^2)\cap A=\emp$. Then $\Pi$ tries to show that $\Sigma$ fails, by choosing $\dlt_1<\bt_1$ and  $T_1^1,T_1^2$ to witness that either $S_1^1$ or $S^2_1$ is not $\bt_1$-stationary. Then $\Sigma$ tries to show that either  $T_1^1$ or $T_1^2$ is not $\dlt_1$-stationary, and so on. The game will end when either player chooses $0$. A slight adjustment has to be made as we don't want to have $\gm$ as a parameter in $\Theta$.

\begin{displaymath}
\begin{split}
\exists i  \leq lh(\vec{X})~\neg\Big[ &\big(X_i=\la\bt_i, S_i^1,S_i^2\ra\text{ with } S_i^1,S_i^2\text{ unbounded}\big)\\
&\wedge \big(d_{\bt_1}(S_1^1)\cap d_{\bt_1}(S_1^2)\cap A=\emp\big) \wedge \big(\bt_{i+1}<\dlt_i\big) \\
&\hspace{-6.5em} \wedge \big(d_{\bt_{i+1}}(S_{i+1}^1)\cap d_{\bt_{i+1}}(S_{i+1}^2)\cap T_i^1=\emp\vee d_{\bt_{i+1}}(S_{i+1}^1)\cap d_{\bt_{i+1}}(S_{i+1}^2)\cap T_i^2=\emp\big)\Big]\\
 \vee ~ \forall i \leq lh(\vec{X})\Big[& \big(Y_i=\la\dlt_i, T_i^1,T_i^2\ra \text{ with } S_i^1,S_i^2\text{ unbounded}\big)\wedge (\dlt_{i}<\bt_i)\\
& \wedge \big(d_{\dlt_{i}}(T_{i}^1)\cap d_{\dlt_{i}}(T_{i}^2)\cap S_i^1=\emp\vee d_{\dlt_{i}}(T_{i}^1)\cap d_{\dlt_{i}}(T_{i}^2)\cap S_i^2=\emp\big)\Big]\\
\vee ~\forall  i<lh(\vec{X}) &\Big[\bt_i\neq 0\wedge\dlt_i\neq 0\wedge\bt_{i+1}\neq 0\Big].\\
\end{split}
\end{displaymath}

The final disjunct ensures that $\Sigma$ cannot ``cheat" by choosing $\bt_1\geq\gm$ - \ie for a $\gm$-stationary $A$, $\Pi$ has a winning strategy in $G_\gm(\kp,\Theta,A)$ even if $\bt_1\geq\gm$. Recall $\Sigma$ must choose odd ordinals $\alp_{i+1}<\alp_i$ with $\alp_1<\gm$. In the first round, $\Pi$ chooses $Y_1=\la\alp_1,A,A\ra$ if $\bt_1\geq\gm$, noting $\alp_1<\gm$ so $\alp_1<\bt_1$. The game then proceeds as if starting from round 2, but this time $\Sigma$ cannot ``cheat". The final disjunct means that it doesn't matter that $\Pi$'s first move doesn't satisfy the second disjunct, unless $\Sigma$ gets down to choosing $\bt_i=0$ and $S_i^1,S_i^2${ with }$ S_i^1,S_i^2\text{ unbounded and } d_{0}(S_{i}^1)\cap d_{0}(S_{i}^2)\cap T_{i-1}^1=\emp$ or $d_{0}(S_{i}^1)\cap d_{0}(S_{i}^2)\cap T_{i-1}^2=\emp$ - but then $\Sigma$ could have won fair and square, without starting with $\bt_1\geq\gm$. 
\qed

\begin{defn}
A cardinal $\kp$ is $\Pi^1_\gm${\em -indescribable} (respectively $\Sigma^1_\gm${\em -indescribable}) if for every $\Delta_0$ formula $\phi$ with 3 free variables and every parameter $A\sbs\kp$ such that $\Pi$ (resp. $\Sigma$) wins the game $G_\gm(\kp,\phi, A)$ we have that there is some $\alp<\kp$ such that $\Pi$ (resp. $\Sigma$) also wins the game $G_\gm(\alp,\phi, A\cap\alp)$.
\end{defn}

\begin{rem}
Note that as with the finite levels if a cardinal is $\Pi^1_\gm$-indescribable then it is also $\Sigma^1_{\gm+1}$-indescribable, and indeed if it is $\Pi^1_\gm$-indescribable for all $\gm<\lmb$ then it is $\Sigma^1_\lmb$-indescribable.
\end{rem}

A consequence of Proposition \ref{good} is that any $\Pi^1_\gm$-indescribable cardinal is $\gm$-reflecting, and any $\Sigma^1_\gm$-indescribable is $\eta$-reflecting for every $\eta<\gm$.

\begin{lem}\label{univ}
For each $\gm$, we have a universal $\Pi^1_\gm$ sentence.
\end{lem}
{\sc Proof:}
This is straightforward as we only really need a universal $\Delta^1_0$ formula. Fix a recursive enumeration $\la\phi_n:n\in\om\ra$ of $\Delta_0$ formulae with 3 free variables. Let $\Psi$ be the $\Delta_0$ formula such that:
$$\Psi(\alp,\vec{X}_m,\vec{Y_m},\la A,n\ra)\leftrightarrow \phi_n(\alp,\vec{X}_m,\vec{Y}_m,A)$$
Then $\Pi$ wins $G_\gm(\alp,\Psi,\la A,n\ra)$ iff $\Pi$ wins $G_\gm(\alp,\phi_n,A)$, so ``$\Pi$ wins $G_\gm(\alp,\Psi,X)$" is universal.
\qed

Thus being $\Pi^1_\gm$-indescribable is $\Pi^1_{\gm+1}$-expressable, and being $\Sigma^1_\gm$-indescribable is $\Pi^1_\gm$-expressable.
The $\Pi^1_\gm$-indescribability filter is defined in the obvious way:
\begin{defn}\label{gmindescribabilityfilter}
The $\Pi^1_\gm$ {\em indescribability filter} on $\kp$, $\F^\gm_\kp$ is the filter generated by sets of the form:
$$\{\alp<\kp:\Pi\text{ wins } G_\gm(\kp,\phi,A) \text{ and }G_\gm(\alp,\phi,A\cap\alp)\}$$
\end{defn}

\begin{rem}
As being $\gm$-stationary is $\Pi^1_\gm$ expressible, and $\Pi^1_\gm$ is closed under conjunction, each $F\in\F^\kp_\gm$ is $\gm+1$-stationary.
\end{rem}

The argument is essentially that of Levy for finite $n$. We give the proof here for the sake of completeness.

\begin{lem}\label{normalgamma}
If $\kp$ is $\Pi^1_\gm$-indescribable, then the $\Pi^1_\gm$-indescribable filter is normal and $\kp$-complete. 
\end{lem}

{\sc Proof:}
Fix an ordinal $\gm$ and $\kp>\gm$ with $\kp$ $\Pi^1_\gm$-indescribable, and let $X\sbs\kp$ be positive with respect to $\F^\gm_\kp$ and $f:X\rightarrow\kp$ be regressive. Suppose for a contradiction that for each $\bt<\kp$ there is an element of $\F^\gm_\kp$ that avoids $f^{-1}``\{\bt\}$. For each $\bt<\kp$ choose $n_\bt$ and $A_\bt\sbs\kp$ witnessing this, \ie such that $\Pi$ wins $G_\gm(\kp,\phi_{n_\bt},A_\bt)$ but if $f(\alp)=\bt$ then $\Sigma$ wins $G_\gm(\alp,\phi_{n_\bt},A_\bt)$. Thus for each $\alp\in X$ we have 
\begin{equation}
\Sigma\text{ wins }G_\gm(\alp,\phi_{n_{f(\alp)}},A_{f(\alp)})\tag{$*$}
\end{equation}
\textbf{Claim:} $``\forall\bt<\kp~ \Pi\text{ wins } G_\gm(\bt,\phi_{n_\bt},A_\bt)"$ is $\Pi^1_\gm$. 

Set $A=\{\la\bt,0,\dlt\ra:\dlt\in A_\bt\}\cup\{\la\bt,1,n_\bt\ra\}$, so that $A\sbs\kp$ codes $A_\bt$ and $n_\bt$ for each $\bt\in\kp$. 
The following $\Delta_0$ sentence $\Phi$ with parameter $A$ produces an appropriate game: $X_0$ is / is not a pair $X_0=\la \bt,X_0'\ra$ with $\bt\in\kp,X_0'\sbs\kp$ and/or $\Psi(\bt,{X_0'}^\frown\vec{X_m}, \vec{Y}_m,\la n_{\bt},A_\bt\ra)$, with ``is \dots and" if $\gm$ is odd and ``is not \dots or" if $\gm$ is even, and $\psi$ is universal as in Lemma \ref{univ}.

Now we have $\Pi$ wins $G_\gm(\kp,\Phi,A)$ so by the $\Pi^1_\gm$-indescribability of $\kp$, for some $\alp\in X$ we have $\Pi$ wins $G_\gm(\alp,\Phi,A\cap\alp)$. Thus, by definition of $\Phi$ and the fact that $f$ is regressive, $\Pi$ wins $G_\gm(\alp,\phi_{n_{f(\alp)}},A_{f(\alp)})$. But this contradicts $(*)$.

Now we have introduced these notions, we apply them to prove the analogue of Lemma \ref{traceclub} for $\Pi^1_\gm$. The proof is essentially the same as before, except that in the proof of Lemma \ref{traceclub} we used the fact proven in \cite{BaMS} that in $L$, any regular cardinal which reflects $n$-stationary sets, and hence any which admits $n+1$-stationary sets, is $\Pi^1_n$-indescribable. As we have not yet shown this for $\gm$-stationary sets and $\Pi^1_\gm$-indescribability, this lemma must be proven inductively along with Theorem \ref{squaregm}.

\begin{lem}\label{traceclub2}$(V=L)$
Let $\gm$ be an ordinal and assume that for all $\gm'<\gm$ we have that any regular cardinal which is $\gm'$-reflecting is $\Pi^1_{\gm'}$-indescribable. If $\kp>\gm$ is $\Pi^1_{\gm}$-indescribable then for any limit $\nu>\kp$ with $L_\nu$ $\Pi^1_\gm$-correct over $\kp$ we have $\mathcal{s}^{\Pi^1_\gm}(L_\nu,p,\kp)$ is $\gm$-club in $\kp$.
\end{lem}

{\sc Proof:}
We proceed as in \ref{traceclub}. For each $\bt\in\trace(L_\nu,p,\kp)$ set $N_\bt=L_\nu\{p\cup\bt\cup\{\kp\}\}$, and $\pi_\bt :N_\bt\cong L_{\nu_\bt}$. Note $\pi_\bt(\kp)=\bt$. We first show that for $\gm$-stationary many $\bt$ we have $L_{\nu_\bt}$ is $\Pi^1_\gm$-correct.

Let $A=\la A_\alp:\alp<\kp\wedge lim(\alp)\ra$ enumerate the subsets of $\kp$ which occur in the filtration such that for any $\bt$ in the trace, $\P(\kp)\cap N_\bt$ is just some initial segment of $A$. As $\kp$ is regular we have on a club $D$ that  $\P(\kp)\cap N_\bt=\{A_\alp:\alp<\bt\}$. Fix some ordering $\la\phi_n(v_0,v_1,v_2):n\in\om\ra$ of $\Delta_0$ formulae with all free variables displayed. Then for each limit $\alp<\kp$ we set 
\begin{displaymath}
C_{\alp+n} = \left \lbrace 
\begin{array}{l}
\{\bt<\kp: \Pi\text{ wins }G_\gm(\bt,\phi_n,A_\alp\cap\bt) \} \,\,\,\,\,\, \text{ if  } \Pi \text{ wins }G_\gm(\bt, \phi_n,A_\alp)\\
\text{ $\kp$ otherwise.}\\
\end{array} \right. 
\end{displaymath}

Now setting $$C=D\cap\triangle_{\alp<\kp}C_\alp$$
we have that $C$ is $\gm$-stationary as it is the diagonal intersection of elements of the $\Pi^1_\gm$-indescribability filter, which is normal by Lemma \ref{normalgamma}.
We claim for each $\bt\in C$ that $L_{\nu_\bt}$ is $\Pi^1_\gm$-correct. So suppose $\bt\in C$. Then $\bt\in D$ so $\P(\kp)\cup N_\bt=\{A_\alp:\alp<\bt\}$ and $\bt\in \triangle_{\alp<\kp}C_\alp$ so for each $\alp<\bt$ we have if   $\Pi\text{ wins }G_\gm(\kp,\phi_n,A_\alp) $ then $\Pi\text{ wins }G_\gm(\bt,\phi_n,A_\alp\cap\bt) $. Now we have $\pi^{-1}:L_{\nu_\bt}\rightarrow L_\nu$ is elementary, $L_\nu$ is $\Pi^1_\gm$-correct over $\kp$ and $\pi(\kp)=\bt$,  $\pi(A_\alp)=A_\alp\cap\bt$ for $\alp<\bt$. Thus if $\psi$ is $\Pi^1_\gm$ then
$L_{\nu_\bt}\vDash \psi(\bt,X)$ is just $L_{\nu_\bt}\vDash \Pi\text{ wins }G_\gm(\bt,\phi_n,A_\dlt\cap\bt) $ for some $n<\om$ and $\dlt<\bt$. Hence we have $L_\nu\vDash \Pi\text{ wins }G_\gm(\kp, \phi_n, A_\dlt)$ and by $\Pi^1_\gm$-correctness $\Pi\text{ does indeed win }G_\gm(\kp,  \phi_n, A_\dlt)$. Then as $\bt\in\triangle_{\alp<\kp}C_\alp$ and $\dlt<\bt$ we have $\bt\in C_{\dlt+n}$, \ie $\Pi\text{ wins }G_\gm(\bt,\phi_n,A_\dlt\cap\bt) $.

 As for $\gm'<\gm$, $\Pi^1_{\gm'}$ sentences are also $\Pi^1_\gm$, we have the other direction by induction. 
Thus $L_{\nu_\bt}$ is $\Pi^1_\gm$-correct.
So we have shown that $\mathcal{s}^{\Pi^1_n}(L_\nu,p,\kp)$ is $\gm$-stationary. (Again, this part of the argument goes through in $V$, but we need $L$ to get $\gm$-stationary closure.)

We show that $\mathcal{s}^{\Pi^1_n}(L_\nu,p,\kp)$ is $\gm$-stationary closed. Suppose $\bt<\kp$ with $\mathcal{s}^{\Pi^1_\gm}(L_\nu,p,\kp)$ $\gm$-stationary below $\bt$. As in \ref{traceclub}, $\bt$ must be regular. Now as $\bt$ has a $\gm$-stationary subset it must be $\gm'$-stationary reflecting for every $\gm'<\gm$, and as we are in $L$ and $\bt$ is regular by our hypothesis we have $\bt$ is $\Sigma^1_{\gm}$-indescribable.  As $\mathcal{s}(L_\nu,p,\kp)$ is unbounded below $\bt$ we have $\bt\in\trace(L_\nu,p,\kp)$, so if $\bt\notin \mathcal{s}^{\Pi^1_\gm}(L_\nu,p,\kp)$ we must have $L_{\nu_\bt}$ is not $\Pi^1_\gm$-correct. Thus for some $\Delta_0$ formula $\varphi$ with $\Pi\text{ wins }G_\gm(\kp,\phi,X)$ we have $\Sigma$ wins $G_\gm(\bt,\phi,X\cap\bt)$. This means by the remark above that for some $Y\sbs\bt$ and $\gm'<\gm$ we have $\Pi$ wins $G_{\gm'}(\bt,\phi', \la X\cap\bt,Y\ra)$. By the $\Pi^1_{\gm'}$-indescribability of $\bt$, $\{\alp<\bt:\Pi$ wins $G_{\gm'}(\alp,\phi',\la X\cap\alp,Y\cap\alp\ra)\}$ contains an $\gm'$-club. But if $\Pi$ wins $G_{\gm'}(\alp,\phi',\la X\cap\alp,Y\cap\alp\ra)$ then $\Sigma$ wins $G_\gm(\alp,\phi,X\cap\alp)$ and so $\alp\notin \mathcal{s}^{\Pi^1_\gm}(L_\nu,p,\kp)$.

 But this cannot be, as $\mathcal{s}^{\Pi^1_\gm}(L_\nu,p,\kp)$ is $\gm$-stationary below $\bt$. Thus we must have $\Pi\text{ wins }G_\gm(\bt,\phi,X\cap\bt)$ and so $L_{\nu_\bt}$ is $\Pi^1_\gm$-correct\footnote{The other direction follows by upwards absoluteness of $\Sigma_1$ formulae.} and $\bt \in  \mathcal{s}^{\Pi^1_\gm}(L_\alp,p,\kp)$. So $\mathcal{s}^{\Pi^1_\gm}(L_\nu,p,\kp)$ is $\gm$-stationary closed, and hence $\gm$-club.
\qed


\subsubsection{$\Bon^{\gm}$ and $\Bon^{<\gm}$}

Definition \ref{squarendef} (defining a $\Bon^n$-sequence) can easily be stated for general ordinals $\gm$. However, this is not the type of $\Bon$-sequence we shall need to take Theorem \ref{squaren} into the transfinite. The problem here is the limit levels. To continue inductively along the ordinals, we first need to show that, if a regular cardinal is $\Sigma^1_\om$-indescribable (\ie $\Pi^1_n$-indescribable for every $n$), but not $\Pi^1_\om$-indescribable, then it does not reflect $\om$-stationary sets, and similarly for further limit ordinals. To witness this, we need to define a new type of $\Bon$-sequence.

\begin{defn}
A $\Bon^{<\gm}$ \emph{sequence} on $\Gm\sbs LimOrd\cap\kp$ is a sequence $\la (C_\alp, \eta_\alp):\alp\in \Gm\ra$ such that for each $\alp$:
\begin{enumerate}
\item $\eta_\alp<\gm$ and $C_\alp$ is an $\eta_\alp$-club subset of $\alp$
\item for every $\bt\in d_{\eta_\alp}(C_\alp)$ we have $\bt\in \Gm$ with $\eta_\alp=\eta_\bt$ and $C_\bt=C_\alp\cap\bt$
\end{enumerate}

We say a $\Bon^{<\gm}$-sequence $\la C_\alp:\alp\in \kp\ra$ \emph{avoids} $A\subset\kp$ if for all $\alp\in\kp$ we have $A\cap d_{\eta_\alp}(C_\alp)=\emp$.

We say $S'=\la (C'_\alp, \eta'_\alp):\alp\in \Gamma\ra$ is a \emph{refinement} of $S=\la (C_\alp,\eta_\alp):\alp\in \Gamma\ra$ iff for each $\alp$ we have $\eta'_\alp=\eta_\alp$ and $C'_\alp\sbs C_\alp$.
\end{defn}

Equipped with this definition, we can see that a $\Bon^{<\gm}$-sequence avoiding some $\gm$-stationary set $A$ witnesses that $A$ does not reflect (\ie $d_\gm(A)=\emp$), even in the case $\gm$ is a limit. This is because, for $\alp\in\kp$, $C_\alp$ is an $\eta_\alp$-club avoiding $A\cap\alp$ and thus $A$ is not $\eta_\alp+1$-stationary in $\alp$, and as $\eta_\alp+1\leq\gm$ we have $\alp\notin d_\gm(A)$. As this works for both limit and successor $\gm$, we shall use only $\Bon^{<\gm}$ in our proof of Theorem \ref{squaregm} and thus deal with the limit and successor cases simultaneously. The analogue of Theorem \ref{squaren} for the infinite ordinals is a corollary: restricting the domain of a $\Bon^{<\gm+1}$-sequence to the ordinals in $d_\gm(\kp)$ gives a $\Bon^\gm$-sequence.


\subsection{The Main Result: $\Bon^{<\gm}$ at a Non-$\Pi^1_\gm$-Indescribable}

We are now ready to state and prove the main theorem of this chapter:

\begin{thm}$(V=L)$\label{squaregm}
Let $\gm<\kp$ be an ordinal and $\kp$ be $\Sigma^1_\gm$-indescribable (\ie $\Pi^1_\eta$-indescribable for every $\eta<\gm$) but not $\Pi^1_{\gm}$-indescribable, and let $A\subseteq\kp$ be $\gm$-stationary. Then there is $E_{A}\sbs A$ and a $\Bon^{<\gm}$-sequence $S$ on $\kp$ such that $E_{A}$ is $\gm$-stationary in $\kp$ and $S$ avoids $E_{A}$.
Thus $\kp$ is not $\gm$-reflecting.
\end{thm}

\begin{col}\label{Cor3.25}
$(V=L)$ A regular cardinal is $\Pi^1_\gm$-indescribable iff it reflects $\gm$-stationary sets.
\end{col}
{\sc Proof of Theorem:}
We assume $ V=L$.
The proof is an induction on $\gm$ making repeated use of Lemma \ref{traceclub2}. Thus we fix $\gm$ and assume we have that for any $\gm'<\gm$, if $\alp$ is $\gm'$-reflecting and regular then $\alp$ is $\Pi^1_{\gm'}$-indescribable.

As before, we produce the $\Bon^{<\gm}$-sequence in two steps: first we define $$S'=\la (C'_\alp,\eta_\alp):\alp\in Reg\cap \kp\ra$$ which is a $\Bon^{<\gm}$-sequence below $\kp$, then we set
$$E_{A}=\{\alp\in A\cap Reg:C'_\alp\cap A\text{ is not $\eta_\alp$-stationary in $ \alp$ or }d_{\eta_\alp}(C'_\alp)=\emp \}$$
and we construct a refinement of $S'$ which avoids $E_A$.
To streamline notation for this proof, we shall  abbreviate $\trace^{\Pi^1_\gm}$ to $\trace^\gm$.

\emph{Constructing $S'$}:

As $\kp$ is not $\Pi^1_{\gm}$-indescribable we can fix a $\Delta_0$ formula $\phi(v_0,v_1,v_2)$ and $Z\sbs\kp$ such that:
$$ \Pi \text{ wins } G_{\gm}(\kp,\phi, Z)$$
but for all $\alp<\kp$
$$ \Sigma \text{ wins } G_{\gm}(\alp,\phi, Z\cap\alp)$$
and so for each $\alp<\kp$
$$\exists \eta<\gm, X\sbs\alp \text{ such that } \Pi \text{ wins } G_{\eta}(\alp,\phi', \la Z\cap\alp,X\ra).$$


Let $\alp\in \kp$. Let $\bar{\eta}_\alp$ be be maximal such that $A\cap\alp$ is $\bar{\eta}_\alp$-stationary (or $0$, if $A$ is not even unbounded), and let $D_\alp$ be the $<_L$ least $\bar{\eta}_\alp$-club avoiding $A\cap\alp$. We have $\Sigma \text{ wins } G_{\gm}(\alp,\phi, Z\cap\alp)$, so we set $Y_\alp$ to be $<_{L}$-least subset of $\alp$ such that $\Sigma$ wins with first move $Y_\alp$, if $\gm$ is odd, and with first move $\la\eta,Y_\alp\ra$ with $\eta$ minimal, if $\gm$ is even. We split into cases and define $\eta_\alp$ and $X_\alp$:\\

\textbf{Case (i)}: $\gm$ is odd,  $A\cap \alp$ is $\gm-1$-stationary, and if $\bar{\eta}_\alp=\gm-1$ we have $Y_\alp<_L D_\alp$.
Set $\eta_\alp=\gm-1$ and $X_\alp=Y_\alp$.\\

\textbf{Case (ii)}: If $\gm$ is even, let $\eta$ be least (odd) ordinal such that $\Sigma$ wins the game $G_\gm(\alp,\phi, Z)$ with first ordinal move $\eta$ (choosing $\alp_1=\eta$), \ie, 
$$\Sigma \text{ wins }G_{\eta+1}(\alp,\phi, Z) \text{ but for any for any even $\eta'<\eta$, } \Pi\text{ wins }G_{\eta'}(\alp,\phi, Z).$$  
 We put $\alp$ under case (ii) if $A\cap\alp$ is $\eta$-stationary and if $\bar{\eta}_\alp=\eta$ then $Y_\alp<_L D_\alp$. In this case, set $\eta_\alp=\eta$ and $X_\alp=Y_\alp$.\\

\textbf{Case (iii)}: Otherwise: set $\eta_\alp=\bar{\eta}_\alp$ and set $X_\alp=D_\alp$. Note that if $\alp$ falls into this case then $\eta_\alp<\gm$, as if $\bar{\eta}_\alp\geq\gm$ then $A\cap\alp$ is $\gm$-stationary so the conditions of (i) or (ii) are fulfilled. \\

We take $\nu_\alp>\alp$ to be the least limit ordinal such that $L_{\nu_\alp}$ is $\Pi^1_{\eta_\alp}$ correct over $\alp$ and $X_\alp,A\cap\alp,Z\cap\alp\in L_{\nu_\alp}$. We set $p_\alp=\{X_\alp,A\cap\alp,Z\cap\alp\}$ and then we set:

\begin{displaymath}
C'_\alp = \left \lbrace 
\begin{array}{l}
\trace^{{\eta_\alp}}(L_{\nu_\alp},p_\alp,\alp) \text{\,\,\, if this is $\eta_\alp$-club in }\alp\\
\text{ an arbitrary non-reflecting $\eta_\alp$-stationary set otherwise}\\
\end{array} \right. 
\end{displaymath}

\noindent This is well defined as we know $\alp$ is $\eta_\alp$ stationary and a cardinal, and so if the trace is not $\eta_\alp$-club we must have that $\alp$ is not $\Pi^1_{\eta_\alp}$-indescribable so a non-reflecting set can be found. 

We set $S'=\la (C'_\alp,\eta_\alp):\alp\in \kp\cap Reg\ra$. 

\begin{claim}
$S'$ is a $\Bon^{<\gm}$-sequence.
\end{claim}

{\sc Proof:}
It is immediate from the definition that each $C'_\alp$ is $\eta_\alp$-club, so we just need to show that we have the coherence property. So let $\alp<\kp$ be regular and suppose $C'_\alp$ is defined as the trace (otherwise $d_{\eta_\alp}(C'_\alp)=\emp$ so coherence is trivial). Let $\bt\in d_{\eta_\alp}(C'_\alp)$. Let $N_\bt=L_{\nu_\alp}\{p_\alp\cup\bt\cup\{\alp\}\}$ and $\pi:N_\bt\cong L_{\bar{\nu}_\bt}$ be the collapsing map. We need to show $\eta_\bt=\eta_\alp$, $\pi``p_\alp=p_\bt$ and $\nu_\bt=\bar{\nu}_\bt$. 

Clearly $\pi(A\cap\alp)=A\cap\bt$ and $\pi(Z\cap\alp)=Z\cap\bt$. Let $X=\pi(X_\alp)=X_\alp\cap\bt$. 

If $\alp$ falls into case (i) then $\gm$ is odd and
 $$L_{{\nu}_\alp}\vDash ``A\cap\alp\text{ is $\eta_\alp$-stationary and for any $D<_L X$ if $D$ is $\eta_\alp$-club then } A\cap d_{\eta_\alp}(D)\neq\emp"$$ 
 and so by elementarity $$L_{\bar{\nu}_\bt}\vDash ``A\cap\bt\text{ is $\eta_\alp$-stationary and for any $D<_L X$ if $D$ is $\eta_\alp$-club then } A\cap d_{\eta_\alp}(D)\neq\emp"$$
thus by $\Pi^1_{\eta_\alp}$ correctness of $L_{\bar{\nu}_\bt}$ we have that we are again in case (i) with $\eta_\bt=\eta_\alp$.
 
If $\alp$ falls into case (ii), we have 
$$L_{{\nu}_\alp}\vDash ``\Sigma \text{ wins } G_{\eta_\alp+1}(\alp,\phi, \la Z\cap\alp\ra) \text{ with first move }X_\alp "$$
\ie $$L_{{\nu}_\alp}\vDash``\Pi \text{ wins } G_{\eta_\alp}(\alp,\phi', \la Z\cap\alp,X_\alp\ra)".$$
Now by elementarity,
\begin{align*}
L_{\bar{\nu}_\bt}  \vDash  & ``\Pi\text{ wins } G_{\eta_\alp}(\bt,\phi', \la Z\cap\bt,X\ra), A\cap\bt\text{ is $\eta_\alp$-stationary} \\
& \text{and for any $D<_L X$ if $D$ is $\eta_\alp$-club then } A\cap d_{\eta_\alp}(D)\neq\emp."
\end{align*}
As this is a $\Pi^1_{\eta_\alp}$ statement we have $\Pi \text{ wins } G_{\eta_\alp}(\bt,\phi', \la Z\cap\bt,X\ra)$ by $\Pi^1_{\eta_\alp}$-correctness of $L_{\bar{\nu}_\bt}$, and thus $\Sigma \text{ wins } G_{\eta_\alp+1}(\bt,\phi, Z\cap\alp)$ with first move $X$. We also have $A\cap\bt$ is $\eta_\alp$-stationary etc. and so $\bt$ falls into case (ii) and $\eta_\bt\leq \eta_\alp$.

For even $\eta \leq\eta_\alp$ we have $L_{\nu_\alp}\vDash \Pi \text{ wins } G_{\eta}(\alp,\phi, \la Z\cap\alp\ra)$ and so by elementarity and $\Pi^1_{\eta_\alp}$-correctness, for such $\eta$ we have $\Pi \text{ wins } G_{\eta}(\bt,\phi, \la Z\cap\bt\ra)$ and so we must have $\eta<\eta_\bt$. Thus $\eta_\bt=\eta_\alp$ in this case. 
 
Finally, if $\alp$ falls into case (iii) then we have $\eta_\alp=\bar{\eta}_\alp$ and $X=\pi(D_\alp)$ so by definition of the case and elementarity,
\begin{align*}
L_{\bar{\nu}_\bt}\vDash & ``A\text{ is $\eta_\alp$-stationary and $X$ is $\eta_\alp$-club avoiding }A \text{ and for all $Y<_L  X$,} \\
& \Sigma \text{ does not win } G_{\eta_\alp+1}(\bt,\phi, \la Z\cap\bt\ra)\text{ with first move }Y."
\end{align*}
 Thus by $\Pi^1_{\eta_\alp}$-correctness we are in case (iii) and $\eta_\alp=\eta_\bt$. 

Now in all cases $\eta_\bt=\eta_\alp$, so we set $\eta=\eta_\alp=\eta_\bt$.
To show $X=X_\bt$: we have already seen that $\bt$ falls into the same case as $\alp$. For case (i) or (ii) we have 
$$L_{\nu_\alp}\vDash \forall U\sbs\alp \,[U<_L X_\alp\rightarrow \Sigma \text{ wins } G_{\eta}(\alp,\phi', \la Z\cap\alp, U\ra)]$$
Then by elementarity we have: 
$$L_{\bar{\nu}_\bt}\vDash \forall U\sbs\bt \,[U<_L X\rightarrow \Sigma \text{ wins } G_{\eta}(\bt,\phi', \la Z\cap\bt, U\ra)]$$
and so by absoluteness
$$\forall U\sbs\bt\,[ U<_L X \rightarrow L_{\bar{\nu}_\bt}\vDash \Sigma \text{ wins } G_{\eta}(\bt,\phi', \la Z\cap\bt, U\ra)]$$
Thus by $\Pi^1_\eta$-correctness of $L_{\bar{\nu}_\bt}$ we do not have have $X_\bt<_L X$. Also $$L_{\bar{\nu}_\bt}\vDash \Pi \text{ wins } G_{\eta}(\bt,\phi', \la Z\cap\bt, X\ra)$$ and so by $\Pi^1_\eta$-correctness $X_\bt=X$. 
In case (iii), we can repeat the same argument but instead of $\phi$ take the sentence giving us that $A\cap\alp$ is not $\eta_\alp+1$-stationary.

It remains to show that $\nu_\bt=\bar{\nu}_\bt$. We already have that $p_\bt\sbs  L_{\bar{\nu}_\bt}$ and that  $L_{\bar{\nu}_\bt}$ is $\Pi^1_\eta$-correct over $\bt$, so we only need to show the minimality requirement.
Now for each limit ordinal $\gm>\alp$ with $\gm<\nu_\alp$ we have:
\begin{displaymath}
\begin{split}
\Theta(\gm): &\big(X_\alp\notin L_\gm\big)\vee\big( A\cap\alp\notin L_\gm\big)\vee \big(Z\cap\alp\notin L_\gm\big)\\
& \vee \exists U\notin L_\gm~\exists n\in\om~\exists \eta'<\eta~\big(U\text{ is minimal with }\psi_{\eta'}(n, U)\big)
\end{split}
\end{displaymath}
where $\psi_\eta(.,.)$ is the universal $\Pi^1_{\eta}$ sentence. 
We finish off as before.
By $\Pi^1_\eta$-correctness we have for each $\gm<\nu_\alp$ that $L_{\nu_\alp}\vDash\Theta(\gm)$ and thus $$L_{\nu_\alp}\vDash\forall\gm\Theta(\gm)\Rightarrow  L_{\bar{\nu}_\bt}\vDash\forall\gm\Theta(\gm)\Rightarrow \forall\gm<\bar{\nu}_\bt\,\,\,L_{\bar{\nu}_\bt}\vDash\Theta(\gm)$$
So by $\Pi^1_\eta$-correctness of $L_{\bar{\nu}_\bt}$ we have $\Theta(\gm)$ for each limit $\gm<\bar{\nu}_\bt$, and $\nu_\bt=\bar{\nu}_\bt$.
\qed

We can now define a part of our $\Bon^{<\gm}$-sequence that will satisfy Theorem \ref{squaregm}. This will be a refinement of $S'$, so the $\eta_\alp$'s that were defined above will stay the same. We set $$\Gm^1=\{\alp\in \kp:C'_\alp\cap A\text{ is $\eta_\alp$-stationary in }\alp\}.$$ If $\alp\in\Gm^1$ we set 
\begin{displaymath}
C_\alp = \left \lbrace 
\begin{array}{l}
d_{\eta_\alp}(C'_\alp\cap A)\,\,\,\text{ if $d_{\eta_\alp}(C'_\alp\cap A)$ is $\eta_\alp$-stationary}\\
\text{an arbitrary non-reflecting $\eta_\alp$-stationary set otherwise.}\\
\end{array} \right. 
\end{displaymath}

\noindent It is clear that $\la (C_\alp,\eta_\alp):\alp\in\Gm^1\ra$ is a $\Bon^{<\gm}$-sequence on $\Gm^1$. Setting $$E=E_A=\{\alp\in A\cap Reg:C'_\alp\cap A\text{ is not $\eta_\alp$-stationary in $ \alp$ or }d_{\eta_\alp}(C'_\alp)=\emp \}$$

we show that each such $C_\alp$ avoids $E$. If $d_{\eta_\alp}(C_\alp)\neq\emp$ then for $\bt\in d_{\eta_\alp}(C_\alp)$ we have $d_{\eta_\alp}(C_\alp'\cap A)$ is $\eta_\alp$-stationary below $\bt$, so $C'_\bt=C'_\alp\cap\bt$ with $\eta_\alp=\eta_\bt$ and $d_{\eta_\bt}(C_\bt)\neq\emp$, and by Lemma \ref{nstatunion} $C'_\alp\cap A=C'_\bt\cap A$ is $ \eta_\alp=\eta_\bt$-stationary below $\bt$.

Now we need to define $C_\alp$ for $\alp\in \Gm^2=\{\alp\in \kp:C'_\alp\cap A\text{ is not $\eta_\alp$-stationary in }\alp\}=\kp\without\Gm^1$. For such $\alp$ we shall find $C_\alp\sbs C'_\alp$ such that (i) $C_\alp$ avoids $A$ and hence $E$ and (ii) for $\bt\in d_{\eta_\alp}(C_\alp)$ we have $C'_\bt\cap A$ is not $\eta_\alp$-stationary (\ie $\bt\in\Gm^2$) and  $C_\bt=C_\alp\cap\bt$ (we shall already have $\eta_\bt=\eta_\alp$ as $C_\alp\sbs C'_\alp$).
Once we have (i) and (ii) it is easy to see that $S=\la C_\alp:\alp\in\kp\ra$ will satisfy Theorem \ref{squaregm}, and it will only remain to show that $E$ is $\gm$-stationary.

For $\alp\in\Gm^2$ we take $\rho_\alp\geq\nu_\alp$ minimal limit ordinal such that $C'_\alp\in L_{\rho_\alp}$ and $L_{\rho_\alp}$ is $\Pi^1_{\eta_\alp}$-correct over $\alp$. Note that if $D_\alp$ is the $<_L$ least $\eta$-club, for some $\eta<\eta_\alp$, that avoids $A\cap C'_\alp$ then $D_\alp\in L_{\rho_\alp}$ by $\Pi^1_{\eta_\alp}$-correctness. We now set 
\begin{displaymath}
C_\alp = \left \lbrace 
\begin{array}{l}
D_\alp\cap\trace^{{\eta_\alp}}(L_{\rho_\alp},\{C'_\alp\}\cup p_\alp,\alp) \text{\,\,\, if this is $\eta_\alp$-club in }\alp\\
\text{ an arbitrary non-reflecting $\eta_\alp$-stationary set otherwise.}\\
\end{array} \right. 
\end{displaymath}

Then we have $C_\alp\sbs d_{\eta_\alp}(C'_\alp)\sbs C'_\alp$ and it is clear that we have (i): $C_\alp$ avoids $A$. For $\bt\in d_{\eta_\alp}(C_\alp)$ we have $\bt\in d_{\eta_\alp}(C'_\alp)$ hence $C'_\bt=C'_\alp\cap\bt$. Then we have $A\cap C'_\bt$ is not $\eta_\alp$-stationary in $\bt$ by elementarity and $\Pi^1_{\eta_\alp}$-correctness, so $\bt\in\Gm^2$. Also $D_\bt\cap\trace^{\eta_\alp}(L_{\rho_\bt},\{C'_\bt\}\cup p_\bt,\bt)=\bt\cap D_\alp\cap\trace^{\eta_\alp}(L_{\rho_\alp},\{C'_\alp\}\cup p_\alp,\alp)=C_\alp\cap\bt$ so we are in the first case of the definition and $C_\bt=C_\alp\cap\bt$. This gives (ii) and so we have a $\Bon^{<\gm}$-sequence on the regulars avoiding $E_A$.

We have thus far only given our $\Bon^{<\gm}$-sequence on the regulars below $\kp$. To deal with singulars we  use 
Jensen's global $\Bon$-sequence just as before.  
It remains to show that $E$ is $\gm$-stationary. 
\begin{defn}
We define $H\sbs A$ by letting $\alp\in H$ iff $\alp\in A$ and there is $\mu_\alp >\alp$ and $q$ is a parameter from $L_{\mu_\alp}$ such that:
\begin{enumerate}
\item  $L_{\mu_\alp}$ is $\Pi^1_{\eta_\alp}$-correct over $\alp$
\item $\mu_\alp<\nu_\alp$
\item $A\cap\mathcal{s}^{{\eta_\alp}}(L_{\mu_\alp},q,\alp)=\emp$
\end{enumerate} 
\end{defn}

\begin{lem}
$H\sbs E$.
\end{lem}
{\sc Proof:}
Suppose $\alp\notin E$, so $C'_\alp\cap A$ is $\eta_\alp$-stationary  and $C'_\alp=\trace^{\eta_\alp}(L_{\nu_\alp},p_\alp,\alp)$. Let $\mu<\nu_\alp$ and $q\in L_\mu$. Then using lemma \ref{tracemono} for some $\bt<\alp$ we have $q\in L_{\nu_\alp}\{\bt\cup p_\alp\cup\{\alp\}\}$ and so $\trace^{\eta_\alp}(L_\mu,q,\alp)\without\bt\supseteq\trace^{\eta_\alp}(L_{\nu_\alp},p_\alp,\alp)=C'_\alp$. Thus $\alp\notin H$.\\
\qed

\begin{lem}
$H$ is $\gm$-stationary.
\end{lem}
{\sc Proof:}
Let $\eta<\gm$ and let $C\sbs\kp$ be $\eta$-club. Take $\mu>\kp$ minimal such that $C\in L_\mu$ and $L_\mu$ is $\Pi^1_\eta$-correct. Set $D=\trace^{\eta}(L_\mu,\{C,Z\},\kp)$ and note that $D\sbs d_\eta(C)\sbs C$, and $D$ is $\eta$-club by Lemma \ref{traceclub2}. Take $\dlt=min (D\cap A)$. Set $\mu_\dlt$ to be the ordinal such that  $L_\mu\{\{C,Z,\kp\}\cup\dlt\}\cong L_{\mu_\dlt}$.
 We show that $\dlt\in H$, with $\mu_\dlt$ and $\{C\cap\dlt,Z\cap\dlt\}$ witnessing this. 
 
\textbf{Claim 1}: $\eta_\dlt=\eta$\\
First to see $\eta_\dlt\geq \eta$. If $\dlt$ falls into case (i) then $\eta_\dlt=\gm-1$ so this is trivial. If $\dlt$ falls into case (iii) then we know that $A\cap\dlt$ is not $\eta_\dlt+1$-stationary, but $L_{\mu_\dlt}\vDash A\cap\dlt\text{ is $\eta$ stationary}$ and so by $\Pi^1_\eta$-correctness we must indeed have $A\cap\dlt$ is $\eta$-stationary and hence $\eta\leq\eta_\dlt$.

If $\dlt$ falls into case (ii) then $\eta_\dlt$ is odd and we know $\Sigma \text{ wins }G_{\eta_\dlt+1}(\alp,\phi, Z\cap\dlt)$, but for any even ordinal $\eta'\leq\eta$  we have  $L_{\mu}\vDash\Pi \text{ wins } G_{\eta'}(\kp,\phi, Z) $ and so by elementarity and $\Pi^1_\eta$-correctness of $L_{\mu_\dlt}$, we have $\Pi \text{ wins } G_{\eta'}(\dlt,\phi, Z\cap\dlt)$.
Thus we must have $\eta<\eta_\dlt+1$, \ie $\eta\leq\eta_\dlt$.

Now we have to show $\eta_\dlt\leq \eta$. First, suppose $\dlt$ is $\Pi^1_\eta$-indescribable. Then by Lemma \ref{traceclub2} we have $\trace ^{\eta}(L_{\mu_\dlt},\{C\cap\dlt,Z\cap\dlt\},\dlt)$ is $\eta$-club below $\dlt$. But $\trace ^{\eta}(L_{\mu_\dlt},\{C\cap\dlt,Z\cap\dlt\},\dlt)=\dlt\cap\trace ^{\eta}(L_{\mu},\{C,Z\},\kp)=D\cap\dlt$ so we must have $A$ is not $\eta+1$ stationary below $\dlt$. But in each of the three cases $A$ is $\eta_\dlt$-stationary, so  $\eta_\dlt\leq \eta$.


Now, if $\dlt$ is not $\Pi^1_\eta$-indescribable we know that $\dlt$ is not $\eta+1$-stationary, so again we have that $A$ is not $\eta+1$-stationary and as above $\eta_\dlt\leq \eta$.

\textbf{Claim 2}: $\mu_\dlt<\nu_\dlt$\\
We have that for every $X\in \P(\dlt)\cap L_{\mu_\dlt}$ 
$$L_{\mu_\dlt}\vDash \Sigma \text{ wins } G_{\eta}(\dlt,\phi', \la Z\cap\dlt, X\ra) .$$
As $L_{\mu_\dlt}$ is $\Pi^1_\eta$-correct this means $$ \forall X\in \P(\dlt)\cap L_{\mu_\dlt}, ~ \Sigma \text{ wins } G_{\eta}(\dlt,\phi', \la Z\cap\dlt, X\ra).$$ But we know by the choice of $\nu_\dlt$ that 
 $$ \exists X\in \P(\dlt)\cap L_{\nu_\dlt}, ~\Pi \text{ wins } G_{\eta}(\dlt,\phi', \la Z\cap\dlt, X\ra) .$$
 Hence we must have $\nu_\dlt>\mu_\dlt$
 
\textbf{Claim 3}: $A\cap\mathcal{s}^{\eta_\alp}(L_{\mu_\dlt},\{C\cap\dlt, Z\cap\dlt\},\dlt)=\emp$\\
We have (Lemma \ref{coltrace}) $A\cap\mathcal{s}^{\eta_\alp}(L_{\mu_\dlt},\{C\cap\dlt, Z\cap\dlt\},\dlt)=A\cap\dlt\cap\mathcal{s}^{\eta_\alp}(L_{\mu},\{C, Z\},\kp)=A\cap\dlt\cap\mathcal{s}^{\eta}(L_{\mu},\{C, Z\},\kp)=A\cap\dlt\cap D=\emp$.
\qedtwo{Theorem \ref{squaregm}}\\


We can now state Lemma \ref{traceclub2} without the assumption that for all $\gm'<\gm$ any regular cardinal which is $\gm'$-reflecting is $\Pi^1_{\gm'}$-indescribable, as this is a consequence of Theorem \ref{squaregm}.

\begin{lem}\label{traceclub3}$(V=L)$
Let $\gm$ be an ordinal and $\kp>\gm$ be $\Pi^1_{\gm}$-indescribable. Then for any limit $\nu>\kp$ with $L_\nu$ $\Pi^1_\gm$-correct over $\kp$ we have $\mathcal{s}^{\Pi^1_\gm}(L_\nu,p,\kp)$ is $\gm$-club in $\kp$.
\end{lem}

We finish by deriving a $\Bon^\gm$-sequence from a $\Bon^{<\gm+1}$-sequence, showing that Theorem \ref{squaren} and its generalisation to infinite $\gm$ is indeed an easy corollary of Theorem \ref{squaregm}. 

\begin{prop}
If $S=\la (C_\alp,\eta_\alp):\alp\in\Gm\ra$ is a $\Bon^{<\gm+1}$-sequence on $\Gm\sbs\kp$, then $\la C_\alp:\alp\in\Gm\cap d_\gm(\kp)\ra$ is a $\Bon^\gm$-sequence.
\end{prop}

{\sc Proof:}
Let $\alp\in\Gm\cap d_\gm(\kp)$. Then $C_\alp$ is $\eta_\alp$-club and $\eta_\alp\leq\gm$ and $\alp$ is $\gm$-stationary, so $C_\alp$ is $\gm$-club. Suppose $\bt\in d_\gm(C_\alp)$. Then as $\eta_\alp\leq\gm$, $\bt\in d_{\eta_\alp}(C_\alp)$, so $C_\bt=C_\alp\cap\bt$. Also $\bt\in d_\gm(\kp)$, so we have coherence.
\qed

\begin{col}
$(V=L)$
Let $\gm<\kp$ and $\kp$ be a  $\Pi^1_\gm$- but not $\Pi^1_{\gm+1}$-indescribable cardinal, and let $A\subseteq\kp$ be $\gm+1$ stationary. Then there is $E_{A}\sbs A$ and a $\Bon^\gm$-sequence $S$ on $\kp$ such that $E_{A}$ is $\gm+1$-stationary in $\kp$ and $S$ avoids $E_{A}$.
\end{col}
{\sc Proof:}
This follows easily from the preceding proposition - just restricting the domain of the $\Bon^{<\gm+1}$-sequence from Theorem \ref{squaregm} to $d_\gm(\kp)$ gives a $\Bon^\gm$-sequence avoiding $E_A$ as defined there.
\qed

\subsection{The $\gm$-club Filters and Non-Threaded $\Bon^\gm$-sequences} 


In this subsection we look in more detail at the $\gm$-club filter and what we can prove from certain assumptions about it. In the first subsection we shall revisit the results of 2.2 and generalise our result from there to $\gm$-stationarity. In the second subsection we generalise the notion of non-threaded $\Bon$, and show that this $\Bon^\gm(\kp)$ must fail at any $\gm+1$-reflecting cardinal where the $\gm$-club filter is normal. In the final section we show that with a certain requirement on the generalised club filters $\gm$-stationarity is downward absolute to $L$. If we make the stronger (and easier to state) assumption that for any ordinal $\eta$ the $\eta$-club filter is normal on any $\eta$-reflecting cardinal, then we have that for any ordinal $\gm$, $\gm$-stationarity is downward absoluteness to $L$ at any regular cardinal (Corollary \ref{dasimple}).



\subsubsection{The $\gm$-Club and $\Pi^1_\gm$-Indescribability Filters} 

Here we look in more detail at the $\gm$-club filter, revisiting the material from 2.2.1 in the light of the definition of $\Pi^1_\gm$-indescribability given in 3.3.

\begin{prop}
If $\Cf^\gm(\kp)$ is normal, then for any $\eta<\gm$ we have $\Cf^\eta(\kp)$ is also normal. 
\end{prop}

{\sc Proof:}
This is because each $\eta$-club is $\gm$-club, and so by the normality of $\Cf^\gm(\kp)$, the diagonal intersection of $\eta$-clubs must be $\gm$-club and hence $\eta$-stationary. The $\eta$-stationary closure is automatic.
\qed

\begin{prop}\label{nolimitnormality}
If $\gm$ is a limit ordinal and $\kp>\gm$ then $\bigcup_{\eta<\gm}\Cf^\eta(\kp)$ is not $\gm$ complete.
\end{prop}

{\sc Proof:}
Unless $\kp$ is $\gm$-stationary $\bigcup_{\eta<\gm}\Cf^\eta(\kp)$ is not even a filter, so suppose $\kp$ is $\gm$-stationary.
For $\eta<\gm$ set $C_\eta=d_\eta(\kp)$. Then each $C_\eta$ is $\eta$-club and $\bigcap_{\eta<\gm}C_\eta=d_\gm(\kp)$. By Proposition \ref{prop} we have that for each $\eta<\kp$, $\kp\without d_{\eta+1}(\kp)$ is $\eta+1$-stationary and hence $d_\gm(\kp)\notin \Cf^\eta(\kp)$. Thus $d_\gm(\kp)\notin\bigcup_{\eta<\gm}\Cf^\eta(\kp)$.

Recall (Definition \ref{gmindescribabilityfilter})
that the $\Pi^1_\gm$ indescribability filter on $\kp$, $\F^\gm(\kp)$ is the filter generated by sets of the form
$\{\alp<\kp:\Pi\text{ wins } G_\gm(\kp,\phi,A) \text{ and }G_\gm(\alp,\phi,A\cap\alp)\}$.

\begin{lem}\label{filtersV2}
If $\kp$ is $\Pi^1_\gm$-indescribable then the $\gm$-club on $\kp$ filter is a subset of $\F^\gm(\kp)$ and hence it is normal.
\end{lem}

{\sc Proof:}
The proof is essentially the same as Lemma \ref{filtersV1} for finite $\gm$. Firstly, any $\gm$-club is in $\F^\gm(\kp)$. Suppose $C$ is $\gm$-club. ``$C$ is $\gm$-club" is $\Pi^1_\gm$ and so reflects to a set in the $\Pi^1_\gm$-indescribability filter on $\kp$. But this is the set of $\alp<\kp$ such that $C\cap\alp$ is $\gm$-club, \ie $d_\gm(C)$. As $d_\gm(C)\sbs C$ we have that $C\in\F^\gm(\kp)$. Now suppose $\la C_\alp:\alp<\kp\ra$ is a sequence of $\gm$-clubs and set $C=\triangle_{\alp<\kp}C_\alp$. By the normality of $\F^\gm(\kp)$ and the fact that each $C_\alp\in\F^\gm(\kp)$, we have that $C\in \F^\gm(\kp)$ and hence $C$ is $\gm+1$-stationary. $\gm$-closure is easily verified.
\qed

\begin{col}(Fodor's Lemma for $\gm$-stationary sets)
If $\kp$ is $\Pi^1_\gm$-indescribable and $A\subset\kp$ is $\gm$-stationary, then for any regressive function $f:A\rightarrow \kp$ there is an $\gm$-stationary $B\subseteq A$ such that $f$ is constant on $B$.
\end{col}

The following generalisation of Theorem \ref{filtersL1} is a consequence of Lemma \ref{traceclub2} and Theorem \ref{squaregm}.
\begin{col}\label{filtersL2}
If $V=L$ the $\gm$-club filter coincides with the $\Pi^1_\gm$-indescribability filter at any $\Pi^1_\gm$-indescribable cardinal.
\end{col}

{\sc Proof:}
Suppose $\Pi$ wins $G_\gm(\kp,\phi,X)$. Then $\trace^{\Pi^1_\gm}(L_{\kp^+},\{X\},\kp)$ is $\gm$-club by Lemma \ref{traceclub3}, and for each $\alp\in\trace^{\Pi^1_\gm}(L_{\kp^+},\{X\},\kp)$, we have  $\Pi$ wins $G_\gm(\alp,\phi,X\cap\alp)$. Thus $\{\alp<\kp:\Pi$ wins $G_\gm(\alp,\phi,X\cap\alp)\}$ is in the $\gm$-club filter and in general $F^\gm(\kp)$ is included in the $\gm$-club filter. By \ref{filtersV2} we have the reverse inclusion.
\qed

\subsubsection{Splitting Stationary Sets}\label{sectionsplitting2}

In this section we extend the results of 2.2.2 to split $\gm+1$-stationary sets. We should note here that the methods of section 2.2.2 and what follows do not allow us to show that $\gm$-stationary sets can be split when $\gm$ is a limit ordinal. This is because by Proposition \ref{nolimitnormality}, if $\gm$ is a limit ordinal then the filter corresponding to the $\gm$-stationary sets, $\bigcup_{\eta<\gm}\Cf^\eta(\kp)$, is not $\gm$ complete, and our results require $\kp$ completeness.

\begin{lem}\label{split2}
If $\Cf^{\gm}(\kp)$ is $\kp$ complete then any $\gm+1$-stationary subset of $\kp$ is the union of two disjoint $\gm+1$-stationary sets.
\end{lem}

This lemma is proven in the same way as Lemma \ref{split}, the key points to enable the generalisation being Lemma \ref{filtersV2} and the fact that a measurable $\kp$ is $\Pi^1_\gm$-indescribable for any $\gm<\kp$.

{\sc Proof:}
Let $S$ be $\gm+1$-stationary in $\kp$ and suppose $S$ is not the union of two disjoint $\gm+1$-stationary sets. Define 
$$F=\{X\subset\kp:X\cap S \text{ is $\gm+1$-stationary}\}$$

\emph{Claim:} $F$ is a $\kp$ complete ultrafilter.\\
Upwards closure is clear. Intersection follows by the fact that $S$ cannot be split, as well as $X\in F\Rightarrow \kp\without X\notin F$. That $X\notin F\Rightarrow \kp\without X\in F$ follows by the definition of $\gm$-stationary, and $\kp$ completeness follows from the $\kp$ completeness of the $\gm$-club filter.

\emph{Claim:} $F$ is normal.\\
As we have shown that $F$ is a $\kp$ complete ultrafilter on $\kp$, we have that $\kp$ is measurable. 
Now, all measurables are $\Pi^2_1$-indescribable (Lemma \ref{measurables}), and it is easy to see that $\Pi$ winning the game $G_\gm(\kp,\phi,A)$ for $\phi$ a $\Delta_0$ formula and $A\sbs\kp$ is $\Pi^2_1$ expressible over $V_\kp$, as in describing the game we only quantify over finite sequences of subsets of $\kp$.  Hence $\kp$ is $\Pi^1_{\gm}$-indescribable and so by Lemma \ref{filtersV2}, $\Cf^{\gm}(\kp)$ is normal. Let $\la X_\alp:\alp<\kp\ra$ be a sequence of sets in $F$. Then each $S\without X_\alp$ is in the non-$\gm+1$-stationary ideal on $\kp$, so $X_\alp\cup (\kp\without S)\in\Cf^{\gm}(\kp)$. Now by the normality of $\Cf^{\gm}(\kp)$, we have for $X:=\bigtriangleup_{\alp<\kp}X_\alp\cup(\kp\without S)$ that $X\in\Cf^{\gm}(\kp)$, and so $X\cap S$ is $\gm+1$-stationary. But $X=\{\alp<\kp:\forall\bt<\kp\, \,\alp\in X_\bt\cup(\kp\without S)\}=\{\alp<\kp:\alp\in\kp\without S\vee \forall\bt<\kp\,\,\alp\in X_\bt\}=(\kp\without S)\cup\bigtriangleup_{\alp<\kp}X_\alp$. So $X\cap S= \bigtriangleup_{\alp<\kp}X_\alp\cap S$, and thus $\bigtriangleup_{\alp<\kp}X_\alp\cap S$ is $\gm+1$-stationary and hence in $F$.\\

Now we have that $F$ is a normal measure, so the second part of Lemma \ref{measurables} gives us that for any $R\sbs V_\kp$ and  formula $\phi$ that is $\Pi^2_1$, if $\la V_\kp,\in R\ra \vDash \phi$ then $$\{\alp<\kp:\la V_\alp,\in, R\cap V_\alp\ra\vDash \phi\}\in F.$$ Setting $R=S$ and $\phi=``S$ is $\gm+1$-stationary$"$ we can conclude
$$\{\alp<\kp:S\cap\alp\text{ is $\gm+1$-stationary}\}\in F.$$

Then by definition of $F$ we have $A:=\{\alp\in S:S\cap\alp\text{ is $\gm+1$-stationary}\}$ is $\gm+1$-stationary.
But by Proposition $\ref{prop}$ we have  $A':=\{\alp\in S:S\cap\alp \text{ is not $\gm+1$-stationary}\}$ is $\gm+1$-stationary. This contradicts our assumption on $S$ as $A'$ and $A$ are two disjoint, $\gm+1$-stationary subsets of $S$.
\qed

\begin{thm}\label{splittinggm}
If $\kp$ is weakly compact and $C^{\gm}(\kp)$ is $\kp$-complete then any $\gm+1$-stationary subset of $\kp$ is can be split into $\kp$ many disjoint $\gm+1$-stationary sets.
\end{thm}

{\sc Proof:} Let $S\subseteq\kp$ be $\gm+1$-stationary. We construct a tree $T$ of $\gm+1$-stationary subsets of $S$ ordered by $\supseteq$, using Lemma \ref{split2}: each $\gm+1$-stationary set can be split in two. $T$ will have height $\kp$ and levels of size $<\kp$ so as $\kp$ has the tree property it has a $\kp$ length branch, from which we can construct a partition of $S$.

We define $T$ inductively such that each level (i) consists of disjoint $\gm+1$-stationary sets, (ii) has size $<\kp$, and (iii) is non-empty. Let $T_0=\{S\}$. 
Now let $\alp<\kp$ and suppose we have defined $T_\bt$ for each $\bt<\alp$ such that (i)-(iii) hold. 
If $\alp$ is a successor, say $\bt+1$, for each set $A\in T_\bt$ we use Lemma \ref{split2} to choose $A'\subseteq A$ such that $A'$ and $A\without A'$ are both $\gm+1$-stationary, and take $T_{\bt+1}=\{A',A\without A':A\in T_\bt\}$. Clearly, (i)-(iii) are preserved. 

If $\alp$ is a limit, we define 
$$T_\alp=\{\bigcap b:b\text{ is a branch in $T_{<\alp}$ and $\bigcap b$ is $\gm+1$-stationary}\}.$$  
Now (i) is clear and as $|\{\bigcap b:b\text{ is a branch in } T_{<\alp}\}|\leq 2^{|T_{<\alp}|}<\kp$ we also have (ii). To demonstrate (iii) first note that by the construction of $T$ we have $S=\bigcup \{\bigcap b:b\text{ is a branch in }T_{<\alp}\}$ (For each $a \in S$, $\{A\in T_\alp:a\in A\}$ is clearly a branch although it may have height $<\alp$). Now we know $|\{\bigcap b:b\text{ is a branch in }T_{<\alp}\}|<\kp$ so by $\kp$-completeness of $\Cf^{\gm}(\kp)$ we must have at least one branch $b$ such that $\bigcap b$ is $\gm+1$-stationary. But then $b$ must have height $\alp$, for if $b$ had height $\bt<\alp$, then $b$ must have limit height so by definition $\bigcap b\in T_\bt$ which would make $b$ not maximal. So $\bigcap b\in T_\alp$, and hence (iii) holds (and the definition of $T_\alp$ makes sense).

Now by the tree property $T$ has a $\kp$ branch. Take $B$ to be such a branch. For $A\in B$ let $A^+$ denote the immediate successor of $A$ in $B$, and note that $A\without A^+$ is disjoint from any set in $B$ succeeding $A$. Then $\{A\setminus A^+:A\in B\}$ is a partition of $S$ into $\kp$ many disjoint $\gm+1$-stationary sets.
\qed

\begin{col}
Let $\kp$ be $\Pi^1_{\gm}$-indescribable, $\gm\geq 1$. Then any $\gm+1$-stationary subset of $\kp$ can be split into $\kp$ many disjoint $\gm$-stationary sets.
\end{col}

{\sc Proof:}
A $\Pi^1_{\gm}$-indescribable cardinal $\kp$ is of course weakly compact, and by Lemma \ref{filtersV2}, the $\gm$-club filter on $\kp$ is normal and hence $\kp$-complete.
\qed

It is also easy to see that if $\kp$ is inaccessible and $C^{\gm}_\kp$ is $\kp$-complete then for any $\alp<\kp$ we can split any $\gm+1$-stationary set into $\alp$ many $\gm+1$-stationary pieces.


\subsubsection{Non-Threaded $\Bon^\gm$-sequences}

In this section we generalise a nice folklore result\footnote{With thanks to Philipp L\"ucke for pointing out this result as something which may generalise.} that relates \emph{non-threaded $\Bon$} (Definition \ref{nonthreadeddef}) to simultaneous stationary reflection: if we have a non-threaded $\Bon$-sequence on a cardinal $\kp$ then there are two stationary subsets of $\kp$ which are not both stationary in any $\alp<\kp$. Thus non-threaded $\Bon$ fails at any $2$-stationary cardinal. We start by recalling the definitions.

\begin{defn}\label{squaredef} 
A $\Bon^\gm$\emph{-sequence} (for $\gm\geq 0$) on $\kp$ is a sequence $\la C_\alp:\alp\in d_\gm(\kp)\ra$ such that for each $\alp$:
\begin{enumerate}
\item $C_\alp$ is a $\gm$-club subset of $\alp$
\item for every $\bt\in d_\gm(C_\alp)$ we have $C_\bt=C_\alp\cap\bt$
\end{enumerate}

We say a $\Bon^\gm$-sequence $\la C_\alp:\alp\in d_\gm(\kp)\ra$ \emph{avoids} $A\subset\kp$ if for all $\alp$ we have $A\cap d_\gm(C_\alp)=\emp$.
\end{defn}

The generalisation of non-threaded $\Bon$ is straight-forward:

\begin{defn}\label{nonthreadeddef}
$\Bon^\gm(\kp)$ is the statement that there is a $\Bon^\gm$-sequence on $\kp$ which has no thread, \ie  there is no $C\sbs\kp$, such that $C$ is $\gm$-club and for all $\alp\in d_\gm(C)$ we have $C_\alp=C\cap\alp$. Such a sequence is called a $\Bon^\gm(\kp)$-sequence.
\end{defn}

\begin{thm}\label{nonthreaded}
Suppose $\kp$ is a regular $\gm$-reflecting cardinal and the $\gm$-club filter on $\kp$ is normal. Let $S=\la C_\alp:\alp\in d_\gm(\kp)\ra$  be a $\Bon^\gm$-sequence. Then the following are equivalent:
\begin{enumerate}
\item $S$ is a $\Bon^\gm(\kp)$-sequence, \ie $S$ has no thread.
\item For any $\gm+1$-stationary set $T$ there are $\gm+1$-stationary $S_0$, $S_1\sbs T$ such that for any $\alp\in d_\gm(\kp)$ we have $d_\gm(C_\alp)\cap S_0=\emp$ or $d_\gm(C_\alp)\cap S_1=\emp$
\end{enumerate}
Thus $\Bon^\gm(\kp)$ implies $\kp$ is not $\gm+1$-reflecting.
\end{thm}

We give the proof for $\gm>0$, following the proof given in \cite{lambiehanson}, though the generalisation is not straight-forward. For $\gm=0$ see \cite{lambiehanson} Proposition 27 or the following but adding a ``$-1$-club", where $C\sbs\kp$ is $-1$-club if $C$ is an end-segment of $\kp$.

{\sc Proof:}

We start with $(2)\rightarrow(1)$, so let $S_0$ and $S_1$ satisfy (2) and assume for a contradiction that $C$ is an $\gm$-club subset of $\kp$ which threads $S$.  Then we have $\alp<\bt<\dlt$ in $d_\gm(C)$ such that $\alp\in S_0$ and $\bt\in S_1$. Now $C_\dlt=C\cap\dlt$ so we have $\alp\in C_\dlt\cap S_0$ and $\bt\in C_\dlt\cap S_1$ - but this is a contradiction.\\

Now suppose $S$ has no thread and $T\sbs\kp$ is $\gm+1$-stationary. We split into two cases.\\
\textbf{Case 1:}
There is an $\eta<\gm$ and an $\eta$-club set $D$ such that $\{\alp<\kp:d_\gm(C_\alp)\cap D=\emp\text{ or } \alp\notin T\}$ contains an $\gm$-club. 

If we define $T'=\{\alp\in T\cap D:d_\gm(C_\alp)\cap D=\emp\}$ then $T'$ is $\gm+1$-stationary as it is the intersection of a $\gm$-club with $T$.   Firstly, suppose for each $\alp\in d_\gm(\kp)$ we have $|d_\gm(C_\alp)\cap T'|\leq 1$. Then for any pair $S_0$ and $S_1$ of disjoint $\gm+1$-stationary subsets of $T'$  and any $\alp$ we must have $d_\gm(C_\alp)\cap S_0=\emp$ or $d_\gm(C_\alp)\cap S_1=\emp$, so we're done. 

If this is not the case, we show that $S$ avoids $T'$ and hence we can use Lemma \ref{split2} to split $T'$ into 2 disjoint $\gm+1$-stationary sets, which will give (2). So suppose for some $\alp\in d_\gm(\kp)$ that $|d_\gm(C_\alp)\cap T'|\geq 2$. Take $\bt_0<\bt_1$ both in $d_\gm(C_\alp)\cap T'$. Then we have $C_{\bt_1}\cap\bt_0=C_\alp\cap\bt_0=C_{\bt_0}$, so $\bt_0\in d_\gm(C_{\bt_1})$ - but this is a contradiction with the definition of $T'$ since $\bt_1\in T'$ but $\bt_0\in D$.\\

\textbf{Case 2:}
For any $\eta<\gm$ and any $\eta$-club $D$ we have $\{\alp\in T:d_\gm(C_\alp)\cap D\neq\emp\}$ is $\gm+1$-stationary. 

For $\alp\in d_\gm(\kp)$ we set $$S_0^\alp=\{\bt\in T\without \alp:\alp\notin d_\gm(C_\bt)\}$$ and $$S_1^\alp=\{\bt\in T\without \alp:\alp\in d_\gm(C_\bt)\}.$$ Then it is clear that for any given $\alp$ either $S_0^\alp$ or $S_1^\alp$ is $\gm+1$-stationary.\\

{\em Claim}
There is some $\alp\in d_\gm(\kp)$ such that $S_0^\alp$ and $S_1^\alp$ are both $\gm+1$-stationary.

{\sc Proof of Claim:}
Suppose not. Then for each $\alp$ we have either $S_0^\alp$ is not $\gm+1$-stationary or $S_1^\alp$ is not $\gm+1$-stationary. Set $A=\{\alp\in\kp:S_0^\alp\text{ is not $\gm+1$-stationary}\}$. We claim that $A$ is $\gm$-stationary. So let $D$ be $\eta$-club for some $\eta<\gm$. Then as we are in Case 2,  $\{\alp\in T:d_\gm(C_\alp)\cap D\neq\emp\}$ is $\gm+1$-stationary. Define a regressive function on $T$ by $\alp\mapsto min(d_\gm(C_\alp)\cap D)$. By normality of the $\gm$-club filter, we can fix $T'\sbs T\cap D$ and $\alp$ such that for any $\bt\in T'$ we have $min(d_\gm(C_\bt)\cap D)=\alp$. Then $S_1^\alp\supseteq T'$ so $S^\alp_1$ is $\gm+1$-stationary and so we have $\alp\in A\cap D$. As $\eta<\gm$ and $D$ were arbitrary and we have shown that $A\cap D\neq \emp$ we can conclude that $A$ is $\gm$-stationary.

Now let $\alp<\alp'$ be elements of $A$ and $D^\alp$, $D^{\alp'}$ by $\gm$-clubs witnessing, respectively, that $S_0^\alp$ and $S_0^{\alp'}$ are not $\gm+1$-stationary. Fix $\bt\in D^\alp\cap D^{\alp'}\cap T$.
Then $\alp\in d_\gm(C_\bt)$ so $C_\alp=C_\bt\cap\alp$. Similarly, $C_{\alp'}=C_\bt\cap\alp'$ so we have $C_\alp=C_{\alp'}\cap\alp$. Setting $C=\bigcup_{\alp\in A}C_\alp$ we have that $C$ is $\gm$-club: as each $C_\alp$ is $\gm$-closed $C$ is $\gm$-closed and as the $\gm$-stationary union of $\gm$-stationary sets, $C$ must be stationary. But then $C$ is a thread through $S$ - contradiction.
\qedtwo{{\em Claim}}\\

Now let $\alp$ be such that $S_0^\alp$ and $S_1^\alp$ are both $\gm+1$-stationary. Set $S_0=S_0^\alp$ and $S_1=S_1^\alp$ and note that for any $\bt\in S_0\cup S_1$ we have $\bt>\alp$. We show these sets satisfy (2). Suppose $\dlt\in d_\gm(\kp)$. If $\bt\in d_\gm(C_\dlt)\cap S_0$ and $\bt'\in d_\gm(C_\dlt)\cap S_1$ then $\alp\in C_\bt$ and $\alp\notin C_{\bt'}$ by definition of $S_0$ and $S_1$, but $C_\bt\cap\bt'=C_\dlt\cap\bt\cap\bt'=C_{\bt'}\cap\bt$. As $\alp<\bt\cap\bt'$ this gives us a contradiction. Hence $C_\dlt$ must avoid either $S_0$ or $S_1$. \\ \mbox{ } \hfill
\qedtwo{Theorem \ref{nonthreaded}}


\section{Downward Absoluteness}

The downward absoluteness to $L$ of a cardinal being $1$-reflecting was proven by Magidor in \cite{magidorrss} \S 1 (the theorem was stated there for $\kp=\om_2$, but the proof is the same for any regular $\kp$). There, the only assumption needed on $\kp$ was regularity. To generalise this and get $\gm$-reflecting cardinals in $L$, we shall require some extra assumptions, as we see below. The proof works inductively, and we shall need a slightly stronger statement than the downward absoluteness of $\kp$ being $1$-reflecting - the downward absoluteness of $2$-stationarity for sets in $L$. 
At higher levels of stationarity the first assumption we require is the normality of a certain club filter - in the case of $1$-reflecting cardinals the club filter was guaranteed to be normal by assuming $\kp$ is regular.  we shall also need a further assumption that ``many" cardinals below $\kp$ have the properties guaranteeing downward absoluteness of lower levels of stationarity. 

The proof will be split into three cases - limit ordinals, successors of limit ordinals, and double successors. The proof for double successors case is based on Magidor's proof (which essentially gives downward absoluteness of $2$-stationarity) though there is more work to be done for the higher levels as we do not have absoluteness of $\gm$-clubs for $\gm>0$. This will be seen particularly in the last part of the proof. The limit stages are straightforward, but for the successors of limits we need a slight variation on the notion of normality. The following definition gives us the appropriate notion.

\begin{defn}\label{normaldef}
For $\gm>0$ we call a cardinal $\kp$ \emph{$\gm$-normal} if $\kp$ is $\gm$-reflecting and for any $\gm+1$-stationary $S\sbs\kp$ and regressive function $f:S\rightarrow\kp$, $f$ is constant on a $\gm$-stationary set.
\end{defn}

For successor ordinals this reduces to normality of the appropriate club filter:

\begin{prop}\label{normalnormal}
If $\gm>\eta$ then $\kp$ is $\gm$-normal implies $\Cf^\eta(\kp)$ is normal. If $\gm=\eta+1$ and $\kp$ is $\gm$-reflecting then we have the converse: $\Cf^\eta(\kp)$ being normal implies $\kp$ is $\eta+1$-normal.
\end{prop}

Thus $\kp$ is $1$-normal iff the club filter is normal iff $\kp$ is regular.

{\sc Proof:}
Suppose $\gm>\eta$ and $\kp$ is $\gm$-normal. Let $\la C_\alp:\alp<\kp\ra$ be a sequence of $\eta$-club subsets of $\kp$. Then $\diag_{\alp<\kp}C_\alp$ is $\eta$-stationary closed so if $\diag_{\alp<\kp}C_\alp$ is not $\eta$-club then it is not $\eta$-stationary. Let $\eta'<\eta$ and $C$ be a $\eta'$-club avoiding $\diag_{\alp<\kp} C_\alp$. Then setting $C'_\alp=C_\alp\cap C$ we have $C'_\alp$ is $\eta$-club and $\diag_{\alp<\kp} C'_\alp=\emp$. Define $f:\kp\rightarrow\kp$ by $f(\alp)$ is the least $\bt$ such that $\alp\notin C'_\bt$. Then as $\diag_{\alp<\kp} C'_\alp=\emp$, $f$ is regressive on $\kp$ and so by the $\gm$-normality of $\kp$ there is some $\bt<\kp$ such that $f^{-1}(\bt)$ is  a $\gm$-stationary set. But this contradict $C'_\bt$ being $\eta$-club.
\qed

Note that for a {limit} ordinal $\gm$, the domain of $f$ must be $\gm+1$-stationary, so this requirement is weaker than $\bigcup_{\eta<\gm} \Cf^\eta$ being normal (which is always false for limit $\gm$, see Proposition \ref{nolimitnormality}), but stronger than each $\Cf^\eta$ being normal (the latter can occur when $\kp$ is not $\gm$-reflecting).
The notion gets stronger as $\gm$-increases:
\begin{prop}\label{normalincreasing}
If $\gm>\eta$ and $\kp$ is $\gm$-normal then $\kp$ is $\eta$-normal. 
\end{prop}

{\sc Proof:}
Fix $\gm>\eta$ and suppose $\kp$ is $\gm$-normal. By Proposition \ref{normalnormal} we have $\Cf^\eta(\kp)$ is normal, and hence by Fodor's Lemma, for any $\eta+1$-stationary $S$ any regressive $f:S\rightarrow \kp$ is constant on an $\eta+1$-stationary, and hence an $\eta$-stationary, set.
\qed

\begin{prop}\label{pinormal}
If $\kp$ is $\Pi^1_\gm$-indescribable then 
\begin{enumerate}
\item $\kp$ is $\gm$-normal, and
\item $\{\lmb<\kp:\lmb \text{ is $\eta$-normal for all }\eta<\gm\}$ is $\gm+1$-stationary in $\kp$.
\end{enumerate}
\end{prop}

{\sc Proof:}
If $\kp$ is $\Pi^1_\gm$-indescribable then $\Cf^\gm(\kp)$ is normal by Lemma \ref{filtersV2}, so $\kp$ is $\gm$-normal. We have that being $\Sigma^1_\gm$-indescribable is $\Pi^1_\gm$ expressible so $$\{\lmb<\kp:\lmb \text{ is $\Sigma^1_\gm$-indescribable}\}$$ is in the $\Pi^1_\gm$ indescribability filter, and so by Lemma \ref{filtersV2} is $\gm+1$-stationary. As a $\Sigma^1_\gm$-indescribable cardinal is $\Pi^1_\eta$-indescribable for all $\eta<\gm$, by part (1),  $$\{\lmb<\kp:\lmb \text{ is $\Sigma^1_\gm$-indescribable}\}\sbs\{\lmb<\kp:\lmb \text{ is $\eta$-normal for all }\eta<\gm\}$$ and thus the latter is $\gm+1$-stationary.
\qed

We now state our main theorem.

\begin{thm}\label{da}
Let $\gm$ be an ordinal and $\kp>\gm$ a regular cardinal  such that 
\begin{enumerate}
\item for all $\eta<\gm$ $\kp$ is $\eta$-normal and
\item setting $$A_\gm=\{\alp<\kp:\text{\emph{for all $\eta$ with $1<\eta+1<\gm$,~$\alp$ is $\eta$-normal}}\}$$ we have that  $A_\gm$ is $\gm$-stationary in $\kp$. 
\end{enumerate}
Then $(\kp$ is $\gm$-stationary$)^L$ and hence $(\kp$ is $\Sigma^1_\gm$-indescribable$)^L$.
Furthermore, for any $S\sbs\kp$ such that $S\in L$ and $S\cap A_\gm$ is $\gm$-stationary we have $(S$ is a $\gm$-stationary subset of $\kp)^L$. Consequently, if $A_\gm$ is $\eta$-club for some $\eta<\gm$ then $\gm$-stationarity at $\kp$ is downward absolute to $L$.

\end{thm}

A few notes before we begin the proof. Assumption (2) is needed so that there are enough $\alp<\kp$ where we can apply the inductive hypothesis. For $\gm=\eta+2$, $A_\gm=\{\alp<\kp:\alp\text{ is $\eta$-normal}\}$ and if $\eta$ is a limit ordinal and $\gm=\eta$ or $\gm=\eta+1$ then $A_\gm=\{\alp<\kp:\text{for all $\eta'<\eta$~$\alp$ is $\eta'$-normal}\}$. For $\gm=2$ requirement (1) reduces to regularity and $(2)$ is vacuous. For $\gm=3$ requirement (1) is just the $1$-club filter is normal, and (2) that the regular cardinals are $3$-stationary - but this is a consequence of (1) as we shall see below. 

{\sc Proof:}
We prove this by induction, so suppose the theorem is true for any $\eta<\gm$ and $\kp$ satisfies the assumptions $(1)$ and $(2)$ of the theorem. First suppose $\gm$ is a limit ordinal . Then 
\begin{displaymath}
\begin{split}
A_\gm&=\{\alp<\kp:\text{\emph{for every $\eta+1<\gm$ ~ $\alp$ is $\eta$-normal}}\}\\
&=\{\alp<\kp:\text{\emph{for every $\eta<\gm$ ~ $\alp$ is $\eta$-normal}}\}\\
&=\bigcap_{\eta<\gm} A_\eta
\end{split}
\end{displaymath}
and so each $A_\eta$ is $\gm$-stationary and hence $\eta$-stationary. Thus we have $\kp$ satisfies, for all $\gm'<\gm$ $(1)$  for all $\eta<\gm'$ $\kp$ is $\eta$-normal and $(2)$ $A_{\gm'}$ is $\gm$-stationary.
Also, if $S\sbs\kp$ is such that $S\cap A$ is $\gm$-stationary, then for any $\eta<\gm$, $S\cap A_\eta$ is $\gm$- (and hence $\eta$-) stationary. Therefore for any such $S$ the inductive hypothesis gives that, for each $\eta<\gm$, $(S$ is $\eta$-stationary$)^L$. But then $(S$ is $\gm$-stationary$)^L$ by definition.

Now suppose $\gm=\eta+1$ and $\eta$ is a limit ordinal. By $(1)$, $\kp$ is $\eta$-normal. Suppose we have $S\sbs\kp$ with $S\in L$ such that $S\cap A_\gm$ is $\gm$-stationary. To show that $S$ is $\eta+1$-stationary in $L$, let $B\sbs\kp$ with $B\in L$ such that  for every $\alp\in S$
$$L\vDash B\cap\alp\text{ is not $\eta$-stationary}.$$
For each $\alp\in S$ set $\eta_\alp$ to be the least ordinal such that $B\cap\alp$ is not $\eta_\alp$ stationary. As $\eta$ is a limit ordinal, $\eta_\alp<\eta$ for all $\alp\in S\cap A_\gm$. Therefore, by the $\eta$-normality of $\kp$ there is some $\dlt<\eta$ such that setting $X_\dlt=\{\alp\in S:\eta_\alp=\dlt\}$ we have $X_\dlt\cap A_\gm$  is $\eta+1$-stationary. Fix such a $\dlt$, and note $X_\dlt\in L$. 

As $\dlt<\eta$, by the inductive hypothesis $\kp$ is $\dlt$-reflecting in $L$ and $X$ is $\dlt+1$-stationary in $L$. Now working in $L$, suppose $B$ were $\dlt$-stationary. Then $d_\dlt(B)$ would be $\dlt$-club below $\kp$, so we would have some $\alp\in d_\dlt(B)\cap X$. But this is a contradiction as for any $\alp\in X$, $B\cap\alp$ was not $\dlt$-stationary. So we have $B$ is not $\dlt$-stationary in $L$. Thus
 $$L\vDash B\text{ is not $\eta$-stationary in }\kp$$
and so $S$ is $\eta+1=\gm$-stationary in $L$.

Finally, if $\gm$ is not a limit ordinal or a successor of a limit ordinal, we have more work to do. The proof will proceed roughly as Magidor's proof for the case $\gm=2$ ($\gm=1$ is the downward absoluteness of stationarity, which is obvious). There will be several points of departure from Magidor's proof, as we shall have to take into account the non-absoluteness of $\eta$-clubs and $\eta$-stationary sets. To simplify the presentation we shall assume that $\gm>2$.

Let $\gm=\eta+2$ with $\eta\geq 1$. Let $\kp$ be regular and satisfy $(1)$ and $(2)$ of the theorem and $S\sbs\kp$ be such that $S\cap A_\gm$ is $\gm$-stationary with $S\in L$. Then $(1)$ gives us that $\Cf^{\eta}(\kp)$ is normal. 


\begin{claim}
The regular cardinals are $1$-club below $\kp$.
\end{claim} 
{\sc Proof:}
As $\eta\geq1$ and we have assumed the $\eta$-club filter on $\kp$ is normal, $\Cf^1(\kp)$ must be normal. Suppose the singulars were $2$-stationary below $\kp$. Let $f:Sing\cap\kp\rightarrow\kp$ with $f(\alp)=cof(\alp)$. Then $f$ is regressive so the normality of $\Cf^1(\kp)$ would give a $2$-stationary set $A$ of ordinals all having the same cofinality. But this is impossible, as for any $\alp\in A$ taking a club of order type $cf(\alp)$ its limit points would all have smaller cofinality, so $A\cap\alp$ could not be stationary.
\qed

We now show that we can assume for each $\alp$ in $S$, we have that in $L$, $\alp$ is regular and $\eta$-reflecting.

Now let $S'=\{\alp\in S:(\alp$ is regular and $\eta$-reflecting$)^L\}$. Clearly $S'\in L$. So we just need to show $S'\cap A_\gm$ is $\gm$-stationary in $V$. We have $S\cap A_\gm$ is $\eta+2$-stationary and as $\kp$ is $\eta+1$-reflecting $d_{\eta+1}(A_\gm)$ is $\eta+1$-club. Also, by the above claim, the regulars are $1$-club, so setting $$S^*=S\cap A_\gm\cap d_{\eta+1}(A_\gm)\cap Reg$$ we have $S^*$ is $\eta+2=\gm$-stationary. We show $S^*\sbs S'$, so let $\alp\in S^*$. Then $\alp$ is regular and $\eta$-normal.  Also as $\alp\in d_{\eta+1}(A_\gm)$ we have $A_\gm$, and hence $A_\eta$, is $\eta+1$-stationary in $\alp$. Thus $\alp$ satisfies the assumptions of the theorem for $\eta+1$ and so by the inductive hypothesis $\alp$ is $\eta$-reflecting and regular in $L$.

From now on we assume $S=S'$, so for all $\alp\in S$, $(\alp$ is regular and $\eta$-reflecting$)^L$, and $S^*=S\cap A_\gm\cap d_{\eta+1}(A_\gm)\cap Reg$ is a $\gm$-stationary subset of $S$. Furthermore, for any $\alp\in S^*$ we have by the inductive hypothesis that for any $T\sbs \alp$ such that $T\cap A_\gm$ is $\eta+1$-stationary, $(T$ is $\eta+1$-stationary$)^L$. 

We now want to show that $S$ is $\eta+2$-stationary in $L$, so we let $B\sbs\kp$ with $B\in L$ such that  for every $\alp\in S$
$$L\vDash B\cap\alp\text{ is not $\eta+1$-stationary}.$$
we shall show $$L\vDash B\text{ is not $\eta+1$-stationary in }\kp.$$

Now working in $L$: For $\alp\in S$  we have  $\alp$ is $\eta$-reflecting and $B$ is not $\eta+1$-stationary, so we can find an $\eta$-club $D\sbs\alp$ which avoids $B$. Let $D_\alp$ be the minimal such set in the canonical well-ordering of $L$. Now we can take $\nu_\alp$ to be minimal such that $B\cap\alp, D_\alp\in L_{\nu_\alp}$ and $L_{\nu_\alp}$ is $\Pi^1_\eta$-correct for $\alp$. (In Magidor's proof $\eta$-stationary correctness was not required, because there it was $0$-stationary correct; and a set being unbounded is absolute for transitive models.) It is clear that $|\nu_\alp|=|\alp|$. Also, if  $\alp$ is regular then any $L_{\bt}$ for $\bt>\alp$ will be correct about $d_{\eta}$ below $\alp$, and thus $D_\alp$ can be uniformly defined within any level of $L$ which is $\Pi^1_\eta$-correct for $\alp$, contains $B\cap\alp$ and sees that $B\cap\alp$ is not $\eta+1$-stationary.

Set $M_\alp=\la L_{\nu_\alp},\in,\alp,B\cap\alp,D_\alp\ra.$ 
Now following Magidor's proof, as the structure $M_\alp$ is no bigger than $\alp$ in $L$ it is isomorphic to a structure $N_\alp=\la\alp, E_\alp, \mu_\alp, B'_\alp, D'_\alp\ra$. We take $N_\alp$ minimal and let $h_\alp $ be the inverse collapse, $h_\alp:M_\alp\rightarrow N_\alp$. Set $f_\alp=h_\alp\restr \alp$.
Now we go back to working in $V$.
\begin{lem}
There is an $\eta$-club $G\sbs\kp$ such that for any $\alp,\bt\in G\cap S^*$ if $\alp<\bt$ then $\la N_\alp, f_\alp\ra\prec\la N_\bt,f_\bt\ra$.
\end{lem}

{\sc Proof:}
Let $\La_N$ be the language $\La_\in$ with added constant symbols $\mu, B, D$ and function symbol $f$. First we show that for any particular formula $\varphi\in\La_N$ with $k$ free variables and any $\alp_1,\dots,\alp_k$ ordinal parameters from $\kp$ we have (where $\mu, B, D$ and $f$ are interpreted  in the obvious way as $\mu_\bt, B'_\bt,D'_\bt$ and $f_\bt$) that either  $$X_{\varphi(\vec{\alp})}=\{\bt<\kp:\la N_\bt,f_\bt\ra\vDash\varphi(\vec{\alp}) \text{ or } \bt\notin S^*\}$$ 
or $$X_{\neg\varphi(\vec{x})}=\{\bt<\kp: \la N_\bt,f_\bt\ra\vDash\neg\varphi(\vec{\alp}) \text{ or } \bt\notin S^*\}$$ contains an $\eta$-club.  The lemma will then follow by taking the diagonal intersection over parameters and formulae. 

Fix $\varphi$ etc. Working in V, assume for contradiction that setting $$A_1:=\{\bt\in S: \alp_1,\dots,\alp_k<\bt, N_\bt\vDash\varphi(\vec{\alp})\}$$ $$A_2:=\{\bt\in S: \alp_1,\dots,\alp_k<\bt, N_\bt\vDash\neg\varphi(\vec{\alp})\}$$ we have $A_1\cap S^*$ and $A_2\cap S^*$ are both $\eta+1$-stationary, and $A_1,A_2\in L$. As  $S^*$ is $\eta+2$-stationary, we can find a regular $\alp\in S^*$ such that $A_1\cap S^*$ and $A_2\cap S^*$  are both $\eta+1$-stationary in $\alp$. Thus we have $\alp$ is satisfies the assumptions $(1)$ and $(2)$ for $\eta+1$, and so by the inductive hypothesis $\alp$ is $\Pi^1_{\eta}$-indescribable in $L$, and as $A_{\eta+1}\sps A_\gm\sps S^*$ we have $A_1\cap\alp, A_2\cap\alp$ are also $\eta+1$-stationary in $L$. Without loss of generality suppose $ N_\alp\vDash\varphi(\vec{\alp})$. We now work in $L$ and show that on an $\eta$-club below $\alp$ we have $N_\bt\vDash\varphi(\vec{\alp})$ and hence we have a contradiction with $A_2$ being $\eta+1$-stationary.

So, let $\rho>\nu_\alp$ be an admissible ordinal such that $L_\rho$ is $\Pi^1_\eta$-correct for $\alp$. Now $\alp$ is $\Pi^1_{\eta}$-indescribable in $L$ so by Lemma \ref{traceclub3}, we have 
$$\trace^{\Pi^1_\eta}(L_\rho,\{B\cap\alp, D\cap\alp,\alp_1,\dots,\alp_k\},\alp)$$
is $\eta$-club below $\alp$. Let $\bt\in\trace^{\Pi^1_\eta}(L_\rho,\{B\cap\alp, D\cap\alp,\alp_1,\dots,\alp_k\},\alp)$, let $\dlt$ be such that $L_\dlt\cong L_\rho\{\bt\cup\{\alp,B\cap\alp,D\cap\alp,\alp_1,\dots,\alp_k\}\}$, and $\pi$ be the collapsing map. It is clear that $\pi(B\cap\alp)=B\cap\bt$ and $\dlt$ is admissible. As the $D_\bt$'s were uniformly definable from $B\cap\bt$ at $\Pi^1_\eta$-correct levels of $L$, we have $\pi(D_\alp)=D_\bt$.

As $\rho>\nu_\alp$, we have 
$$L_\rho\vDash\exists\gm ~L_\gm\text{ is $\Pi^1_\eta$-correct for $\alp$ and }B\cap\alp, D_\alp\in L_\gm.$$
Thus we have $$L_\dlt\vDash\exists\gm~ L_\gm\text{ is $\Pi^1_\eta$-correct for $\bt$ and }B\cap\bt, D_\bt\in L_\gm$$ and because $L_\dlt$ is $\Pi^1_\eta$-correct this statement hold in $L$, so we have $\nu_\bt<\dlt$. Thus $\pi(N_\alp)=N_\bt$, and so $N_\bt\vDash\varphi(\vec{\alp})$.
Hence, on an $\eta$-club (the $\Pi^1_\eta$ trace) below $\alp$ we have $N_\bt\vDash\varphi(\vec{\alp})$ and hence we have a contradiction with $A_2$ being $\eta+1$-stationary in $\alp$. Thus we must have either $A_1\cap S^*$ is not $\eta+1$-stationary in $\kp$ or $A_2\cap S^*$ is not $\eta+1$-stationary in $\kp$. But this means either $X_{\varphi(\vec{\alp})}$ or $X_{\neg\varphi(\vec{\alp})}$ contains an $\eta$ club, which is what we wanted to show.

Now we conclude by taking the diagonal intersection. Let $\la\phi_n:n\in\om\ra$ list all the formulae in the language $\La_\N$. Let $\la \vec{\alp}_\dlt:\dlt<\kp\wedge lim(\dlt)\ra$ enumerate $^{<\om}\kp$ in order type $\kp$. On a club of $C$ we have $\la\vec{\alp}_\dlt:\dlt<\alp\ra$ lists all of $^{<\om}\alp$. Then for limit $\dlt<\kp$ we can set $C_{\dlt+n}=\kp$ if $\phi$ does not have $lh(\vec{\alp}_\dlt)$ free variables. Otherwise we set $C_{\dlt+n}=X_{\varphi_n(\vec{\alp_\dlt})}$ if this contains an $\eta$-club and if not $C_{\dlt+n}=X_{\neg\varphi_n(\vec{\alp_\dlt})}$, which must then contain an $\eta$-club. Then $C\cap\diag_{\bt<\kp} C_\bt$ contains an $\eta$-club, and for any $\alp<\bt$ with $\alp,\bt\in S^*\cap C\cap \diag_{\bt<\kp} C_\bt$, and for any $\phi\in \La_N$ and $\vec{\alp}\in {^{<\om}\alp}$ we have that $$\la N_\alp,f_\alp\ra\vDash\varphi(\vec{\alp})\text{ iff }\la N_\bt,f_\bt\ra\vDash\varphi(\vec{\alp})$$ and thus $\la N_\alp, f_\alp\ra\prec\la N_\bt,f_\bt\ra$. (Although each $X_{\varphi(\vec{x})}\in L$, as containing a $\eta$-club is not absolute between $V$ and $L$, the appropriate sequence of $X_{\varphi(\vec{x})}$ to take the diagonal intersection of and obtain $G$ need not be in $L$, and hence we cannot guarantee $G\in L$ at this stage.)
\qed

To finish the proof we want to show $B$ is not $\eta+1$-stationary in $L$, so we shall produce an $\eta$-club set $D$ with $D\cap B=\emp$. For each $\alp<\bt\in S^*\cap G$ we have $M_\alp\prec M_\bt$. Let $M$ be the direct limit  $\varinjlim \la M_\alp\ra_{\alp\in G\cap S^*}$. As each $M_\alp\vDash V=L$ we have $M\vDash V=L$, and thus $M$ is just some $L_\rho$. Also setting  $D=\bigcup_{\alp\in S^*\cap G}D_\alp$, it is clear that $M=\la L_\rho,\in,\kp, B,D\ra $ (remember $B=\bigcup_{\alp\in S}B_\alp$). we shall show that in $L$, $d_{\eta}(D)\cap B=\emp$ before checking that $D$ is stationary. 

It is clear that for any $\alp<\bt$ from $G\cap S^*$ we have $D_\bt$ is an end extension of $D_\alp$. So if $D$ is $\eta$-stationary below some $\alp<\kp$ then for some (any) $\bt>\alp$ with $\bt\in S^*\cap G$ we have $D_\bt\cap\alp=D\cap\alp$ so by definition of $D_\bt$ we have $\alp\notin B$ as required. It remains to show that $D$ is $\eta$-stationary. We want to use the fact that an $\eta$-stationary union of $\eta$-stationary sets is $\eta$-stationary - but to apply this we need to find an $\eta$-stationary set $H\in L$ with $D_\alp=D\cap\alp$ for each $\alp\in H$. (Note we cannot do this in $V$ and use the inductive hypothesis to show $D$ is $\eta$-stationary in $L$, because the $D_\alp$'s need not be $\eta$-stationary in $V$.) The obvious candidate for this is $G\cap S^*$, but as noted above, we needn't have $G\cap S^*\in L$. 

\begin{claim}
There is $H\in L$ with $ G\cap S^*\sbs H$ and for each $\alp\in H$, $M_\alp\prec M$.
\end{claim}

{\sc Proof:}

Let $H=S\cap\trace^{\Pi^1_\eta}(L_\rho,\{B,D\},\kp)$. Note that we do not assume $L_\rho$ is $\Pi^1_\eta$-correct, so we do not yet know that $H\neq\emp$.

First we show that for $\alp$ in $H$ we have $L_\rho\{\alp,B,D\}\cong L_{\nu_\alp}$ with $\pi(B)=B\cap\alp$ and $\pi(D)=D_\alp$ and hence $M_\alp\prec M$. Fix $\alp\in H$ and let $\dlt$ be such that $L_\dlt\cong L_\rho\{\alp\cup\{B,D\}\}$, and $\pi$ be the collapsing map. Clearly $\pi(B)=B\cap\alp$ and as $D$ was defined from $B$ in the same way that $D_\alp$ was defined from $B\cap\alp$, and $L_\dlt$ is $\Pi^1_\eta$-correct and sees that $B\cap\alp$ is not $\eta+1$-stationary, we must have $\pi(D)=D_\alp$. Thus to show $\dlt=\nu_\alp$ we only need to verify minimality of $L_\dlt$. Let $\Theta(\bt)$ be the following statement, where $\psi_\gm$ is the universal $\Pi^1_\gm$ formula. 


\begin{displaymath}
\begin{split}
\Theta(\bt): \quad\quad\forall \gm>\bt\big[&\big(B\cap\bt \notin L_\gm\big)\vee\big( D_\bt\notin L_\gm\big)\\
& \vee \exists U\notin L_\gm~\exists n\in\om~\exists \eta'<\eta~\big(U\text{ is minimal with }\psi_{\eta'}(\bt,n, U)\big)\big].
\end{split}
\end{displaymath}

Then if $L_\nu$ is $\Pi^1_\eta$-for $\bt$ and $L_\nu\vDash \Theta(\bt)$ we must have $L_\nu$ is the minimal $\Pi^1_\eta$-correct level containing $B\cap\bt$ and $D_\bt$.

Now for any $\bt\in G\cap S^*$
$$L_{\nu_\alp}\vDash \Theta(\bt)$$
so by elementarity $$L_\rho\vDash\Theta(\kp)$$ and again by elementarity $$L_\dlt\vDash\Theta(\dlt).$$ 
Thus $\dlt$ is minimal so $\dlt=\nu_\alp$ and $M_\alp\prec M$.

Now assume $\alp\in S^*\cap G$. We show $L_{\nu_\alp}\cong L_\rho\{\alp\cup\{B,D,\kp\}\}$ and so as $L_{\nu_\alp}$ is $\Pi^1_\eta$-correct $\alp\in H$. Let $j:M_\alp\prec M$ 
By Lemma \ref{tracecont} $L_{\nu_\alp}\{\{B\cap\alp,D_\alp\}\cup\alp+1\}=L_{\nu_\alp}$ and thus $$L_{\nu_\alp}\cong L_\rho\{j``L_{\nu_\alp}\}=L_\rho\{\{j(B\cap\alp),j(D_\alp)\}\cup j``\alp+1\}=L_\rho\{\{B,D,\kp\}\cup\alp\}$$ and we're done.
\qed

Now as $H\sps S^*\cap G$ we have in $V$ that $H$ is $\eta$-stationary in $\kp$, and hence by the induction $L\vDash H$ is $\eta$-stationary in $\kp$. Also, for every $\alp\in H$ we have $M_\alp\prec M$ and hence $D_\alp=D\cap\alp$ and $D_\alp$ is $\eta$-stationary. Then in $L$, $D$ is the $\eta$ stationary union of $\eta$-stationary sets, so $L\vDash D$ is $\eta$-stationary in $\kp$. As $d_{\eta}(D)$ avoids $B$ we have
$$L\vDash B\text{ is not $\eta+1$-stationary in }\kp$$
and so we're done.
\qed

This gives us the following equiconsistency result, showing that assumption on the $\gm$-club filters can be as strong as $\Pi^1_\gm$-indescribability.

\begin{col}
Fix an ordinal $\gm$. It is consistent that there is a regular cardinal $\kp>\gm$ such that 
\begin{enumerate}
\item for all $\eta<\gm+1$, $\kp$ is $\eta$-normal and
\item $\{\alp<\kp:\text{\emph{for all $\eta$ with $0<\eta<\gm$,~$\alp$ is $\eta$-normal}}\}$ is $\gm+1$-stationary in $\kp$. 
\end{enumerate}
if and only if it is consistent that there is a $\Pi^1_\gm$-indescribable cardinal.
\end{col}

{\sc Proof:}
Assuming $\kp$ satisfying (1) and (2) as above, Theorem \ref{da} tells us that in $L$, $\kp$ is $\Sigma^1_{\gm+1}$-indescribable and hence $\Pi^1_\gm$-indescribable. Thus if the existence of a regular cardinal satisfying $(1)$ and $(2)$ is consistent with $ZFC$ so is the existence of a $\Pi^1_\gm$-indescribable cardinal. Proposition \ref{pinormal} gives the converse.
\qed

If we make a stronger assumption on $\kp$ we can get the following much simpler (though weaker) statement of downward absoluteness:

\begin{col}\label{dasimple}
Assume that for any ordinal $\gm$ and any regular $\gm$-reflecting cardinal $\kp$, $\Cf^\gm(\kp)$ is normal. Then if $\kp$ is regular and $S\sbs\kp$ is $\eta$-stationary with $S\in L$ we have $(S\text{ is $\eta$-stationary in }\kp)^L$.
\end{col}

We shall see in the next section that we need the stronger statement to get the downward absoluteness of $\gm$-ineffability.

\section{Some applications: Ineffability and $\diamondsuit$ Principles} 

\label{Chapter6} 

In this chapter we give generalisations of ineffability and $\diamondsuit$ principles using $\gm$-stationary sets. 
As we shall see, many of the old results follow through in this context. We shall also detail the relations between the different levels of these generalised principles. 


\subsection{$\gm$-Ineffables}

We start by defining a new, natural generalisation of ineffability, and exploring its basic properties.

\begin{defn}
A regular, uncountable cardinal $\kp$ is $\gm$-ineffable for $\gm<\kp$ iff, whenever $f:[\kp]^2\rightarrow 2$, there is a $\gm$-stationary set $X\subseteq \kp$ such that $|f``[X]^2|=1$.
\end{defn}

\begin{rem} It is clear that $\gm$-ineffability implies $\bt$-ineffability if $\bt<\gm$. Note that $1$-ineffable reduces to the ordinary definition of ineffability, and $0$-ineffable is weakly compact. 
\end{rem}

From now on we shall assume $\gm\geq 1$. The following theorem gives a useful characterisation of $\gm$-ineffability, well-known for $\gamma =1$.

\begin{thm}\label{ineff}
Let $\kp$ be a regular, uncountable cardinal. Then $\kp$ is $\gm$-ineffable iff whenever $\la A_\alp:\alp<\kp\ra$ is such that for all $\alp$, $A_\alp\subseteq\alp$, there is a set $A\subseteq\kp$ such that $\{\alp<\kp:A_\alp=A\cap\alp\}$ is $\gm$-stationary in $\kp$.
\end{thm}

{\sc Proof:}
($\Rightarrow$) Suppose $\kp$ is $\gm$-ineffable and let  $\la A_\alp:\alp<\kp\ra$ be as above. Define a function $h:[\kp]^2\rightarrow 2 $ by $h(\{\alp,\bt\})=0$ iff, assuming $\alp<\bt$, there is some $\dlt$ such that $A_\alp\cap \dlt\subsetneq A_\bt\cap\dlt$. By $\gm$-ineffibility, let $X$ be $\gm$-stationary such that $|h``[X]^2|=1$. Suppose first $h``[X]^2=\{0\}$. Then for $\alp,\bt\in X$ with  $\alp<\bt$ we have $A_\alp\cap\dlt\subsetneq A_\bt\cap\dlt$ for some $\dlt$.

For each $\nu<\kp$, let $\alp_\nu$ be least in $X$ such that $\alp_\nu\geq\nu$ and for all $\bt\in X$ with $\bt>\alp_\nu$, we have $A_\bt\cap\nu=A_{\alp_\nu}\cap\nu$.
This is possible as if $\la\bt_i:i<\kp\ra$ enumerates $X\setminus\nu$ we have $sup\{\alp<\kp:A_\nu\cap\alp\subseteq A_{\bt_i}\cap\alp\}$ is strictly increasing, so will pass $\nu$. Then after this $ A_{\bt_i}\cap\nu$ is subset-increasing, so as $\kp$ is regular it must eventually be constant.

Let $ C=\{\dlt\in\kp:\forall \nu<\dlt \,\, \alp_\nu<\dlt\}$, which is closed unbounded. Thus $Y=X\cap C\cap LimOrd$ is $\gm$-stationary in $\kp$.  Now for $lim(\nu)$ we have $\alp_\nu$ is the least $\alp\in X$ such that $\alp \geq sup_{\eta<\nu}\alp_\eta$, so if $\nu\in Y$ then $\alp_\nu=\nu$. Hence for any $\alp\in Y$ and $\bt\in X$ with $\bt>\alp$ we have $A_\bt\cap\alp= A_\alp$. Thus setting $A= \bigcup_{\alp\in Y}A_\alp$ we have $Y\subseteq \{\alp\in\kp: A\cap\alp=A_\alp\}$ so we're done.

Now suppose $h``[X]^2=\{1\}$. Then for $\alp<\bt$ in $X$ we have $A_\alp=A_\bt$ or $A_\alp\cap\dlt\supsetneq A_\bt\cap\dlt$ for some $\dlt$. For each $\nu<\kp$ we can define $\alp_\nu$ exactly as before, for similar reasons, and the construction goes through in the same way.

($\Leftarrow$) Let $f:[\kp]^2\rightarrow 2$ be given. For $\alp<\kp$ set $$A_\alp=\{\bt<\alp:f(\{\bt,\alp\})=1\}$$ Then by assumption there is $A\subseteq\kp$ such that $$X=\{\alp<\kp:A\cap\alp=A_\alp\}$$ is $\gm$-stationary. It is easy to see that $f``[A\cap X]^2=\{1\}$ and $f``[X\setminus A]^2=\{0\}$. But $X=(X\cap A) \cup (X\setminus A)$, so one of these must be $\gm$-stationary and we're done.
\qed

We cannot fully generalise the implication from ineffability to $\Pi^1_2$-indescribability, the reason being that the $\gm$-club filter need not coincide with the $\Pi^1_\gm$-indescribability filter outside of $L$. We do however have that any $\gm$-ineffable is $\gm+1$-reflecting (Theorem \ref{ineffablereflecting} below).

\begin{thm}
If $V=L$ or $\gm\in \{1,2\}$, and $\kp$ is $\gm$-ineffable then $\kp$ is $\Pi^1_{\gm+1}$-indescribable.
\end{thm}

The case for $\gm=1$ is given in \cite{De} VII.2.2.3. We do an induction essentially following that proof.

{\sc Proof:}
Assume the theorem holds for any $\eta<\gm$. Suppose for a contradiction that $\phi$ is $\Delta_0$, $A\sbs\kp$ and 
$$\forall X\sbs\kp ~\Sigma \text{ wins } G_\gm(\kp,\phi,\la A,X\ra)$$
but for any $\alp<\kp$
$$\exists X\sbs\alp ~\Pi \text{ wins } G_\gm(\alp,\phi,\la A\cap\alp,X\ra)$$
For each $\alp<\kp$ fix some $X_\alp\sbs\alp$ such that $\Pi \text{ wins } G_\gm(\alp,\phi,\la A,X_\alp\ra)$. By $\gm$-ineffability of $\kp$, take $X\sbs\kp$ such that $S:=\{\alp<\kp:X_\alp=X\cap\alp\}$ is $\gm$-stationary. Now we have $\Sigma \text{ wins } G_\gm(\kp,\phi,\la A,X\ra)$, so we can fix $\eta<\gm$ and $Y\sbs\kp$ such that $\Pi$ wins $G_\eta(\kp,\phi',\la A,X,Y\ra)$. Now inductively we have that $\kp$ is $\Pi^1_\eta$-indescribable, so if $V=L$ or $\eta=0,1$, by Theorem \ref{filtersL1} and Lemma \ref{filtersL2} this statement reflects to an $\eta$-club $C$. Let $\alp\in C\cap S$. Then as $\alp\in S$ we have $\Pi \text{ wins } G_\gm(\alp,\phi,\la A,X_\alp\ra)$, \ie for any $\gm'<\gm$ and any $Y'\sbs\alp$, $\Sigma$ wins $G_{\gm'}(\alp,\phi',\la A,X\cap\alp,Y'\ra)$. But $\alp\in C$ so $\Pi$ must win $G_\eta(\alp,\phi',\la A\cap\alp,X\cap\alp,Y\cap\alp\ra)$ so we have a contradiction.
\qed

Outside of $L$, we can still obtain the following:

\begin{thm}\label{ineffablereflecting}
If $\kp$ is $\gm$-ineffable, then $\kp$ is $\gm+1$-reflecting.
\end{thm}

{\sc Proof:}
Suppose $\kp$ is $\gm$-ineffable, $S\sbs\kp$ with $S$ not $\gm+1$-stationary in any $\alp<\kp$. We show $S$ is not $\gm+1$-stationary. For each $\alp\in \kp$ take $C_\alp\sbs\alp$ to be $\gm_\alp$-club avoiding $S$, for some $\gm_\alp\leq\gm$. By $\gm$-ineffability, let $C$ and $\gm'$ be such that $X=\{\alp<\kp:C_\alp=C\cap\alp, \gm'=\gm_\alp\}$ is $\gm$-stationary. We claim $C$ is $\gm'$-club and avoids $S$.
$C$ is $\gm'$-stationary as $X\sbs C$: to see this, let $\alp<\bt\in X$. Then $C_\alp$ is $\gm'$-stationary in $\alp$ by definition of $C_\alp$. Now as $\alp,\bt\in X$, $C_\bt\cap\alp=C_\alp$, and as $C_\bt$ is $\gm'$ stationary closed, thus $\alp\in C_\bt$ and so $\alp\in C$. Also, $C$ is $\gm'$-stationary closed as each $C_\alp$ with $\alp\in X$ is. So $C$ is $\gm'$ club and avoids $S$, and we're done. \qed


In order to prove that $\gm$-ineffability is downwards absolute to $L$, we now look in more detail at the theory of $\gm$-ineffables. The following shows that if $\kp$ is $\gm+1$-ineffable then the $\gm$-club filter on $\kp$ is normal. 


\begin{lem}\label{ineffclubnormal}
If $\kp$ is $\gm$-ineffable then $\kp$ is $\gm$-normal.
\end{lem}

{\sc Proof:}
Let $\kp$ be $\gm$-ineffable. First we show that any $f:\kp\rightarrow\kp$ which is regressive is constant on a $\gm$-stationary set. So let $f$ be such a function and for each $\alp<\kp$ set $A_\alp=\{f(\alp)\}$. By ineffability, there is some $A\sbs\kp$ such that $$X=\{\alp<\kp:A\cap\alp=\{f(\alp)\}\}$$ is $\gm$-stationary. But then for $\alp,\bt\in X$ we must have $f(\alp)=f(\bt)$. 

Now suppose for a contradiction that $S$ is $\gm+1$-stationary and $f:S\rightarrow\kp$ is regressive such that for any $\alp<\kp$, $f^{-1}(\alp)$ is not $\gm$-stationary. For each $\alp<\kp$ take $\eta_\alp<\gm$ and $C_\alp\sbs\kp$ such that $C_\alp$ is $\eta_\alp$-club and avoids $f^{-1}(\alp)$. Then as each $C_\alp$ is $\gm$-stationary closed, $\diag_{\alp<\kp} C_\alp$ is $\gm$-closed. Also, as each $C_\alp$ avoids $f^{-1}(\alp)$ and $f$ is regressive, $\diag_{\alp<\kp} C_\alp$ must avoid $S$. Thus we must have that $\diag_{\alp<\kp} C_\alp$ is not $\gm$-stationary. Let $\eta<\gm$ and $C$ be a $\eta$-club avoiding $\diag_{\alp<\kp} C_\alp$. Then setting $C'_\alp=C_\alp\cap C$ we have $C'_\alp$ is $max\{\eta,\eta_\alp\}$-club and $\diag_{\alp<\kp} C'_\alp=\emp$. Let $f:\kp\rightarrow\kp$ be defined by $f(\alp)$ is the least $\bt$ such that $\alp\notin C'_\bt$. As each $\alp\notin \diag_{\alp<\kp} C'_\alp$ we have that $f$ is well defined and regressive. By above, $f$ is constant on a $\gm$-stationary set $X$, so there is $\dlt<\kp$ such that for each $\alp\in X$, $\alp\notin C'_\dlt$. But this contradicts $C'_\dlt$ being $\eta_\dlt$-club.
\qed

\begin{defn}
We say $S\sbs\kp$ is $\gm${\em -ineffable in } $\kp$ iff whenever $f:[S]^2\rightarrow 2$, there is $X\subseteq S$ such that $|f``[X]^2|=1$ and  $X$ is $\gm$-stationary in $\kp$.
\end{defn}

\begin{rem}
It is easy to see that if some $S\sbs\kp$ is $\gm$-ineffable in $\kp$ then $\kp$ is $\gm$-ineffable.
\end{rem}

It is clear that to be $\gm$-ineffable a subset of $\kp$ must be $\gm$-stationary, but we shall see that such a set must in fact be $\gm+2$-stationary. There are always proper subsets of $\kp$ which are $\gm$-ineffable:

\begin{prop}
If $\kp$ is $\gm$-ineffable and $C$ is $\eta$-club for some $\eta<\gm$ then $C$ is $\gm$-ineffable.
\end{prop}

{\sc Proof:}
Extend $f:[C]^2\rightarrow 2$ to $g:[\kp]\rightarrow 2$ arbitrarily. By $\gm$-ineffability of $\kp$ there is a $\gm$-stationary set $X\sbs\kp$ such that $g$ is constant on $[X]^2$. As $X$ is $\gm$-stationary, $X\cap C$ is also $\gm$ stationary, and as $g$ extended $f$, we must have $f$ is constant on $X\cap C$.
\qed

We also have the analogue of Theorem \ref{ineff} for $\gm$-ineffable subsets:

\begin{thm}
Let $\kp$ be a regular, uncountable cardinal. Then $S\sbs \kp$ is $\gm$-ineffable iff whenever $\la A_\alp:\alp\in S\ra$ is such that for all $\alp\in S$, $A_\alp\subseteq\alp$, there is a set $A\subseteq\kp$ such that $\{\alp\in S:A_\alp=A\cap\alp\}$ is $\gm$-stationary in $\kp$.
\end{thm}

{\sc Proof:}
The proof of Theorem \ref{ineff} works in exactly the same way relativised to $S$.\\
\qed


\begin{prop}
If $S\sbs\kp$ is $\gm$-ineffable then $S$ is $\gm+2$-stationary.
\end{prop}

{\sc Proof:}
This is essentially the same argument as Theorem \ref{ineffablereflecting}. Suppose $S$ is $\gm$-ineffable but not $\gm+2$-stationary. Let $T\sbs\kp$ be $\gm+1$-stationary with $d_{\gm+1}(T)\cap S=\emp$. Then for each $\alp\in S$, $T\cap\alp$ is not $\gm+1$-stationary in $\alp$ so we can find $C_\alp\sbs\alp$ which is $\gm$-club and avoids $T\cap\alp$. Now, using the $\gm$-ineffability of $S$ we can find $C\sbs\kp$ such that $\{\alp\in S:C\cap\alp=C_\alp\}$ is $\gm$-stationary. Then $C$ is $\gm$-club and avoids $T$ - but this contradicts $T$ being $\gm+1$-stationary.
\qed

\begin{lem}\label{ineffbelow}
If $\kp$ is $\gm$-ineffable then the set $$E_\gm=\{\alp<\kp:\alp\text{ is $\eta$-ineffable for every $\eta<\gm$}\}$$ is $\gm$-ineffable. 
\end{lem}

{\sc Proof:}
Let $\kp$ be $\gm$-ineffable. Suppose $\la B_\alp:\alp\in E_\gm\ra$ is such that $B_\alp\sbs\alp$ for each $\alp\in E_\gm$. We show that there is a $B\sbs\kp$ such that $\{\alp\in E_\gm:B\cap\alp=B_\alp\}$ is $\gm$-stationary. For $\alp\in E_\gm$ set $A_\alp=\{\bt+1: \bt\in B_\alp\}$.
For each $\alp\notin E_\gm$ let $\eta_\alp<\gm$ be such that $\alp$ is not $\eta_\alp$-ineffable and let $B_\alp=\la B^\bt_\alp\sbs\bt:\bt<\alp\ra$ be a sequence witnessing that $\alp$ is not $\eta_\alp$-ineffable. Let $A_\alp\sbs\alp$ code the sequence $B_\alp$ and $\eta_\alp$ such that $0\in A_\alp$.

By $\gm$ ineffability we have a $A\sbs\kp$ such that $$X=\{\alp<\kp:A\cap\alp=A_\alp\}$$ is $\gm$-stationary. Now if $0\notin A$ we have $X\sbs E_\gm$ and so we're done - setting $B=\{\alp<\kp:\alp+1\in A\}$ we have for limit $\alp$ that $A\cap\alp=A_\alp$ iff $B\cap\alp=B_\alp$, so $B$ is as required. If $0\in A$ then $X\cap E_\gm=\emp$. We show that this leads to a contradiction. 

So suppose each $\alp\in X$ is not $\eta$-ineffable for some $\eta<\gm$. Then $A$ codes $B=\la B^\bt\sbs\bt:\bt<\kp\ra$ and $\eta$ such that for each $\alp\in X$, $B\restr \alp=B_\alp$ and $\eta_\alp=\eta$.. Using the $\gm$-ineffability of $\kp$ again, let $C\sbs\kp$ be such that $$Y=\{\alp<\kp:C\cap\alp=B^\alp\}$$ is $\gm$-stationary.  As $\kp$ is $\eta$-reflecting $d_\eta(Y)$ is $\eta$-club. But then for $\alp\in d_\eta(Y)\cap X$ we have $\{\bt<\alp:C\cap\bt=B^\bt_\alp\}=\{\bt<\alp:C\cap\bt=B^\bt\}=Y\cap\alp$ and  must therefore be $\eta$-stationary in $\alp$. But this contradicts our assumption that $\la B_\bt^\alp:\bt<\alp\ra$ was a witness that $\alp$ was not $\eta_\alp$-ineffable, as $C\cap\alp$ correctly guesses $B_\bt$ on a $\eta=\eta_\alp$-stationary subset of $\alp$. So we have a contradiction and $E_\gm$ is $\gm$-ineffable.
\qed

We can now show that $\gm$-ineffability is downward absolute to $L$. Recall that (\cite{De} VII.2.2.5) if $\kp$ is ineffable then $(\kp$ is ineffable$)^L$. To generalise this proof we shall use Theorem \ref{da}.
\begin{thm}\label{ineffda}
If $\kp$ is $\gm$-ineffable then $(\kp$ is $\gm$-ineffable$)^L$. 
\end{thm}
{\sc Proof:}
Let $\kp$ be $\gm$-ineffable. First we show that $\kp$ satisfies $(1)$ and $(2)$ of Theorem \ref{da}. By Lemma \ref{ineffclubnormal} we have that $\kp$ is $\gm$-normal and hence $\eta$-normal for any $\eta<\gm$.  By Lemma \ref{ineffbelow} we have $\{\alp<\kp:\text{for all $\eta<\gm$,~$\alp$ is $\eta$-ineffable}\}$ is $\gm$-ineffable. But if $\alp$ is $\eta$-ineffable then $\alp$ is $\eta$-normal so setting $E=\{\alp<\kp:\text{for all $\eta<\gm$,~$\alp$ is $\eta$-normal}\}$ then $E$ is $\gm$-ineffable in $\kp$ and hence $\gm+2$-stationary. Thus $$A_\gm=\{\alp<\kp:\text{\emph{for all $\eta$ with $1<\eta+1<\gm$~$\alp$ is $\eta$-normal}}\}\sps E$$ and so $A_\gm$ is $\gm$-stationary.

Now, using the characterisation of $\gm$-ineffability from Theorem \ref{ineff}, let $\la A_\alp:\alp<\kp\ra$ be a sequence in $L$ with each $A_\alp\sbs\alp$. This is clearly such a sequence in $V$, so as $E$ is $\gm$-ineffable we can find a set $A\sbs\kp$ such that $$X=\{\alp<\kp:A_\alp=A\cap\alp\text{ and for every $\eta<\gm$, $\alp$ is $\eta$-normal}\}$$ is $\gm$-stationary. Then by the weak compactness of $\kp$, as each $A_\alp\in L$ and $X$ is unbounded in $\kp$, we have $A\in L$. Setting $X'=\{\alp<\kp:A_\alp=A\cap\alp\}$, we have $X'\in L$ and $X\cap X'=X$ is $\gm$-stationary.  Thus by Theorem \ref{da} $X'$ is $\gm$-stationary in $L$, and hence $(\kp$ is $\gm$-ineffable$)^L$.
\qed


\subsection{Diamond Principles}

We now turn to generalising diamond ($\diamondsuit$) principles, and relate them to our generalised ineffability. 
Like $\Bon(\kp)$, $\diam_\kappa$ asserts the existence of a sequence of sets $S_\alp\sbs\alp$ for $\alp<\kp$ - in the case of $\diam$ the sequence must ``guess" any subset of $\kp$ sufficiently often. For the original $\diam$, ``sufficiently'' is ``stationarily'', so by altering this to $\gm$-stationarity we can define a new notion: 

\begin{defn}
$\diamondsuit^\gm_\kappa$ is the assertion that there is a sequence $\la S_\alp:\alp<\kp\ra$ such that for any $S\subseteq\kp$ we have $\{\alp<\kp:S_\alp=S\cap\alp\}$ is $\gm$-stationary in $\kp$.
\end{defn}

The original principle is thus $\diamondsuit^1_\kappa$.

\begin{rem} As with ineffability, these principles get stronger as $\gm$ increases: for $\bt<\gm$ we have $\diamondsuit^\gm_\kappa\Rightarrow\diamondsuit^\bt_\kappa$.
\end{rem}

\begin{thm}
If $\kp$ is $\gm$-ineffable then $\diamondsuit^\gm_\kappa$ holds.
\end{thm}

{\sc Proof:}
Define by recursion a sequence $\la(S_\alp,C_\alp, \eta_\alp):\alp<\kp\ra$: let $(S_\alp,C_\alp)$ be any pair of subsets of $\alp$ with $\eta_\alp<\gm$ such that $C_\alp$ is $\eta_\alp$-club in $\alp$ and for any $\bt\in C_\alp$ we have $S_\alp\cap \bt\neq S_\bt$. If there is no such pair, set $S_\alp=C_\alp=\eta_\alp=\emptyset$. We now use $\gm$-inneffability to see that $\la S_\alp:\alp<\kp\ra$ is a $\diamondsuit^\gm_\kappa$-sequence.

By the characterisation of $\gm$-ineffability in Theorem \ref{ineff} and some simple coding (noting that successor levels are irrelevant), we can find $S,C\subseteq\kp$ and $\eta<\gm$ such that $$A=\{\alp<\kp:S\cap\alp=S_\alp\wedge C\cap\alp=C_\alp\wedge\eta_\alp=\eta \}$$ is $\gm$-stationary. 

Let $X\subseteq\kp$ and suppose that $B:=\{\alp<\kp:S_\alp=X\cap\alp\}$ is not $\gm$-stationary in $\kp$. Take $D$ $\gm'$-club in $\kp$ witnessing this, i.e $\gm'<\gm$ and for any $\alp\in D$, $X\cap\alp\neq S_\alp$. By $\gm$-stationarity of $A$ we can pick $\alp<\bt$ both in $A\cap d_{\gm'}(D)$. Then $C_\bt\cap\alp=C\cap\bt\cap\alp=C\cap\alp=C_\alp$. Now as $(X\cap \alp, D\cap\alp,\gm')$ works for the definition of $(S_\alp,C_\alp,\eta_\alp)$ we must have that $C_\alp$ is $\eta$-club in $\alp$. Hence $C_\bt\cap\alp=C_\alp$ is $\eta$-stationary in $\alp$, and so $C_\bt\neq\emp$ and we have by $\eta$-stationary closure that $\alp\in C_\bt$. Then by definition of $(S_\bt,C_\bt,\eta_\bt)$ we must have $S_\alp\neq S_\bt\cap\alp$. But this is a contradiction as $\alp,\bt\in A$ so $S_\alp=S\cap\alp=S_\bt\cap\alp$.
\qed

The $\diam$ principle itself is not a large cardinal notion, and in fact within the constructible universe $\diam_\kappa$ holds at every regular uncountable cardinal $\kp$ (see \cite{De} Chapter IV). Our new $\diamondsuit^\gm_\kappa$ principles also hold in $L$, with minimal assumptions on $\kp$, as we shall now see. 
Recall the following:
\begin{defn}
A model $\la M,\in\ra$ is $\gm$-stationary correct at $\kp$ if for any $S\in\P(\kp)\cap M$, $M\vDash``S$ is $\gm$-stationary in $\kp"$ iff $S$ is $\gm$-stationary in $\kp$.
\end{defn}

\begin{defn}
The \emph{$\gm$-trace} of $M,p$ on $\alp$ is denoted $\mathcal{s}^{\gm}(M,p,\alp)$ and consists of all $\bt<\alp$ such that $\bt\in\mathcal{s}(M,p,\alp)$ and if $\pi:M\{p\cup\bt\cup\{\alp\}\}\cong N$ is the transitive collapse then $N$ is $\gm$-stationary correct for $\bt=\pi(\alp)$.
\end{defn}

As $\gm$ stationarity is $\Pi^1_\gm$ expressible we have, as a corollary to Lemma \ref{traceclub3}, the following:
\begin{lem}
If $V=L$ and $\kp$ is $\Pi^1_\gm$-indescribable then for any limit $\nu>\kp$ with $L_\nu$ $\Pi^1_\gm$-correct over $\kp$ we have $\mathcal{s}^{\gm}(L_\nu,p,\kp)$ is $\gm$-club in $\kp$.
\end{lem}



\begin{prop}
If $V=L$ and $\kp$ is a $\gm$-stationary regular cardinal then $\diamondsuit^\gm_\kappa$ holds. 
\end{prop}
{\sc Proof:}
By recursion we define, for each $\alp<\kp$, an ordinal $\eta_\alp<\gm$ and a pair of subsets of $\alp$ $(S_\alp,C_\alp)$ (the $S_\alp$'s will form our $\diamondsuit^\gm$-sequence).

Assume we have defined $\la S_\bt:\bt<\alp\ra$. Let $\psi(\alp,\eta,C, S)$ be the statement that $\eta$ is an ordinal below $\gm$, and $(S,C)$ is a pair of subsets of $\alp$ with $C$ $\eta$-club in $\alp$ and for any $\bt\in C$, $S\cap\bt\neq S_\bt$.

If there are $\eta$, $S$ and $C$ such that $\psi(\alp,\eta,S,C)$ holds, we take $\eta_\alp$ to be the least such $\eta$, and then $(S_\alp,C_\alp)$ the $<_L$ least pair with $\psi(\alp,\eta_\alp,S_\alp,C_\alp)$. If no such $\eta$ exists set $S_\alp=C_\alp=\eta_\alp=\emp$.

Suppose $\la S_\alp:\alp<\kp\ra$ is not a $\diamondsuit^{\gm}_\kappa$-sequence and take $\eta$ and $(S,C,\eta)$ to be the least witness to this, i.e $\psi(\kp,\eta,S,C)$ holds and if $\eta'<\eta$ or $(S',C')<_L(S,C)$ then $\neg\psi(\kp,\eta',S',C')$. Note that all this can be carried out in $L_{\kp^+}$, and that as $\eta<\gm$ and we are in $L$, $\kp$ is $\Pi^1_\eta$-indescribable and hence by Lemma \ref{traceclub3} the $\eta$ trace forms an $\eta$-club.

Suppose $\alp\in\trace^{\eta}(L_{\kp^+},\kp,\{S,C,\eta\})$ with $\alp>\eta$ and $L_\nu\cong L_{\kp^+}\{\alp,\{S,C,\eta\}\}$.

$L_{\nu}\vDash ``\eta$ is the least ordinal such that there exists a pair $(C',S')$ with $\psi(\alp,\eta,S',C,)$ and $(C\cap\alp,S\cap\alp)$ is the $<_L$ least such pair"

As $\alp$ is in the $\eta$ trace, we have $L_\nu$ is $\eta$-stationary correct so $\psi(\alp,\eta,S\cap\alp,C\cap\alp)$ must hold and $\psi(\alp,\eta',S',C,)$ fail for $\eta'<\eta$ or $(S',C')<_L(S\cap\alp,C\cap\alp)$. But this was the definition of $S_\alp$ so $S_\alp=S\cap\alp$, contradiction.
\qed

A stronger principle than $\diamondsuit$ is $\diamondsuit^*$, which requires sets to be guessed on a club set of $\alp$, but allows for more guesses at each $\alp$. Unlike $\diamondsuit$, $\diamondsuit^*$ is incompatible with ineffability. The original principle is $\diamondsuit^{*1}$ in the following definition.

\begin{defn}
$\diamondsuit^{*\gm}_\kappa$ is the assertion that there is a sequence $\la A_\alp:\alp<\kp\ra$ such that  $A_\alp\subseteq \P(\alp)$ and $|A_\alp|\leq|\alp|$ for each $\alp<\kp$, and for any $X\subseteq\kp$ there is some $\gm'<\gm$ such that $\{\alp<\kp:X\cap\alp\in A_\alp\}$ is in the $\gm'$-club filter on $\kp$.
\end{defn}

 \begin{rem} In contrast to the case for $\diamondsuit$, here we have that at $\gm$-reflecting cardinals $\kp$,  $\diamondsuit^{*{\gm'}}_\kappa$ implies $\diamondsuit^{*\gm}_\kappa$ for $\gm'<\gm$. 
 \end{rem}
 
As for the original principle, $\diamondsuit^{*\gm}$ cannot hold at a $\gm$-ineffable cardinal:
 
 \begin{thm}\label{ineffdiamondstar}
If $\kp$ is $\gm$-ineffable then $\diamondsuit^{*\gm}_\kappa$ fails.
\end{thm} 

{\sc Proof:}
Suppose $\kp$ is $\gm$-ineffable and let $\la A_\alp:\alp<\kp\ra$ be a sequence such that  $A_\alp\subseteq \P(\alp)$ and $|A_\alp|\leq|\alp|$ for each $\alp<\kp$. We find $B\sbs\kp$ such that $\{\alp<\kp:B\cap\alp\notin A_\alp\}$ is $\gm$-stationary. 
For each $\alp<\kp$ let $B_\alp\sbs\alp$ be a set different from each set in $A_\alp$  - we can find such a set as $|A_\alp|=\alp$. Now by $\gm$-ineffability there is $B\sbs\kp$ such that $X=\{\alp<\kp:B_\alp=B\cap\alp\}$ is $\gm$-stationary. But then $B$ is not guessed by $\la A_\alp:\alp<\kp\ra$ on $X$, so $\la A_\alp:\alp<\kp\ra$ cannot be a $\diamondsuit^{*\gm}$-sequence. As $\la A_\alp:\alp<\kp\ra$ was arbitrary, $\diamondsuit^{*\gm}_\kappa$ fails.\\
\qed

In fact, the notion of $\diamondsuit^{\gm}$ is only really of interest for a successor ordinal $\gm$: 

\begin{prop}\label{boring limits}
If $\gm$ is a limit ordinal and $\diamondsuit^{\gm}_\kappa$ holds iff there is some $\gm'<\gm$ such that $\diamondsuit^{\gm'}_\kappa$.
\end{prop}

{\sc Proof:}
Suppose $\gm$ is a limit ordinal and  $\diamondsuit^{\gm}_\kappa$ holds. 
Let $\A=\la A_\alp:\alp<\kp\ra$ a sequence such that  $A_\alp\subseteq \P(\alp)$ and $|A_\alp|\leq|\alp|$ for each $\alp<\kp$. Assume that $\A$ is constructibly closed, in that for each $\alp$, $A_\alp=\P(\alp)\cap L_\nu[A_\alp]$ for some limit ordinal $\nu$ with $\alp<\nu<\alp^+$ - clearly we can always expand $\A$ to satisfy this condition, and $\A$ will remain a  $\diamondsuit^{*\gm}$-sequence. Suppose for each $\gm'<\gm$, $\A$ is not a $\diam^{*\gm'}$-sequence, and take $B_{\gm'}$ to witness this, setting $X_{\gm'}=\{\alp<\kp:B_{\gm'}\cap\alp\notin A_\alp\}$, which is $\gm'$-stationary. Set $X=\bigcup X_{\gm'}$ and let $B$ code (definably and uniformly) each of the $B_{\gm'}$'s. Then $X$ is $\gm$-stationary and we claim $X\sbs\{\alp<\kp:B\cap\alp\notin A_\alp\}$. Suppose $\alp\in X$ and $B\cap\alp\in A_\alp$. Then as the coding was definable, each $B_{\gm'}\cap\alp\in A_\alp$ - contradiction.
\qed

\begin{col}
If $\gm$ is a limit ordinal and $\kp$ is $\gm'$-ineffable for every $\gm'<\gm$ then $\diam^{*\gm}$ fails.
\end{col}

In $L$ we have that, given the precondition of $\gm$-stationarity, the failure of $\diamondsuit^{*\gm}$ characterises the $\gm$-ineffables - but only for successor ordinals $\gm$. 

\begin{thm}\label{ineffdiamL}
Assume $V=L$ and let $\kp$ be a regular cardinal that is $\gm$-stationary with $\gm$ a successor ordinal. Then $\kp$ is not $\gm$-ineffable iff $\diamondsuit^{*\gm}_\kappa$ holds.
\end{thm}

{\sc Proof:}
($\Rightarrow$) Let $\kp$ be a regular uncountable cardinal which is not $\gm$-ineffable. Let $\la A_\alp:\alp<\kp\ra$ be the $<_L$ least sequence such that $A_\alp\subseteq \alp$ for each $\alp$ and for any $A\subseteq\kp$ we have $\{\alp<\kp:A_\alp=A\cap\alp\}$ is not $\gm$-stationary. 
For each $\alp<\kp$ set $M_\alp$ to be the least $M\prec L_\kp$ such that $\alp+1\subseteq M_\alp$ and $\la A_\alp:\alp<\kp\ra\in M$. Let $\sigma_\alp:M_\alp\cong L_{\nu_\alp}$.

Set $S_\alp=\P(\alp)\cap L_{\nu_\alp}$, and note $|L_{\nu_\alp}|=|\nu_\alp|=|\alp|$, so $|S_\alp|\leq|\alp|$. we shall show that if $\eta+1=\gm$ then for any $X\subseteq\kp$, setting $C=\trace^{\eta}(L_{\kp^+},\{X\},\kp)$ we have $C$ is $\eta$-club and $X\cap\alp\in S_\alp$ for all $\alp\in C$.
 Fix such an $X$ 
 and take $\alp\in\trace^{\eta}(L_{\kp^+},\{X\},\kp)$.
 Let $\pi:N_\alp\cong L_{\mu}$. Then $\pi\restr\alp=\alp$, $\pi(\kp)=\alp$ and $\pi(X)=X\cap\alp$, so $X\cap\alp\in L_\bt$. Thus we are done if we can show $\mu\leq\nu_\alp$.
 Suppose to the contrary that $\mu>\nu_\alp$. 
 As $\la A_\alp:\alp<\kp\ra$ is definable in $L_{\kp^+}$ and $\alp\in M_\alp$ we have $\la A_\gm:\gm\leq\alp\ra\in M_\alp$. Then $\sigma_\alp(\la A_\gm:\gm\leq\alp\ra)=\la A_\gm:\gm\leq\alp\ra$ because $\alp+1\subseteq M_\alp$, so $\la A_\gm:\gm\leq\alp\ra\in L_{\nu_\alp}\subseteq L_\bt$.
 Setting $E=\{\gm<\alp:A_\gm=A_\alp\cap\gm\}$ we have $E\in L_\bt$ and
 $\pi^{-1}(E)=\{\gm<\kp:A_\gm=\gm\cap\pi^{-1}(A_\alp)\}$.
First suppose $L_\bt\vDash ``E$ is $\gm$-stationary in $\alp$". Then $\pi^{-1}(E)$ is $\gm$-stationary in $\kp$ by elementarity. Also we have (in $L_{\kp^+}$ and hence in $L$) that $\pi^{-1}(A_\alp)\subseteq\kp$ and $\pi^{-1}(E)=\{\bt<\kp:A_\bt=\bt\cap\pi^{-1}(A_\alp)\}$. But this is a contradiction as  $\la A_\alp:\alp<\kp\ra$ was chosen to witness $\kp$ not being $\gm$-ineffable.
 
 So we must have $L_\bt\vDash ``E$ is not $\gm$-stationary in $\alp$". Thus for some $C\subseteq\alp$, $L_\bt\vDash ``C$ is $\eta$-club in $\alp$ and $C\cap E=\emp$". Then inverting the collapse we get $C'=\pi^{-1}(C)$ is $\eta$-club in $\kp$ and $C'\cap\pi^{-1}( E)=\emp$". As $L_\bt$ is in the $\eta$-trace it is $\eta$-stationary correct for $\alp$. Now $C'\cap \alp=C$ and by $\eta$-stationary correctness, we have that $C'\cap\alp$ is $\eta$-stationary, and hence $\alp\in C'$. Thus by the definition of $E$ we have $A_\alp\neq\alp\cap \pi^{-1}(A_\alp)$. But this is a contradiction as $\pi\restr\alp=id\restr\alp$. Thus we must have $\mu\geq\nu_\alp$ and we're done.

 The converse direction ($\Leftarrow$) follows from Theorem \ref{ineffdiamondstar}.
 \qed

In fact this characterisation can fail for limit $\gm$: Assuming there is an $\om$-ineffable we can take $\kp$ to be the least cardinal which is $n$-ineffable for every $n<\om$. Then $\kp$ is not $\Pi^1_\om$-indescribable as being $n$-ineffable is $\Pi^1_{n+2}$ over $L_\kp$, so being $n$-ineffable for every $n<\om$ is $\Pi^1_\om$. Hence (as we are in $L$), $\kp$ is not $\om$-reflecting, and so $\kp$ is not $\om$-ineffable. Now, for each $n<\om$, as $\kp$ is $n$-ineffable $\diam^{*n}$ fails, so by Proposition \ref{boring limits} $\diam^{*\om}$ also fails.
We can however give a complete description of when $\diam^{*\gm}_\kappa$ holds in $L$:

\begin{col}\label{ineffdiamcol}
If $V=L$ then $\diam^{*\gm}_\kappa$ holds iff $\kp$ is $\gm$-stationary and one of the following:
\begin{enumerate}
\item $\gm$ is a successor and $\kp$ is not $\gm$-ineffable
\item $\gm$ is a limit and $\kp$ is not $\gm'$-ineffable for some $\gm'<\gm$.
\end{enumerate}
\end{col}

{\sc Proof:}
For $\gm$ a successor this is just Theorem \ref{ineffdiamL}. For $\gm$ a limit by Proposition \ref{boring limits} if $\diam^{*\gm}_\kappa$ holds then for some $\gm'<\gm$, $\diam^{*\gm'}_\kappa$ holds so by Theorem \ref{ineffdiamondstar} $\kp$ is not $\gm'$-ineffable. The converse direction is immediate from Theorem \ref{ineffdiamL}, remembering that for $\gm'<\gm$, $\diam^{*\gm'}_\kappa$ implies $\diam^{*\gm}_\kappa$.
\qed

Now we turn to the relationship between $\diam$ and $\diam^*$. In the classical case $\diam^*_\kappa$ implies $\diam_\kappa$ (Kunen. For a proof see, for instance, \cite{Sch14}, 5.38).  If the $\gm$-club filter on $\kp$ is normal then we can generalise this: $\diamondsuit^{\gm+1}_\kappa$ implies $\diamondsuit^{\gm+1}_\kappa$. The proof goes \via \!a weaker principle, $\diamondsuit^{-\gm}$:

\begin{defn}
$\diamondsuit^{-\gm}_\kappa$ is the assertion that there is a sequence $\la A_\alp:\alp<\kp\ra$ such that  $A_\alp\subseteq \P(\alp)$, $|A_\alp|\leq|\alp|$ for each $\alp<\kp$, and for any $X\subseteq\kp$ the set $\{\alp<\kp:X\cap\alp\in A_\alp\}$ is $\gm$-stationary in $\kp$.
\end{defn}
This is a clear weakening of both $\diam^*$ and $\diam$.

\begin{thm}
If the $\gm$-club filter on $\kp$ is normal then $\diam^{\gm+1}_\kappa\leftrightarrow\diam^{-\gm+1}_\kappa$, and hence $\diam^{*\gm+1}_\kappa\rightarrow\diam^{\gm+1}_\kappa$.
\end{thm}

{\sc Proof:}
It is clear that $\diamondsuit^{\gm}\rightarrow\diam^{-\gm}$ and $\diam^{\gm}\rightarrow\diam^{-\gm}\diam^{-}$. Suppose the $\gm$-club filter on $\kp$ is normal and $\diam^{-\gm+1}_\kappa$ holds. Let $\la \A_\alp:\alp<\kp\ra$ be a $\diam^{-\gm+1}$-sequence such that each $A\in\A_\alp$ codes a subset of $\alp\times\alp$. Let $\{A_\alp^\bt:\bt<\alp\}$ enumerate these coded sets in $\A_\alp$. 



For a given $\alp$ and $\bt<\alp$, set $S^\bt_\alp=\{\nu\in\alp:(\nu,\bt)\in A^\bt_\alp\}$. we shall see that for some $\bt<\kp$, $\la S_\alp^\bt:\alp<\kp\ra$ is a $\diam_{\gm+1}$-sequence. Suppose this is not the case. For each $\bt<\kp$ take $X^\bt\sbs\kp$ and $C^\bt$ $\gm$-club such that for $\alp\in C^\bt$, $X^\bt\cap\alp\neq S^\bt_\alp$. Set $$X=\bigcup_{\bt<\kp}X^\bt\times\{\bt\}$$ and $$C=\triangle_{\bt<\kp}C^\bt.$$
By the normality of the  $\gm$-club filter, $C$ is $\gm$-club. Then for $\alp\in C$, if $\bt<\alp$ then $X^\bt\neq S_\alp^\bt$ and so $X\cap(\alp\times\alp)\neq A^\bt_\alp$. Thus $X\notin \A_\alp$. But this contradicts $\la \A_\alp:\alp<\kp\ra$ being a $\diam^{-\gm+1}$-sequence.
\qed

\bibliography{settheory10r}
\bibliographystyle{asl}




\end{document}